\documentclass[10pt,a4paper]{article}
\usepackage{amsmath,amssymb,amsthm,graphics,graphicx,epic,eepic,multicol,ascmac}

\usepackage[all]{xy}

\numberwithin{equation}{section}
\theoremstyle{theorem}
\newtheorem{thm}{Theorem}[section]
\newtheorem{prop}[thm]{Proposition}
\newtheorem{lem}[thm]{Lemma}

\newtheorem{ex}[thm]{Example}

\theoremstyle{definition}
\newtheorem{defn}[thm]{Definition}

\def\Lm{\Lambda}
\def\nd{\noindent}
\def\ovl{\overline}

\begin{document}

\title{An algorithm for Berenstein-Kazhdan decoration functions and
trails for classical Lie algebras}

\author{Yuki Kanakubo\thanks{Faculty of Pure and Applied Sciences, University of Tsukuba,
1-1-1 Tennodai, Tsukuba, Ibaraki 305-8577,
Japan: {y-kanakubo@math.tsukuba.ac.jp}.}, Gleb Koshevoy\thanks{Institute for Information Transmission Problems Russian Academy of Sciences,
Russian Federation :
{koshevoy@cemi.rssi.ru}.} and Toshiki Nakashima\thanks{Division of Mathematics, 
Sophia University, Kioicho 7-1, Chiyoda-ku, Tokyo 102-8554,
Japan: {toshiki@sophia.ac.jp}.} } 
\date{}
\maketitle
\begin{abstract}
For a simply connected connected simple algebraic group $G$,
it is known that a variety $B_{w_0}^-:=B^-\cap U\overline{w_0}U$
has a geometric crystal structure with a positive structure
$\theta^-_{\mathbf{i}}:(\mathbb{C}^{\times})^{l(w_0)}\rightarrow B_{w_0}^-$
for each reduced word $\mathbf{i}$ of the longest element $w_0$ of Weyl group.
A rational function $\Phi^h_{BK}=\sum_{i\in I}\Delta_{w_0\Lambda_i,s_i\Lambda_i}$
on $B_{w_0}^-$ is called a half-potential, where $\Delta_{w_0\Lambda_i,s_i\Lambda_i}$ is a generalized minor.
Computing $\Phi^h_{BK}\circ \theta^-_{\mathbf{i}}$ explicitly,
we get an explicit form of string cone or polyhedral realization of $B(\infty)$
for the finite dimensional simple Lie algebra $\mathfrak{g}={\rm Lie}(G)$.

In this paper, for
an arbitrary reduced word $\mathbf{i}$,
we give an algorithm to compute the summand
$\Delta_{w_0\Lambda_i,s_i\Lambda_i}\circ \theta^-_{\mathbf{i}}$ of $\Phi^h_{BK}\circ \theta^-_{\mathbf{i}}$
in the case 
$i\in I$ satisfies that
for any weight $\mu$ of $V(-w_0\Lambda_i)$ and $t\in I$, it holds
$\langle
h_t,\mu
\rangle\in\{2,1,0,-1,-2\}$. In particular, if $\mathfrak{g}$ is of type ${\rm A}_n$, ${\rm B}_n$, ${\rm C}_n$ or
${\rm D}_n$ then all $i\in I$ satisfy this condition so that
one can completely calculate $\Phi^h_{BK}\circ \theta^-_{\mathbf{i}}$.
We will also prove that our algorithm works in the case 
$\mathfrak{g}$ is of type ${\rm G}_2$.

\end{abstract}

\vspace{-10pt}

\section{Introduction}

In \cite{BK0},
the notion of `geometric crystals' were
defined as irreducible algebraic varieties $X$ equipped with certain
$\mathbb{C}^{\times}$-actions and rational functions,
which are geometric analogs of Kashiwara's crystals
with
Kashiwara operators and $\varepsilon$-functions, weight functions.
If there exists a birational map $\theta:T'\rightarrow X$ called a positive structure
with an algebraic torus $T'$ then one can obtain a crystal via
a tropicalization functor ${\rm Trop}$ (see subsections \ref{GC-trop},\ \ref{Trop-geom}) from the geometric crystal $X$.
By this functor, the torus $T'$ corresponds to the set of cocharacters ${\rm Trop}(T')=X_*(T')={\rm Hom}(\mathbb{C}^{\times},T')$.

Several varieties related to a reductive group $G$ have geometric crystal structures.
As shown in \cite{BK0, BK}, the varieties
$B_{w_0}^-=B^-\cap U\ovl{w_0}U$ and $T\cdot B_{w_0}^-$ have geometric crystal structures,
where $B$, $B^-$ are opposite Borel subgroups, $U\subset B$ is a unipotent radical, $T=B\cap B^-$ is a maximal torus
and $\ovl{w_0}\in {\rm Norm}_G(T)$ is a representative of the longest element $w_0$ in Weyl group $W={\rm Norm}_G(T)/T$. 
The variety $T\cdot B_{w_0}^-$ has a 
positive structure $\theta_{\mathbf{i}}: T\times (\mathbb{C}^{\times})^{l(w_0)}\rightarrow T\cdot B_{w_0}^-$
associated with each reduced word $\mathbf{i}$ of $w_0$ so that
we obtain a crystal $X_*(T\times (\mathbb{C}^{\times})^{l(w_0)})$ by the tropicalization functor.
The {\it Berenstein-Kazhdan decoration function}  $\Phi_{BK}$ on 
$T\cdot B_{w_0}^-$
  is defined as
\begin{equation}\label{BK-def}
\Phi_{BK}=\sum_{i\in I}\frac{\Delta_{w_0\Lambda_i,s_i\Lambda_i}}{\Delta_{w_0\Lambda_i,\Lambda_i}}
+
\sum_{i\in I}\frac{\Delta_{w_0s_i\Lambda_i,\Lambda_i}}{\Delta_{w_0\Lambda_i,\Lambda_i}}.
\end{equation}
Here, $\Lambda_i$ is the $i$-th fundamental weight, for $u$, $v\in W$, the function $\Delta_{u\Lambda_i,v\Lambda_i}$
is a generalized minor (Definition \ref{gen-def}).
Considering the tropicalization of the rational function $\Phi_{BK}$ on $T\cdot B_{w_0}^-$,
one obtains a subcrystal 
\[
\{ z\in X_*(T\times (\mathbb{C}^{\times})^{l(w_0)}) | {\rm Trop}(\Phi_{BK}\circ \theta_{\mathbf{i}})(z)\geq0 \},
\]
which is isomorphic to the disjoint union of all crystal bases $B(\lambda)$ of
the finite dimensional irreducible representations of the quantum group $U_q(^L\mathfrak{g})$ with highest weights $\lambda$ \cite{BK}. 
Here, $^L\mathfrak{g}$ is the Langlands dual Lie algebra of $\mathfrak{g}={\rm Lie}(G)$.
In \cite{GKS20}, relations between $\Phi_{BK}$
and Gross-Hacking-Keel-Kontsevich potentials (\cite{GHKK}) 
on simply connected connected simply laced simple algebraic group
are established. 

In \cite{KN}, a {\it half-potential} $\Phi^h_{BK}=\sum_{i\in I}\Delta_{w_0\Lambda_i,s_i\Lambda_i}$ is introduced,
which is the restriction of the first term of $\Phi_{BK}$ in (\ref{BK-def}) to $B_{w_0}^-$.
Defining a positive structure $\theta^-_{\mathbf{i}}: (\mathbb{C}^{\times})^{N}\rightarrow B_{w_0}^-$ on $B_{w_0}^-$ for each reduced word $\mathbf{i}=(i_1,\cdots,i_N)$ of $w_0$,
we get a subcrystal
\begin{equation}\label{cone1}
\{ z\in X_*((\mathbb{C}^{\times})^{N}) | {\rm Trop}(\Phi^h_{BK}\circ \theta_{\mathbf{i}}^-)(z)\geq0 \},
\end{equation}
which is isomorphic to the crystal base $B(\infty)$ of the negative part $U_q^-(^L\mathfrak{g})\subset U_q(^L\mathfrak{g})$.
Under the identification $X_*((\mathbb{C}^{\times})^{N})=\mathbb{Z}^N$,
the set (\ref{cone1}) coincides with the set of integer points in a string cone (\cite{Lit}) and
in \cite{KN}, it is shown that the crystal (\ref{cone1}) is isomorphic to a polyhedral realization (\cite{NZ}) of $B(\infty)$. 
It is a natural problem to get explicit forms of string cone (or, equivalently, the set (\ref{cone1}) or polyhedral realization)
 and
many researchers have worked on this problem.
In \cite{Lit}, Littelmann gave explicit forms of string cones
for finite dimensional simple Lie algebras and specific reduced words $\mathbf{i}$
called nice decompositions.
In the case $\mathfrak{g}$ is of type ${\rm A}_n$, Gleizer and Postnikov
gave a purely combinatorial rule to compute the inequalities defining the string cone
using the rigorous paths in a graph constructed from wiring diagrams \cite{GP}.
In \cite{GKS16}, a combinatorial expression of string cones via dual Reineke vectors constructed by rhombus tiling tools is given. 
The polyhedral realizations
associated with
a sequence $\iota$ in $I$ are defined for any symmetrizable Kac-Moody Lie algebras.
In \cite{H1,H2,Ka,KaN,KS}, 
explicit forms of polyhedral realizations for finite dimensional simple Lie algebras
or several classical affine Lie algebras are given for sequences $\iota$ satisfying a certain condition.

In our previous paper \cite{KKN}, for $i\in I$ such that the finite dimensional irreducible $\mathfrak{g}$-module $V(\Lambda_i)$ with highest weight $\Lambda_i$
is a minuscule representation, we invented an algorithm to compute all monomials in $\Delta_{w_0\Lambda_i,s_i\Lambda_i}\circ \theta_{\mathbf{i}}^-(t_1,\cdots,t_N)$ explicitly for all reduced words $\mathbf{i}$. 
The algorithm generates an edge-colored directed graph $DG$ with vertices labelled by the monomials in $\Delta_{w_0\Lambda_i,s_i\Lambda_i}\circ \theta_{\mathbf{i}}^-(t_1,\cdots,t_N)$.
In particular, if $\mathfrak{g}$ is of type ${\rm A}_n$ then $V(\Lambda_i)$ is minuscule for all $i\in I$.
Thus,
one obtains an explicit form (\ref{cone1}) of $B(\infty)$
in $X_*((\mathbb{C}^{\times})^{N})=\mathbb{Z}^{N}$ via tropicalization.

The main result of this paper (Theorem \ref{thm2}) allows us, for all reduced words $\mathbf{i}$, to 
get all monomials  
in $\Delta_{w_0\Lm_i,s_i\Lm_i}\circ \theta^-_{\mathbf{i}}(t_1,\cdots,t_N)$ explicitly 
 in the following cases
(the numbering of Dynkin diagram is same as in \cite{Kac}), which covers a significantly wide range of indices $i\in I$
comparing with \cite{KKN}.
\begin{table}[h]
  \begin{tabular}{|c|c|c|c|c|c|c|c|c|c|} \hline
  $\mathfrak{g}$ & ${\rm A}_n$ & ${\rm B}_n$ & ${\rm C}_n$ & ${\rm D}_n$ & ${\rm E}_6$ & ${\rm E}_7$ & ${\rm E}_8$ & ${\rm F}_4$ & ${\rm G}_2$\\ \hline
  $i$ & \text{all }$i\in I$ & \text{all }$i\in I$ & \text{all }$i\in I$ & \text{all }$i\in I$ & $1,2,4$, $5,6$ & $1,5,6,7$  & $1,7$ & $1,4$ & \text{all }$i\in I$\\ \hline
  \end{tabular}
\end{table} 
Due to Theorem \ref{thm2}, we compute an edge-colored directed graph $\overline{DG}$ whose vertices are labelled by the monomials in 
$\Delta_{w_0\Lambda_i,s_i\Lambda_i}\circ \theta_{\mathbf{i}}^-(t_1,\cdots,t_{N})$,
and edges are colored by letters of $\{1,2,\cdots,N\}$.
We only use easy computations of the Weyl group action on simple roots and weights and 
multiplications of Laurent monomials.
In particular, in case of $\mathfrak{g}$ is of classical type (${\rm A}_n$, ${\rm B}_n$, ${\rm C}_n$ or ${\rm D}_n$) or type ${\rm G}_2$,
by the tropicalization, we get an explicit form of the crystal (\ref{cone1}) for any reduced word $\mathbf i$.
For example, in the case $\mathfrak{g}$ is of type ${\rm C}_3$, $\mathbf{i}=(2,3,2,1,2,3,2,3,1)$ and $i=2$,
we get the following graph $\overline{DG}$:
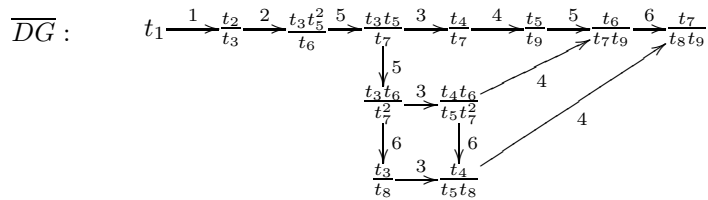
\[
\begin{xy}
(-35,20) *{\overline{DG} :}="0",
(-20,20) *{t_1}="1",
(-10,20) *{\frac{t_2}{t_3}}="2",
(0,20) *{\frac{t_3t_5^2}{t_6}}="3",
(10,20) *{\frac{t_3t_5}{t_7}}="4",
(10,10) *{\frac{t_3t_6}{t_7^2}}="4-1",
(10,0) *{\frac{t_3}{t_8}}="4-2",
(20,10) *{\frac{t_4t_6}{t_5t_7^2}}="5-1",
(20,0) *{\frac{t_4}{t_5t_8}}="5-2",
(20,20) *{\frac{t_4}{t_7}}="5",
(30,20) *{\frac{t_5}{t_9}}="6",
(40,20) *{\frac{t_6}{t_7t_9}}="7",
(50,20) *{\frac{t_7}{t_8t_9}}="8",
\ar@{->} "1";"2"^{1}
\ar@{->} "2";"3"^{2\ }
\ar@{->} "3";"4"^{5\ }
\ar@{->} "4";"5"^3
\ar@{->} "4";"4-1"^{5}
\ar@{->} "5";"6"^{4}
\ar@{->} "4-1";"5-1"^{3}
\ar@{->} "4-1";"4-2"^{6}
\ar@{->} "5-1";"5-2"^{6}
\ar@{->} "5-1";"7"_{4}
\ar@{->} "5-2";"8"_{4}
\ar@{->} "4-2";"5-2"^{3}
\ar@{->} "6";"7"^{5}
\ar@{->} "7";"8"^{6}
\end{xy}
\]
Theorem \ref{thm2} implies
$\Delta_{w_0\Lm_2,s_2\Lm_2}\circ \theta^-_{\mathbf{i}}(t_1,\cdots,t_9)$ is a
linear combination of monomials in $\overline{DG}$ with positive integer coefficients.
Similarly, we get graphs $\overline{DG}$ for $i=1,3$ and can find terms in
$\Delta_{w_0\Lm_i,s_i\Lm_i}\circ \theta^-_{\mathbf{i}}(t_1,\cdots,t_9)$.
Hence, one obtains an
explicit form of the crystal (\ref{cone1}) (see Example \ref{ex-1} for detail). 
The rule to construct the new monomials and connecting edges is similar to the construction of polyhedral realizations given in \cite{NZ}
or action of Kashiwara operators of {\em monomial realizations}
of Kashiwara's crystals \cite{K,Nj}.
We get new monomials and connecting edges one after another
without any complicated combinatorial or representation theoretical tools. 
We also show that in the case $V(\Lambda_i)$ is minuscule, the graphs $DG$ and $\overline{DG}$ coincide. 
The proof is done by computing monomials of generalized minors
by a representation theoretical method of Berenstein-Zelevinsky $\mathbf i$-trails \cite{BZ}.

The paper is organized as follows. Section 2 recalls definitions and important facts on $\mathbf i$-trails.  Section 3 recalls basics on geometric crystals and the Berenstein-Kazhdan potential. Subsection 3.5 introduces some of properties of the half-potentials. Namely, Proposition 3.11 and Lemma 3.13 show that there is unique monomial with non-negative exponents in the half-potential,
which is the unique source of our graph $\overline{DG}$. The monomial of Proposition 3.10 is one of the sinks. 
Section 4 introduces main theorem and two examples of our algorithm for types $\rm C_3$ and $\rm D_4$.
The proof of the main theorem is obtained in Section 5 ($\rm G_2$-type case is not included here). In 
Lemma \ref{lem-a2}, \ref{lem-a3}, we consider the condition (a) of (2) in Theorem \ref{thm2}.
In Lemma \ref{lem-a4}, \ref{lem-a5}, we consider the condition (b) of (2) in Theorem \ref{thm2}.
Section 6 shows that the algorithm from Theorem \ref{thm2} can be applied for $\rm G_2$-type case,
and in the case $V(\Lambda_i)$ is minuscule, the graph $\overline{DG}$ coincides with the graph $DG$ of \cite{KKN}.

\vspace{2mm}

\nd
{\bf Acknowledgements}
Y.K. is supported by JSPS KAKENHI Grant Number JP20J00186.
G.K. thanks to the grant RSF 21-11-00283 for support.
T.N. is supported in part by JSPS KAKENHI Grant Number JP20K03564. We also thank Denis Mironov for writing  computer programs and friendly assistance.

\section{$\mathbf{i}$-trails and generalized minors}

\subsection{Notation}\label{notation}

Let $G$ be a simply connected connected simple algebraic group,
$B,\ B^-\subset G$ its Borel subgroups, $T:=B\cap B^-$ the maximal torus,
$W={\rm Norm}_G(T)/T$ Weyl group, 
$U$, $U^-$ be unipotent radicals
of $B$, $B^-$,
$A=(a_{i,j})$ the Cartan matrix
of $G$ with an index set $I=\{1,2,\cdots,n\}$. 
We define $\mathfrak{g}={\rm Lie}(G)$ with Chevalley generators
$e_i$, $f_i$, $h_i$ ($i\in I$), a Cartan subalgebra $\mathfrak{h}$ and the canonical pairing $\langle, \rangle$ between
$\mathfrak{h}$ and $\mathfrak{h}^*$.
Let
$\Lambda_i$ denote the $i$-th fundamental weight, that is, $\langle h_j, \Lambda_i\rangle=\delta_{j,i}$
and $P=\oplus_{i\in I}\mathbb{Z}\Lm_i$ be the weight lattice, 
$P_+=\oplus_{i\in I}\mathbb{Z}_{\geq0}\Lm_i$ the positive weight lattice,
 $P^*=\oplus_{i\in I}\mathbb{Z}h_i$ the dual weight lattice,
$\{\alpha_i\}$ ($i\in I$) the set of simple roots.
For each $\lambda\in P_+$, let $V(\lambda)$ denote the finite dimensional irreducible $\mathfrak{g}$-module
with highest weight $\lambda$.
Let $U_q(\mathfrak{g})$
be the
quantized universal enveloping algebra with generators $E_i$, $F_i$ $(i\in I)$ and $K_{\lambda}$ ($\lambda \in P$)
and $U_q(\mathfrak{g})^-\subset U_q(\mathfrak{g})$ be the subalgebra generated by $\{F_i\}_{i\in I}$.
It is well-known that $U_q(\mathfrak{g})^-$ 
has the crystal base $(L(\infty),B(\infty))$. 
For two integers $l$, $m\in\mathbb{Z}$ such that $l\leq m$, one sets $[l,m]:=\{l,l+1,\cdots,m-1,m\}$.

\subsection{A birational map}\label{aoe}

One of main object in this paper is a variety $B^-_{w_0}$, where $w_0$ is the longest element in $W$.
Let us define $B^-_{w_0}$ and recall an open embedding $(\mathbb{C}^{\times})^{N}\hookrightarrow B^-_{w_0}$
associated with a reduced word $\mathbf{i}=(i_1,i_2,\cdots,i_N)$ of $w_0$.

First, for $i\in I$ and $t\in\mathbb{C}$, we put
\[
x_i(t):={\rm exp}(te_i),\ y_i(t):={\rm exp}(tf_i)\in G.
\]
There exists the canonical embedding $\phi_i : SL_{2}(\mathbb{C})\rightarrow G$ satisfying
\[
x_i(t)=\phi_i\left(
\begin{pmatrix}
1 & t \\
0 & 1
\end{pmatrix}
\right),\quad
y_i(t)=\phi_i\left(
\begin{pmatrix}
1 & 0 \\
t & 1
\end{pmatrix}
\right).
\]
Using the embedding, 
one puts
\[
t^{h_i}:=
\phi_i\left(
\begin{pmatrix}
t & 0 \\
0 & t^{-1}
\end{pmatrix}
\right)\in T
\]
and
\[
x_{-i}(t):=y_i(t)t^{-h_i}=
\phi_i\left(
\begin{pmatrix}
t^{-1} & 0 \\
1 & t
\end{pmatrix}
\right)\in G
\]
for $i\in I$ and $t\in\mathbb{C}^{\times}$.
One can construct a representative of a simple reflection $s_i\in W={\rm Norm}_G(T)/T$ by
\[
\overline{s_i}:=x_i(-1)y_i(1)x_i(-1)\in {\rm Norm}_G(T)
\]
for each $i\in I$.
For $w\in W$, one can define a representative $\overline{w}\in{\rm Norm}_G(T)$ by the rule
\[
\overline{uv}=\overline{u}\cdot \overline{v} \quad \text{ if } l(uv)=l(u)+l(v),
\]
where $l$ is the length function on $W$.
Now we can define a variety $B^-_{w_0}$ as $B^-_{w_0}:=B^-\cap U\overline{w_0}U $.
One defines a map
$\theta^-_{\mathbf{i}}: (\mathbb{C}^{\times})^{N}\rightarrow G$ associated with a reduced word
$\mathbf{i}=(i_1,\cdots,i_N)$ of $w_0\in W$
by
\begin{equation}\label{bir}
\theta^-_{\mathbf{i}}(t_1,\cdots,t_N):=x_{-i_1}(t_1)\cdots x_{-i_N}(t_N).
\end{equation}
\begin{prop}[\cite{BZ}]\label{OEprop}
The map $\theta^-_{\mathbf{i}}$ is an open embedding from $(\mathbb{C}^{\times})^{N}$ to $B^-_{w_0}$.
\end{prop}

\subsection{Generalized minors and $\mathbf{i}$-trails}

Let $G_0:=U^-TU\subset G$ denote the open subset whose elements $x\in G_0$ are uniquely decomposed as
$x=[x]_-[x]_0[x]_+$ with some $[x]_-\in U^-$, $[x]_0\in T$ and $[x]_+\in U$. 

\begin{defn}[\cite{FZ}]\label{gen-def}
For $u,v\in W$ and $i\in I$, the generalized minor $\Delta_{u\Lambda_i,v\Lambda_i}$ is defined as
the regular function on $G$
such that
\[
\Delta_{u\Lambda_i,v\Lambda_i}(x)=([\overline{u}^{-1}x\overline{v}]_0)^{\Lambda_i}
\]
for any $x\in\overline{u}G_0\overline{v}^{-1}$.
Here, for $t\in \mathbb{C}^{\times}$ and $j\in I$, we define $(t^{h_j})^{\Lambda_i}=(t^{\Lambda_i(h_j)})$ and extend it to
the group homomorphism $T\rightarrow \mathbb{C}^{\times}$.
\end{defn}

For calculations of generalized minors, one can use $\mathbf{i}$-{\it trails} \cite{BZ}. Here in this subsection, we take
$\mathbf{i}=(i_1,\cdots,i_l)$ as a sequence of indices from $I$. Let us review pre-$\mathbf{i}$-trails and $\mathbf{i}$-trails.

\begin{defn}\label{pretrail}
For a finite dimensional representation $V$ of $\mathfrak{g}$, two weights $\gamma$, $\delta$ of $V$ and
a sequence $\mathbf{i}=(i_1,\cdots,i_l)$ of indices from $I$,
a sequence
$\pi=(\gamma=\gamma_0,\gamma_1,\cdots,\gamma_l=\delta)$ 
is said to be a {\it pre-$\mathbf{i}$-trail} from $\gamma$ to $\delta$ if $\gamma_1,\cdots,\gamma_{l-1}\in P$ and
for $k\in[1,l]$, it holds $\gamma_{k-1}-\gamma_k=c_k\alpha_{i_k}$ with some nonnegative integer $c_k$.
\end{defn}
\nd
We can easily check that for $k\in[1,l]$, it holds
\begin{equation}\label{ck}
c_k=\frac{\gamma_{k-1}-\gamma_k}{2}(h_{i_k}).
\end{equation}

\begin{defn}[\cite{BZ}]\label{def-tr}
We consider the setting of Definition \ref{pretrail}.
If a pre-$\mathbf{i}$-trail $\pi$ from $\gamma$ to $\delta$ satisfies the condition
\begin{itemize}
\item $e^{c_1}_{i_1}e^{c_2}_{i_2}\cdots e^{c_l}_{i_l}$ is a non-zero linear map from $V_{\delta}$ to $V_{\gamma}$,
\end{itemize}
then $\pi$ is said to be an {\it $\mathbf{i}$-trail} from $\gamma$ to $\delta$,
where $V=\oplus_{\mu} V_{\mu}$ is the weight decomposition of $V$.
\end{defn}

\nd
For a pre-$\mathbf{i}$-trail $\pi=(\gamma_0,\gamma_1,\cdots,\gamma_l)$ and $k\in[1,l]$, we put
\begin{equation}\label{dk}
d_k(\pi):=\frac{\gamma_{k-1}+\gamma_k}{2}(h_{i_k}).
\end{equation}
One obtains $d_k(\pi)=c_k+\gamma_k(h_{i_k})\in\mathbb{Z}$ by (\ref{ck}).
If $\gamma_{k-1}=s_{i_k}\gamma_k$ then $d_k(\pi)=0$.

\begin{lem}[\cite{KKN}]\label{i-trail-lem}
Let $\gamma$, $\delta$ be weights of a finite dimensional representation $V$ of $\mathfrak{g}$.
Let $\mathbf{i}=(i_1,\cdots,i_l)$ be a sequence of indices from $I$ and $\pi=(\gamma_0,\gamma_1,\cdots,\gamma_l)$, $\pi'=(\gamma_0',\gamma_1',\cdots,\gamma_l')$
be two pre-$\mathbf{i}$-trails from $\gamma$ to $\delta$.
If $d_k(\pi)=d_k(\pi')$ for all $k\in[1,l]$ then $\pi=\pi'$.
\end{lem}

\nd
For a sequence $\mathbf{i}=(i_1,\cdots,i_l)$ of indices from $I$ and $t_1,\cdots,t_l\in\mathbb{C}^{\times}$,
just as in (\ref{bir}), we set 
\[
\theta^-_{\mathbf{i}}(t_1,\cdots,t_l):=
x_{-i_1}(t_1)\cdots x_{-i_l}(t_l)\in G.
\]
Then the following theorem holds:
\begin{thm}[\cite{BZ}]\label{trail-thm}
For $u$, $v\in W$ and $i\in I$, it holds
\[
\Delta_{u\Lambda_i,v\Lambda_i}(\theta^-_{\mathbf{i}}(t_1,\cdots,t_l))
=\sum_{\pi}C_{\pi}t_1^{d_1(\pi)}\cdots t_l^{d_l(\pi)},
\]
where $C_{\pi}$ is a positive integer and $\pi$ runs over all
$\mathbf{i}$-trails from $-u\Lambda_i$ to $-v\Lambda_i$ in $V(-w_0\Lambda_i)$. 
\end{thm}

By this theorem and Lemma \ref{i-trail-lem}, for each monomial $M$ in 
$\Delta_{u\Lambda_i,v\Lambda_i}(\theta^-_{\mathbf{i}}(t_1,\cdots,t_l))$, there uniquely exists a corresponding $\mathbf{i}$-trail
$\pi$ from $-u\Lambda_i$ to $-v\Lambda_i$
satisfying $M=t_1^{d_1(\pi)}\cdots t_l^{d_l(\pi)}$.

\section{Geometric crystals}

\subsection{Crystals}

Let us recall the definition of crystals following \cite{K2}.
Note that we use a slightly different notation from the original paper.
\begin{defn}[\cite{K2}]
A $6$-tuple $(\mathcal{B},
\{\tilde{e}_i\}_{i\in I},\{\tilde{f}_i\}_{i\in I},
\{\tilde{\gamma}_i\}_{i\in I},\{\tilde{\varepsilon}_i\}_{i\in I},\{\tilde{\varphi}_i\}_{i\in I})$ is called a $\mathfrak{g}$-{\it crystal} 
if 
$\mathcal{B}$ is a set and the maps
$\tilde{\gamma}_i:\mathcal{B}\rightarrow \mathbb{Z}$,
$\tilde{\varepsilon}_i,\tilde{\varphi}_i:\mathcal{B}\rightarrow \mathbb{Z}\sqcup \{-\infty\}$
and $\tilde{e}_i$,$\tilde{f}_i:\mathcal{B}\rightarrow \mathcal{B}\sqcup\{0\}$
($i\in I$) satisfy the following conditions: For $b,b'\in\mathcal{B}$, $i,j\in I$,
\begin{enumerate}
\item[$(1)$] $\tilde{\varphi}_i(b)=\tilde{\varepsilon}_i(b)+\tilde{\gamma}_i(b)$,
\item[$(2)$] $\tilde{\gamma}_j(\tilde{e}_ib)=\tilde{\gamma}_j(b)+a_{j,i}$ if $\tilde{e}_i(b)\in\mathcal{B}$,
\quad $\tilde{\gamma}_j(\tilde{f}_ib)=\tilde{\gamma}_j(b)-a_{j,i}$ if $\tilde{f}_i(b)\in\mathcal{B}$,
\item[$(3)$] $\tilde{\varepsilon}_i(\tilde{e}_i(b))=\tilde{\varepsilon}_i(b)-1,\ \ 
\tilde{\varphi}_i(\tilde{e}_i(b))=\tilde{\varphi}_i(b)+1$\ if $\tilde{e}_i(b)\in\mathcal{B}$, 
\item[$(4)$] $\tilde{\varepsilon}_i(\tilde{f}_i(b))=\tilde{\varepsilon}_i(b)+1,\ \ 
\tilde{\varphi}_i(\tilde{f}_i(b))=\tilde{\varphi}_i(b)-1$\ if $\tilde{f}_i(b)\in\mathcal{B}$, 
\item[$(5)$] $\tilde{f}_i(b)=b'$ if and only if $b=\tilde{e}_i(b')$,
\item[$(6)$] if $\tilde{\varphi}_i(b)=-\infty$ then $\tilde{e}_i(b)=\tilde{f}_i(b)=0$.
\end{enumerate}
Here
$0$ and $-\infty$ are additional elements, which do not belong to $\mathcal{B}$ and $\mathbb{Z}$, respectively.
The maps $\tilde{e}_i$,$\tilde{f}_i$ are called {\it Kashiwara operators}.
We say
a crystal $\mathcal{B}$ is {\it free} if the Kashiwara operators
$\tilde{e}_i$ $(i\in I)$ are bijections $\tilde{e}_i:\mathcal{B}\rightarrow \mathcal{B}$. 
\end{defn}

\subsection{Geometric crystals and tropicalization functor}\label{GC-trop}

In this subsection, first, let us review the definition of geometric crystals, which is introduced as a
geometric analog of crystals \cite{BK}. Note that
the following definitions of $\varepsilon_i$, $\varphi_i$ are different from \cite{BK}.
Applying flip $\varepsilon_i\mapsto \varepsilon_i^{-1}$, $\varphi_i\mapsto \varphi_i^{-1}$, they coincide
with the following one. We follow the notation in \cite{KN, KKN}.

\begin{defn}[\cite{BK}]
Let $X$ be an irreducible algebraic variety over $\mathbb{C}$ with
rational functions $\gamma_i$, $\varepsilon_i$ on $X$ and
unital rational $\mathbb{C}^{\times}$ actions
$\ovl{e}_i : \mathbb{C}^{\times}\times X\rightarrow X$ on $X$
for $i\in I$.
The quadruple $(X,\{\ovl{e}_i\}_{i\in I}, \{\gamma_i\}_{i\in I}, \{\varepsilon_i\}_{i\in I})$ is said to be a
$\mathfrak{g}$-{\it geometric crystal} if the following holds:
\begin{enumerate}
\item $(\{1\}\times X)\cap {\rm dom}(\ovl{e}_i)$ is open dense in $\{1\}\times X$ for each $i\in I$, where
${\rm dom}(\ovl{e}_i)$
is the domain of definition of $\ovl{e}_i : \mathbb{C}^{\times}\times X\rightarrow X$.
We write $\ovl{e}_i^c(x):=\ovl{e}_i(c,x)$ for $c\in \mathbb{C}^{\times}$ and $x\in X$.
\item $\gamma_j(\ovl{e}_i^c(x))=c^{a_{i,j}}\gamma_j(x)$ for any $i,j\in I$ and $c\in \mathbb{C}^{\times}$.
\item For any $i,j\in I$ and $c_1,c_2\in\mathbb{C}^{\times}$, it holds
\begin{eqnarray*}
\ovl{e}_i^{c_1}\ovl{e}_j^{c_2}&=&\ovl{e}_j^{c_2}\ovl{e}_i^{c_1} \quad \text{ if }a_{i,j}=0,\\
\ovl{e}_i^{c_1}\ovl{e}_j^{c_1c_2}\ovl{e}_i^{c_2}
&=&\ovl{e}_j^{c_2}\ovl{e}_i^{c_1c_2}\ovl{e}_j^{c_1} \quad \text{ if }a_{i,j}=a_{j,i}=-1,\\
\ovl{e}_i^{c_1}\ovl{e}_j^{c_1^2c_2}\ovl{e}_i^{c_1c_2}\ovl{e}_j^{c_2}
&=&\ovl{e}_j^{c_2}\ovl{e}_i^{c_1c_2}\ovl{e}_j^{c_1^2c_2}\ovl{e}_i^{c_1} \quad \text{ if }a_{i,j}=-2,\ a_{j,i}=-1,\\
\ovl{e}_i^{c_1}\ovl{e}_j^{c_1^3c_2}\ovl{e}_i^{c_1^2c_2}\ovl{e}_j^{c_1^3c_2^2}\ovl{e}_i^{c_1c_2}\ovl{e}_j^{c_2}
&=&\ovl{e}_j^{c_2}\ovl{e}_i^{c_1c_2}\ovl{e}_j^{c_1^3c_2^2}\ovl{e}_i^{c_1^2c_2}\ovl{e}_j^{c_1^3c_2}\ovl{e}_i^{c_1} \quad 
\text{ if }a_{i,j}=-3,\ a_{j,i}=-1.
\end{eqnarray*}
\item For any $i,j\in I$ and $c\in\mathbb{C}^{\times}$,
\[
\varepsilon_i(\ovl{e}_i^c(x))=c^{-1}\varepsilon_i(x),\ \quad
\varepsilon_i(\ovl{e}_j^c(x))=\varepsilon_i(x)\ \text{ if }a_{i,j}=0.
\]
\end{enumerate}
We set $\varphi_i(x):=\gamma_i(x)\varepsilon_i(x)$ for $i\in I$.
\end{defn}

\begin{defn}
For
a $\mathfrak{g}$-geometric crystal
$\chi=(X,\{\ovl{e}_i\}_{i\in I}, \{\gamma_i\}_{i\in I}, \{\varepsilon_i\}_{i\in I})$ and a rational function $f:X\rightarrow \mathbb{C}$,
the pair $(\chi,f)$ is said to be a $\mathfrak{g}$-upper (resp. lower) half-decorated geometric crystal if we have
\begin{equation}\label{hldec}
f(\ovl{e}_i^c(x))=f(x)+(c^{-1}-1)\varepsilon_i(x)\quad
(\text{resp. } f(\ovl{e}_i^c(x))=f(x)+(c-1)\varphi_i(x))
\end{equation}
for any $i\in I$, $c\in\mathbb{C}^{\times}$ and $x\in X$. Then we say $f$ is an
upper (resp. lower) half-decoration or upper (resp. lower) half-potential.
\end{defn}

Next, we recall a functor Trop, which makes a connection between crystals and geometric crystals.
Let
$^t\mathbb{Z}:=\mathbb{Z}\sqcup\{-\infty\}$ be the tropical semi-field with the multiplication `$+$'
and summation `min'. We define a map $\mathbb{V}:\mathbb{C}(x)\rightarrow$$^t\mathbb{Z}$ by
\[
\mathbb{V}(f(x)):=
\begin{cases}
-{\rm deg}(f(x^{-1})) & {\rm if}\ f\not\equiv0,\\
-\infty & {\rm if}\ f\equiv0,
\end{cases}
\]
where for polynomials $g(x)=\sum^m_{j=0}a_jx^j\in\mathbb{C}[x]$ and $h(x)=\sum^{m'}_{j=0}b_jx^j\in\mathbb{C}[x]\setminus \{0\}$
such that $a_m\neq0$, $b_m\neq0$,
we set ${\rm deg}(g(x)/h(x)):=m-m'$. For an algebraic torus $T'$ over $\mathbb{C}$, let $X^*(T'):={\rm Hom}(T',\mathbb{C}^{\times})$ and
$X_*(T'):={\rm Hom}(\mathbb{C}^{\times},T')$ denote the sets of characters and cocharacters, respectively. 

\begin{defn}
\begin{enumerate}
\item[$(i)$] For an algebraic torus $T'$,
a rational function $f:T'\rightarrow \mathbb{C}$ is said to be positive
if $f$ is in the form
\[
f=\frac{g}{h}
\]
with some regular functions $g=\sum_{\mu\in X^*(T')}a_{\mu}\mu$ and 
$h=\sum_{\mu\in X^*(T')}b_{\mu}\mu (\neq0)$
such that all $a_{\mu},\ b_{\mu}$ are nonnegative integers. 
\item[$(ii)$]
Let $f:T'\rightarrow T''$ be a rational map between
algebraic tori $T'$, $T''$. 
If $\xi\circ f : T'\rightarrow \mathbb{C}$ is a positive rational function for arbitrary $\xi\in X^*(T'')$
then the map $f$ is said to be positive.
\end{enumerate}
\end{defn}
For positive rational functions $f_1,f_2\in\mathbb{C}(x)\setminus\{0\}$ on $\mathbb{C}^{\times}$, it follows
\[
\mathbb{V}(f_1\cdot f_2)=\mathbb{V}(f_1)+\mathbb{V}(f_2),
\qquad
\mathbb{V}(f_1/ f_2)=\mathbb{V}(f_1)-\mathbb{V}(f_2),
\]
\[
\mathbb{V}(f_1+ f_2)={\rm min}(\mathbb{V}(f_1),\mathbb{V}(f_2)).
\]
\begin{defn}
Let $T'$, $T''$ be two algebraic tori, $f:T'\rightarrow T''$ be a rational map. 
Let $\langle,\rangle$ denote the pairing between $X^*(T'')$ and $X_*(T'')$.
We define a map
$\widehat{f}:X_*(T')\rightarrow X_*(T'')$ as
\[
\langle\chi,\widehat{f}(\xi)\rangle = \mathbb{V}(\chi\circ f\circ\xi)
\]
for $\chi\in X^*(T''), \xi\in X_*(T')$.
\end{defn} 

Let $\mathcal{T}_+$ denote the category whose objects are algebraic tori over $\mathbb{C}$
and morphisms are positive rational maps and
$\mathfrak{S}\mathfrak{e}\mathfrak{t}$ denote the category of sets.
A functor ${\rm Trop} : \mathcal{T}_+ \rightarrow \mathfrak{S}\mathfrak{e}\mathfrak{t}$ is defined
as
\[
T'\mapsto X_*(T'), \quad (f:T'\rightarrow T'') \mapsto (\widehat{f} : X_*(T')\rightarrow X_*(T'')).
\]
The functor ${\rm Trop}$ is said to be a {\it tropicalization}. For instance, if $T'=(\mathbb{C}^{\times})^3$,
$T''=\mathbb{C}^{\times}$ and $f:T'\rightarrow T''$ is defined as
\[
f(x_1,x_2,x_3)=x_1^3+2\frac{x_1x_2}{x_3}+4\frac{x_2^2x_3^3}{x_1^2}
\]
then ${\rm Trop}(f):X_*((\mathbb{C}^{\times})^3)\rightarrow X_*(\mathbb{C}^{\times})$ is given by
\[
{\rm Trop}(f)(z_1,z_2,z_3)={\rm min}\{3z_1,z_1+z_2-z_3,2z_2+3z_3-2z_1\}.
\]
where we identify $X_*((\mathbb{C}^{\times})^3)$, $X_*(\mathbb{C}^{\times})$ with $\mathbb{Z}^3$, $\mathbb{Z}$, respectively.
In this way, the product $\times$, sum $+$ and division $\div$ in $f$ correspond to
the sum $+$, ${\rm min}$ and minus $-$ in ${\rm Trop}(f)$ respectively.

\subsection{Tropicalizations of geometric crystals}\label{Trop-geom}

\begin{defn}
Let $\chi=(X,\{\ovl{e}_i\}_{i\in I}, \{\gamma_i\}_{i\in I}, \{\varepsilon_i\}_{i\in I})$ be
a $\mathfrak{g}$-geometric crystal, $T'$ an algebraic torus and $\theta:T'\rightarrow X$ a birational map.
We suppose the following (i), (ii):
\begin{enumerate}
\item[(i)] For each $i\in I$, rational functions $\gamma_i\circ \theta:T'\rightarrow \mathbb{C}$ and
$\varepsilon_i\circ \theta:T'\rightarrow \mathbb{C}$ are positive.
\item[(ii)] For each $i\in I$, the rational map $\ovl{e}_{i,\theta}:\mathbb{C}^{\times}\times T'\rightarrow \mathbb{C}^{\times}$
defined by $(c,t)\mapsto \theta^{-1}\circ \ovl{e}_i^c\circ\theta(t)$ is positive.
\end{enumerate}
Then we say $\theta$ is a {\it positive structure} on $\chi$.
If $f$ is an upper (resp. lower) half decoration of $\chi$ and a positive structure $\theta$ satisfies the condition
\begin{enumerate}
\item[(iii)] the rational function $f\circ \theta : T'\rightarrow \mathbb{C}$ is positive
\end{enumerate}
then $\theta$ is said to be a positive structure on upper (resp. lower) half-decorated geometric crystal $(\chi,f)$.
\end{defn}

The following theorem gives a connection between geometric crystals with positive structures
and crystals:

\begin{thm}[\cite{BK0}]
Let $\chi=(X,\{\ovl{e}_i\}_{i\in I}, \{\gamma_i\}_{i\in I}, \{\varepsilon_i\}_{i\in I})$ be
a $\mathfrak{g}$-geometric crystal with a positive structure $\theta:T'\rightarrow X$.
We define
\[
\tilde{e}_i:={\rm Trop}(\ovl{e}_{i,\theta}):\mathbb{Z}\times X_*(T')\rightarrow X_*(T'),
\]
\[
\tilde{\gamma}_i:={\rm Trop}(\gamma_i\circ \theta):X_*(T')\rightarrow\mathbb{Z},\quad 
\tilde{\varepsilon}_i:={\rm Trop}(\varepsilon_i\circ \theta):X_*(T')\rightarrow\mathbb{Z}.
\]
Then the $6$-tuple $(X_*(T'),\{\tilde{e}_i(1,\cdot)\}_{i\in I},\{\tilde{e}_i(-1,\cdot)\}_{i\in I},
\{\tilde{\gamma}_i\}_{i\in I},
\{\tilde{\varepsilon}_i\}_{i\in I},
\{\tilde{\varphi}_i\}_{i\in I})$ is a free $^L\mathfrak{g}$-crystal.
Here, $^L\mathfrak{g}$ is the finite dimensional simple Lie algebra whose Cartan matrix
is $^tA$ and $\tilde{\varphi}_i:=\tilde{\gamma}_i+\tilde{\varepsilon}_i$.
We write $\tilde{e}_i^z(x)=\tilde{e}_i(z,x)$ and $\tilde{f}_i^z(x)=\tilde{e}_i(-z,x)$ for $i\in I$, $z\in\mathbb{Z}_{\geq0}$ and $x\in X_*(T')$
and write $\tilde{e}_i^1=\tilde{e}_i$, $\tilde{f}_i^1=\tilde{f}_i$.
\end{thm}

Let 
$\theta : T'\rightarrow X$ be a positive
structure on
a $\mathfrak{g}$-upper (resp. lower) decorated geometric crystal 
$(\chi,f)$.
We consider a subset
\[
\tilde{B}_{\theta,f}:=\{x\in X_*(T') | {\rm Trop}(f\circ \theta)(x)\geq0\} \subset X_*(T').
\]
Defining $\tilde{e}_i(x)=0$ (resp. $\tilde{f}_i(x)=0$) if $\tilde{e}_i(x)\notin \tilde{B}_{\theta,f}$
(resp. $\tilde{f}_i(x)\notin \tilde{B}_{\theta,f}$), we see that
$\tilde{B}_{\theta,f}$ has a $^L\mathfrak{g}$-crystal structure
\begin{equation}\label{gccrystal}
\mathbb{B}_{\theta,f}:=(\tilde{B}_{\theta,f},\{\tilde{e}_i|_{\tilde{B}_{\theta,f}}\}_{i\in I},
\{\tilde{f}_i|_{\tilde{B}_{\theta,f}}\}_{i\in I},\{\tilde{\varepsilon}_i|_{\tilde{B}_{\theta,f}}\}_{i\in I},
\{\tilde{\varphi}_i|_{\tilde{B}_{\theta,f}}\}_{i\in I},
\{\tilde{\gamma}_i|_{\tilde{B}_{\theta,f}}\}_{i\in I}).
\end{equation}

\begin{prop}[\cite{KN}]
In the above setting,
$\mathbb{B}_{\theta,f}$ is an upper (resp. lower) normal crystal,
that is,
\[
\tilde{\varepsilon}_i(x)={\rm max}\{n\geq0|\tilde{e}_i^n(x)\neq0\},\quad
(\text{resp. } \tilde{\varphi}_i(x)={\rm max}\{n\geq0|\tilde{f}_i^n(x)\neq0\}).
\]

\end{prop}
\nd
By this proposition, the condition (\ref{hldec}) on an upper/lower half-decorated geometric crystal corresponds to
the upper/lower normality of the corresponding crystal.

\subsection{Geometric crystal structure on $B^-_{w_0}$}

Defining maps
\[
\gamma_i: B^-_{w_0}\rightarrow\mathbb{C}^{\times},\quad
\varepsilon_i: B^-_{w_0}\rightarrow\mathbb{C}^{\times},\quad
\ovl{e}_i:\mathbb{C}^{\times}\times B^-_{w_0}\rightarrow B^-_{w_0}
\]
on $B^-_{w_0}=B^-\cap U\overline{w_0}U$,
we get a $\mathfrak{g}$-geometric crystal $(B^-_{w_0},\{\ovl{e}_i\}_{i\in I},\{\gamma_i\}_{i\in I},\{\varepsilon_i\}_{i\in I})$ \cite{BK}.
We do not recall the definition
of maps $\gamma_i$, $\varepsilon_i$ and 
$\ovl{e}_i$ since they are not needed in this paper. 
For the definition of maps, refer to Sect.3 of our previous paper \cite{KKN}.
Let us define a regular function $\Phi^{h}_{\rm BK}$ on $B^-_{w_0}$ as follows:
\[
\Phi^{h}_{\rm BK}:=
\sum_{i\in I} \Delta_{w_0\Lm_i,s_i\Lm_i}.
\]
In \cite{KN}, we proved that
the function $\Phi^{h}_{\rm BK}$ is an upper half-decoration on the geometric crystal $B^-_{w_0}$.

Recall an open embedding
$\theta^-_{\mathbf{i}}:(\mathbb{C}^{\times})^{N}\hookrightarrow B^-_{w_0}$ in
Proposition \ref{OEprop}, which gives a positive structure on $(B^-_{w_0},\Phi^{h}_{\rm BK})$.
Thus, one obtains a crystal $\mathbb{B}_{\theta^-_{\mathbf{i}},\Phi^{h}_{\rm BK}}$ 
as in (\ref{gccrystal}):
\[
\tilde{B}_{\theta^-_{\mathbf{i}},\Phi^{h}_{\rm BK}}:=
\{z\in X_{*}((\mathbb{C}^{\times})^{N})| {\rm Trop}(\Phi^{h}_{\rm BK}\circ \theta^-_{\mathbf{i}})(z)\geq0 \},
\]
\begin{equation}\label{boldb}
\mathbb{B}_{\theta^-_{\mathbf{i}},\Phi^{h}_{\rm BK}}=
(\tilde{B}_{\theta^-_{\mathbf{i}},\Phi^{h}_{\rm BK}},\{\tilde{e}_i\}_{i\in I},
\{\tilde{f}_i\}_{i\in I},
\{\tilde{\varepsilon}_i\}_{i\in I},
\{\tilde{\varphi}_i\}_{i\in I},
\{\tilde{\gamma}_i\}_{i\in I}).
\end{equation}
Here, we omitted the notation of restrictions $|_{\tilde{B}_{\theta^-_{\mathbf{i}},\Phi^{h}_{\rm BK}}}$
for $\tilde{e}_i$, $\tilde{f}_i$, $\tilde{\varepsilon}_i$, $\tilde{\varphi}_i$ and $\tilde{\gamma}_i$.

\begin{thm}[\cite{KN}]\label{thm1a}
For each reduced word $\mathbf{i}$ of the longest element $w_0$,
the set
$\mathbb{B}_{\theta^-_{\mathbf{i}},\Phi^{h}_{\rm BK}}$
is a $^L\mathfrak{g}$-crystal isomorphic to the crystal
$B(\infty)$.
\end{thm}


\subsection{Properties of $\Phi^{h}_{\rm BK}$}

In this subsection, we recall the properties of the functions $\Delta_{w_0\Lm_i,s_i\Lm_i}$,
which are summand of $\Phi^{h}_{\rm BK}$. We fix a reduced word $\mathbf{i}=(i_1,\cdots,i_N)$
of the longest element $w_0\in W$. Note that 
$\Delta_{w_0\Lm_i,s_i\Lm_i}\circ \theta^-_{\mathbf{i}}(t_1,\cdots,t_N)$ is a Laurent polynomial
with positive integer coefficients
by Theorem \ref{trail-thm}. 

\begin{prop}[\cite{KN}]\label{low-prop}
For any $i\in I$, the Laurent polynomial
$\Delta_{w_0\Lm_i,s_i\Lm_i}\circ \theta^-_{\mathbf{i}}(t_1,\cdots,t_N)$
has a term
\begin{equation}\label{lowest-term}
t_{J} t_{J+1}^{a_{i_{J+1},i}}\cdots t_{N}^{a_{i_N,i}},
\end{equation}
where $J:={\rm max}\{1\leq k\leq N | i_k=i\}$.
\end{prop}

\begin{prop}[\cite{KKN}]\label{high-prop}
For $i\in I$, we take $k\in[1,N]$ such that
$s_{i_N}s_{i_{N-1}}\cdots s_{i_{k+1}}\alpha_{i_k}=\alpha_i$.
Then
the Laurent polynomial
$\Delta_{w_0\Lm_i,s_i\Lm_i}\circ \theta^-_{\mathbf{i}}(t_1,\cdots,t_N)$
has a term
\begin{equation}\label{highest-term}
t_k.
\end{equation}
\end{prop}

\nd
For $j\in[1,N]$, we define
\begin{equation}\label{jpm}
\begin{split}
j^+:={\rm min}\{ l\in[1,N] | i_l=i_j,\ l>j \}\cup\{N+1\},\\
j^-:={\rm max}\{ l\in[1,N] | i_l=i_j,\ l<j \}\cup\{0\}.
\end{split}
\end{equation}
For $m\in\mathbb{Z}_{\geq1}$, one also defines $j^{0+}:=j$, $j^{0-}:=j$ and
\begin{equation}\label{jmp}
j^{m+}:=(j^{(m-1)+})^+,\qquad j^{m-}:=(j^{(m-1)-})^-,
\end{equation}
where we set $(N+1)^+:=N+1$, $0^-:=0$.
For $j\in[1,N]$ such that $j^+\leq N$, we define
\begin{equation}\label{ajdef}
A_{j}:=t_jt_{j^+}\prod_{j<l<j^+}t_l^{a_{i_l,i_j}}.
\end{equation}

\nd
The next proposition is shown in \cite{KKN} under the assumption $V(\Lambda_i)$ is a minuscule representation,
however, the proof is valid for arbitrary $i\in I$.

\begin{prop}\label{trail-prop}
Let $\pi=(\gamma_0,\gamma_1,\cdots,\gamma_N)$, $\pi'=(\gamma_0',\gamma_1',\cdots,\gamma_N')$ be two pre-$\mathbf{i}$-trails
from $-w_0\Lambda_i$ to $-s_i\Lambda_i$. We take integers $c_l$, $c_l'$ $(l=1,2,\cdots,N)$ as
\[
\gamma_{l-1}-\gamma_l=c_l\alpha_{i_l},\quad \gamma_{l-1}'-\gamma_l'=c_l'\alpha_{i_l}
\]
and assume that there exists $j\in[1,N]$ such that $j^+\leq N$ and
\[
c_l'=c_{l}\quad \text{for}\ l\in[1,N]\setminus \{j,j^+\},
\]
\[
c_j'=c_{j}+1,\quad c_{j^+}'=c_{j^+}-1.
\]
Then we have
\[
t_1^{d_1(\pi')}\cdots t_N^{d_N(\pi')}=t_1^{d_1(\pi)}\cdots t_N^{d_N(\pi)}A^{-1}_{j},
\]
where $d_j$ is defined as (\ref{dk}).
\end{prop}

The following lemma is a generalization of Lemma 5.6 in \cite{KKN}:
\begin{lem}\label{lem2}
If  a monomial $M$ appearing in
$\Delta_{w_0\Lm_i,s_i\Lm_i}\circ \theta^-_{\mathbf{i}}(t_1,\cdots,t_N)$
belongs to $\{ \prod^N_{l=1} t_l^{d_l} | d_l\in\mathbb{Z}_{\geq0} \}$ then
\[
M=t_k,
\]
where $k\in[1,N]$ is the same one as in Proposition \ref{high-prop}.
\end{lem}

\nd
{\it Proof.}

Let $\pi=(\gamma_0,\gamma_1,\cdots,\gamma_N)$ be the $\mathbf{i}$-trail from $-w_0\Lambda_i$ to $-s_i\Lambda_i$ in $V(-w_0\Lambda_i)$
corresponding to $M$ with integers $\{c_l\}_{l\in[1,N]}$ such that
$\gamma_{l-1}-\gamma_l=c_l\alpha_{i_l}$. It holds $M=\prod^N_{j=1}t_j^{d_j(\pi)}$ (see Theorem \ref{trail-thm}).
For $l\in[k,N]$, we show
\begin{equation}\label{lem2-pr-1}
\gamma_l=-s_{i_{l+1}}\cdots s_{i_N}s_i\Lambda_i
\end{equation}
by downward induction on $l$. The case $l=N$ is clear because $\gamma_N=-s_i\Lambda_i$.
We assume 
\[
\gamma_{r}=-s_{i_{r+1}}\cdots s_{i_N}s_i\Lambda_i
\]
for $r\in[l+1,N]$. This assumption implies if $r>l+1$ then $\gamma_{r-1}=s_{i_r}\gamma_r=\gamma_r-\gamma_r(h_{i_r})\alpha_{i_r}$ so that
$c_{r}=\frac{\gamma_{r-1}-\gamma_{r}}{2}(h_{i_{r}})=-\gamma_{r}(h_{i_{r}})$. Let $v_{-w\Lambda_i}$ denote
an extremal weight vector in $V(-w_0\Lambda_i)$ of weight $-w\Lambda_i$ for $w\in W$.
Using a property of extremal weight vectors and fact that
$s_{i_{l+1}}\cdots s_{i_N}s_i$ is reduced (Lemma 3.11 of \cite{Kac}), one obtains
\[
c_r={\rm max}\{m\in\mathbb{Z}_{\geq0} | e_{i_r}^{m} e_{i_{r+1}}^{c_{r+1}}\cdots e_{i_{N}}^{c_{N}} v_{-s_i\Lambda_i}\neq0 \} \quad (r\in[l+2,N]).
\]
Since $s_{i_{l+1}}\cdots s_{i_N}s_i$ is reduced, 
it holds
\[
\gamma_{l+1}(h_{i_{l+1}})=
\langle h_{i_{l+1}}, -s_{i_{l+2}}\cdots s_{i_N}s_i\Lambda_i
\rangle 
=
-\langle s_is_{i_N}\cdots s_{i_{l+2}}h_{i_{l+1}}, \Lambda_i
\rangle 
\in\mathbb{Z}_{\leq0}.
\]
By $\gamma_{l+1}(h_{i_{l+1}})=\langle h_{i_{l+1}}, -s_{i_{l+2}}\cdots s_{i_N}s_i\Lambda_i \rangle$
and
$e_{i_{l+2}}^{c_{l+2}}\cdots e_{i_{N}}^{c_{N}} v_{-s_i\Lambda_i}=C v_{-s_{i_{l+2}}\cdots s_{i_{N}}s_i\Lambda_i}$ with some
$C\in\mathbb{C}\setminus\{0\}$, we obtain
\[
-\gamma_{l+1}(h_{i_{l+1}})={\rm max}\{m\in\mathbb{Z}_{\geq0} | e_{i_{l+1}}^me_{i_{l+2}}^{c_{l+2}}\cdots e_{i_{N}}^{c_{N}} v_{-s_i\Lambda_i}\neq0
\}.
\]
Since
$e_{i_{l+1}}^{c_{l+1}} e_{i_{l+2}}^{c_{l+2}}\cdots e_{i_{N}}^{c_{N}} v_{-s_i\Lambda_i}\neq0$, it follows
$c_{l+1}\leq -\gamma_{l+1}(h_{i_{l+1}})$. Combining the assumption $d_{l+1}\geq0$ with
$d_{l+1}=\frac{\gamma_{l}(h_{i_{l+1}})+\gamma_{l+1}(h_{i_{l+1}})}{2}=c_{l+1}+\gamma_{l+1}(h_{i_{l+1}})\leq 0$, we have
$c_{l+1}+\gamma_{l+1}(h_{i_{l+1}})=0$ and $c_{l+1}=-\gamma_{l+1}(h_{i_{l+1}})$.
Therefore, one gets
\[
\gamma_l= {\rm wt}(e_{i_{l+1}}^{c_{l+1}}e_{i_{l+2}}^{c_{l+2}}\cdots e_{i_{N}}^{c_{N}} v_{-s_i\Lambda_i})
=-s_{i_{l+1}}\cdots s_{i_N}s_i\Lambda_i
\]
and (\ref{lem2-pr-1}) is shown. In particular, it holds
\begin{equation}\label{lem2-pr-2}
c_l={\rm max}\{m\in\mathbb{Z}_{\geq0} | e_{i_l}^{m} e_{i_{l+1}}^{c_{l+1}}\cdots e_{i_{N}}^{c_{N}} v_{-s_i\Lambda_i}\neq0 \} \quad (l\in[k+1,N])
\end{equation}
and
$\gamma_k=-s_{i_{k+1}}\cdots s_{i_N}s_i\Lambda_i$.
By
\[
s_{i_{k+1}}\cdots s_{i_N}s_i\Lambda_i
=s_{i_{k+1}}\cdots s_{i_N}(\Lambda_i-\alpha_i)
=s_{i_{k+1}}\cdots s_{i_N}\Lambda_i - \alpha_{i_k}
\]
and
\[
\langle
h_{i_k}, s_{i_{k+1}}\cdots s_{i_N}\Lambda_i \rangle
=\langle s_{i_{N}}\cdots s_{i_{k+1}}h_{i_k}, \Lambda_i \rangle
=\langle h_i, \Lambda_i \rangle=1,
\]
it follows
\begin{equation}\label{pr3-3-3}
s_{i_{k+1}}\cdots s_{i_N}s_i\Lambda_i=s_{i_k}s_{i_{k+1}}\cdots s_{i_N}\Lambda_i,
\end{equation}
which yields
$\gamma_k=-s_{i_{k+1}}\cdots s_{i_N}s_i\Lambda_i=-s_{i_k}s_{i_{k+1}}\cdots s_{i_N}\Lambda_i$.
Hence there exists $C'\in\mathbb{C}\setminus\{0\}$ such that
\[
e_{i_{k+1}}^{c_{k+1}} e_{i_{k+2}}^{c_{k+2}}\cdots e_{i_{N}}^{c_{N}} v_{-s_i\Lambda_i}
=C' v_{-s_{i_k}\cdots s_{i_N}\Lambda_i}
\]
and we obtain $c_k=0$ since $e_{i_{k}}^{c_{k}}e_{i_{k+1}}^{c_{k+1}} e_{i_{k+2}}^{c_{k+2}}\cdots e_{i_{N}}^{c_{N}} v_{-s_i\Lambda_i}\neq0$,
which yields $\gamma_{k-1}=\gamma_k=-s_{i_k}s_{i_{k+1}}\cdots s_{i_N}\Lambda_i$.
Note that $e_{i_1}^{c_1}\cdots e_{i_{k-1}}^{c_{k-1}}e_{i_{k}}^{c_{k}}e_{i_{k+1}}^{c_{k+1}} e_{i_{k+2}}^{c_{k+2}}\cdots e_{i_{N}}^{c_{N}} v_{-s_i\Lambda_i}
=e_{i_1}^{c_1}\cdots e_{i_{k-1}}^{c_{k-1}}e_{i_{k+1}}^{c_{k+1}} e_{i_{k+2}}^{c_{k+2}}\cdots e_{i_{N}}^{c_{N}} v_{-s_i\Lambda_i}$
is the highest weight vector with weight $-w_0\Lambda_i$.
Thus, the values of $c_l$ ($l=1,\cdots,k-1$) are uniquely determined by
\[
c_l=
{\rm max}\{m\in\mathbb{Z}_{\geq0} | e_{i_{l}}^me_{i_{l+1}}^{c_{l+1}}\cdots e_{i_{N}}^{c_{N}} v_{-s_i\Lambda_i}\neq0
\}.
\]
The values of $\{c_r\}_{r\in[k+1,N]}$ and $c_k$ are uniquely determined as (\ref{lem2-pr-2}) and $c_k=0$.
In this way, the integers $\{c_l\}_{l\in[1,N]}$ are uniquely determined from the assumption of our claim.
Hence, the conclusion $M=t_k$ follows. \qed

\section{Main theorem}

In this section, 
let $N:=l(w_0)$ and we fix a reduced word $\mathbf{i}=(i_1,\cdots,i_N)$ of $w_0$ and $i\in I$.
We give an algorithm to compute all monomials in $\Delta_{w_0\Lm_i,s_i\Lm_i}\circ \theta^-_{\mathbf{i}}(t_1,\cdots,t_N)$
(subsection \ref{algo}).
In the algorithm, we will associate some integer vector $(b_1,b_2,\cdots,b_N)$ to each monomial.
Let us define integers $b_1,b_2,\cdots,b_N$ in the first subsection.

\subsection{$\textbf{b}$-integers}

\begin{defn}\label{bint}
Let $M=\prod_{l=1}^{N}t_l^{d_l}$ be a Laurent monomial.
We inductively define integers $\{b_l\}_{l=N,N-1,\cdots,1}$ as
\[
b_N=d_N+s_i\Lambda_i(h_{i_N}),
\]
\[
b_t=d_t+s_i\Lambda_i(h_{i_t})-\sum_{l=t}^{N-1}b_{l+1}a_{i_t,i_{l+1}} \ \ (t=N,N-1,\cdots,1).
\]
\end{defn}

\begin{lem}\label{lem-3}
Let $M=t_k$ in Proposition \ref{high-prop}. Then $b^0_t=b_t$ $(t\in[1,N])$ in the previous lemma are as follows:
\begin{eqnarray*}
b^0_t&=&s_i\Lambda_i(h_{i_t})-\sum_{l=t}^{N-1}b^0_{l+1}a_{i_t,i_{l+1}}\\
&=&\langle h_{i_t}, s_{i_{t+1}}s_{i_{t+2}}\cdots s_{i_{N}}s_i\Lambda_i \rangle \ \ (t=N,N-1,\cdots,k+1),
\end{eqnarray*}
\[
b^0_k=0,
\]
\begin{eqnarray*}
b^0_t&=&s_i\Lambda_i(h_{i_t})-\sum_{l=t}^{N-1}b^0_{l+1}a_{i_t,i_{l+1}}\\
&=&
\langle h_{i_t}, s_{i_{t+1}}s_{i_{t+2}}\cdots s_{i_{N}}\Lambda_i \rangle 
\ \ (t=k-1,k-2,\cdots,1).
\end{eqnarray*}
\end{lem}

\nd
{\it Proof.} For $t=N,N-1,\cdots,k+1$, it follows by $d_t=0$ that
\[
b^0_t=s_i\Lambda_i(h_{i_t})-\sum_{l=t}^{N-1}b^0_{l+1}a_{i_t,i_{l+1}}=\langle h_{i_t}, s_i\Lambda_i -\sum_{l=t}^{N-1}b^0_{l+1}\alpha_{i_{l+1}} \rangle.
\]
Thus,
\begin{equation}\label{lem-3-pr-1}
b^0_t={\rm max}\{l\in\mathbb{Z}_{\geq0} | f_{i_t}^l f_{i_{t+1}}^{b^0_{t+1}}f_{i_{t+2}}^{b^0_{t+2}}\cdots f_{i_{N}}^{b^0_{N}}v_{s_i\Lambda_i}\neq0 \}
\end{equation}
and
\[
b^0_t
=\langle h_{i_t}, s_i\Lambda_i -\sum_{l=t}^{N-1}b^0_{l+1}\alpha_{i_{l+1}} \rangle
=\langle h_{i_t}, s_{i_{t+1}}s_{i_{t+2}}\cdots s_{i_{N}}s_i\Lambda_i \rangle .
\]
Here, $v_{w\Lambda_i}$ is an extremal weight vector of weight $w\Lambda_i$ in $V(\Lambda_i)$ for $w\in W$.

By our assumption $d_k=1$, $s_{i_{N}}\cdots s_{i_{k+2}}s_{i_{k+1}}h_{i_k}=h_i$ and (\ref{lem-3-pr-1}), it holds
\begin{eqnarray*}
b^0_k&=&1+s_i\Lambda_i(h_{i_k})-\sum_{l=k}^{N-1}b^0_{l+1}a_{i_k,i_{l+1}}
=1+
\langle h_{i_k}, s_i\Lambda_i -\sum_{l=k}^{N-1}b^0_{l+1}\alpha_{i_{l+1}} \rangle\\
&=&
1+
\langle h_{i_k}, s_{i_{k+1}}s_{i_{k+2}}\cdots s_{i_{N}}s_i\Lambda_i \rangle\\
&=&
1+\langle h_i,s_i\Lambda_i\rangle=0.
\end{eqnarray*}
For $t=k-1,\cdots,1$, it follows by (\ref{pr3-3-3}), $d_t=0$ and $b^0_k=0$ that 
\begin{eqnarray*}
b^0_t&=&s_i\Lambda_i(h_{i_t})-\sum_{l=t}^{N-1}b^0_{l+1}a_{i_t,i_{l+1}}\\
&=&
\langle h_{i_t}, s_i\Lambda_i -\sum_{l=k}^{N-1}b^0_{l+1}\alpha_{i_{l+1}}
-\sum_{l=t}^{k-2}b^0_{l+1}\alpha_{i_{l+1}}
 \rangle \\
&=&\langle h_{i_t}, s_{i_{k+1}}s_{i_{k+2}}\cdots s_{i_{N}}s_i\Lambda_i
-\sum_{l=t}^{k-2}b^0_{l+1}\alpha_{i_{l+1}}
 \rangle\\
&=&\langle h_{i_t}, s_{i_k}s_{i_{k+1}}s_{i_{k+2}}\cdots s_{i_{N}}\Lambda_i
-\sum_{l=t}^{k-2}b^0_{l+1}\alpha_{i_{l+1}}
 \rangle.
\end{eqnarray*}
Therefore, we can verify
\begin{equation}\label{lem-3-pr-2}
b^0_t={\rm max}\{l\in\mathbb{Z}_{\geq0} | f_{i_t}^l f_{i_{t+1}}^{b^0_{t+1}}\cdots
f_{i_{k-1}}^{b^0_{k-1}} v_{s_{i_k}s_{i_{k+1}}\cdots s_{i_N}\Lambda_i}\neq0 \}
\end{equation}
and
\[
b^0_t=
\langle h_{i_t}, s_{i_{t+1}}\cdots s_{i_{k-1}}s_{i_k}s_{i_{k+1}}s_{i_{k+2}}\cdots s_{i_{N}}\Lambda_i
 \rangle.
\]
\qed

\begin{lem}\label{btct}
For a monomial 
$M=\prod_{l=1}^{N}t_l^{d_l}$ in $\Delta_{w_0\Lm_i,s_i\Lm_i}\circ \theta^-_{\mathbf{i}}(t_1,\cdots,t_N)$,
we take integers
$\{b_l\}_{l=N,N-1,\cdots,1}$ as
in Definition \ref{bint}.
Let $\pi=(\gamma_0,\gamma_1,\cdots,\gamma_N)$ be the $\mathbf{i}$-trail corresponding to $M$
with integers $c_l$ $(l=1,2,\cdots,N)$ such that
$\gamma_{l-1}-\gamma_l=c_l\alpha_{i_l}$. Then we get
\[
c_t=b_t
\]
for $t\in[1,N]$.
\end{lem}
\nd
{\it Proof.}
The proof proceeds by downward induction on $t=N,N-1,\cdots,1$.
By (\ref{ck}), (\ref{dk}) and $\gamma_N=-s_i\Lambda_i$, we see that
\[
c_N=\frac{\gamma_{N-1}-\gamma_N}{2}(h_{i_N})
=\frac{\gamma_{N-1}+\gamma_N}{2}(h_{i_N})-\gamma_N(h_{i_N})
=d_N-\gamma_N(h_{i_N})=d_N+s_i\Lambda_i(h_{i_N})=b_N.
\]
We assume that
for $l=t+1,t+2,\cdots,N$, it holds $b_l=c_l$.
One gets
\begin{equation}\label{btct-1}
c_t=
\frac{\gamma_{t-1}-\gamma_t}{2}(h_{i_t})
=
\frac{\gamma_{t-1}+\gamma_t}{2}(h_{i_t})
-\gamma_t(h_{i_t})
=d_t-\gamma_t(h_{i_t}).
\end{equation}
It follows from the expression
\[
\gamma_t
=\gamma_N+\sum_{l=t}^{N-1} (\gamma_l-\gamma_{l+1})
=\gamma_N+\sum_{l=t}^{N-1}c_{l+1}\alpha_{i_{l+1}},
\]
induction hypothesis
and (\ref{btct-1}) that
\[
c_t=
d_t-\left(\gamma_N(h_{i_t})+\sum_{l=t}^{N-1}c_{l+1}a_{i_t,i_{l+1}}\right)
=d_t+s_i\Lambda_i(h_{i_t})-\sum_{l=t}^{N-1}b_{l+1}a_{i_t,i_{l+1}}=b_t.
\]
\qed

\subsection{Main Theorem}\label{algo}

In this subsection, we take $i\in I$ such that
for any weight $\mu$ of $V(-w_0\Lambda_i)$ and $t\in I$, it holds
\begin{equation}\label{mini-pro}
\langle
h_t,\mu
\rangle\in\{2,1,0,-1,-2\}.
\end{equation}
The above condition is a natural generalization of minuscule representations since
$V(-w_0\Lambda_i)$ is minuscule if and only if $\langle h_t,\mu \rangle \in\{1,0,-1\}$
for all $t\in I$ and weight $\mu$.
The list of $i\in I$ satisfying the above condition (\ref{mini-pro}) is as follows:
\begin{table}[h]
  \begin{tabular}{|c|c|c|c|c|c|c|c|c|c|} \hline
  $\mathfrak{g}$ & ${\rm A}_n$ & ${\rm B}_n$ & ${\rm C}_n$ & ${\rm D}_n$ & ${\rm E}_6$ & ${\rm E}_7$ & ${\rm E}_8$ & ${\rm F}_4$ & ${\rm G}_2$\\ \hline
  $i$ & \text{all }$i\in I$ & \text{all }$i\in I$ & \text{all }$i\in I$ & \text{all }$i\in I$ & $1,2,4$, $5,6$ & $1,5,6,7$  & $1,7$ & $1,4$ & $\phi$ \\ \hline
  \end{tabular}
\end{table} 

\begin{thm}\label{thm2}
For $i\in I$, we take $k\in[1,N]$ such that
$s_{i_N}s_{i_{N-1}}\cdots s_{i_{k+1}}\alpha_{i_k}=\alpha_i$.
Then the set of monomials appearing in
$\Delta_{w_0\Lm_i,s_i\Lm_i}\circ \theta^-_{\mathbf{i}}(t_1,\cdots,t_N)$
coincides with the set of Laurent monomials labelling the vertices of the directed graph $\overline{DG}$ obtained
by the following algorithm:
\begin{enumerate}
\item[(1)] Let $\overline{DG}_0$ be the graph which has only one vertex $t_k$ and no arrow.
We associate the integer vector $(b_1^0,b_2^0,\cdots,b_N^0)$ in Lemma \ref{lem-3} to the monomial $t_k$.
\item[(2)] For each sink $M=\prod_{j=1}^Nt_j^{d_j}$ of $\overline{DG}_l$
associated with a vector $(b_1,b_2,\cdots,b_N)$
and for each $j\in[1,N]$ such that $j^+\leq N$, we add a
vertex $M\cdot A_{j}^{-1}$ and an arrow $M\rightarrow M\cdot A_{j}^{-1}$ to $\overline{DG}_l$
if and only if $d_j>0$, $b_{j^+}>0$ and one of the following conditions holds:
\begin{enumerate}
\item[(a)] $d_{j^+}< d_j$,
\item[(b)] $d_{j^+}= d_j$ and there exists
$p\in\mathbb{Z}_{\geq2}$ such that
\[
d_{j^{m+}}=b_{j^{m+}}=0 \quad(m=2,3,\cdots,p-1)
\]
and $d_{j^{p+}}=-1$, $b_{j^{p+}}=1$.
\end{enumerate}
We associate integers $(b'_1,b'_2,\cdots,b'_N)$ defined by $b'_s=b_s$ $(s\in[1,N]\setminus\{j,j^+\})$, $b_j'=b_j+1$, $b'_{j^+}=b_{j^+}-1$
to $M\cdot A_{j}^{-1}$.
Let $\overline{DG}_{l+1}$ be the directed graph obtained from $\overline{DG}_l$ by this step.
\item[(3)] If there is no monomial
$M=\prod_{j=1}^Nt_j^{d_j}$ in $\overline{DG}_r$ satisfying the condition in (2) then
we stop this algorithm and set $\overline{DG}=\overline{DG}_r$.
\end{enumerate}
\end{thm}
We call the graph $\overline{DG}$ {\it decoration graph}.
Note that two same monomials obtained from different sinks define the same vertex in step (2).
For a Laurent polynomial $F=\sum_{j=1}^s c_jM_j$ with $c_j\in\mathbb{Z}_{\geq1}$, $M_j\in\{\prod_{j=1}^N t_j^{d_j} | d_j\in\mathbb{Z} \}$,
{\it the set of monomials appearing in} $F$ means the set $\{M_j | j\in[1,s]\}$.

We see that the monomial $t_k$ in Proposition \ref{high-prop} is the unique source in $\overline{DG}$.
If $M=t_{J} t_{J+1}^{a_{i_{J+1},i}}\cdots t_{N}^{a_{i_N,i}}$ in Proposition \ref{low-prop} then
no $j\in[1,N]$ satisfies the condition in (2) so that $t_{J} t_{J+1}^{a_{i_{J+1},i}}\cdots t_{N}^{a_{i_N,i}}$ is a sink
in $\overline{DG}$. There are cases the graph $\overline{DG}$ has several sinks (see Example \ref{ex2}).

\subsection{Examples}

\begin{ex}\label{ex-1}
Let $G$ be of type ${\rm C}_3$, $\mathbf{i}=(2,3,2,1,2,3,2,3,1)$, $i=2$.
Following the algorithm in previous theorem, we compute the graph $\overline{DG}_0,\overline{DG}_1,\cdots,
\overline{DG}_r=\overline{DG}$.
\vspace{2mm}

\hspace{-7mm}
\underline{$\overline{DG}_0$} By $s_1s_3s_2s_3s_2s_1s_2s_3\alpha_2=\alpha_2$,
we see that $\overline{DG}_0$ is the graph which has only one vertex $t_1$ and no arrow.
Following Lemma \ref{lem-3}, the monomial $t_1$ is associated with
\[
(b_1^0,b_2^0,\cdots,b_9^0)=(0,0,1,1,0,1,2,1,1).
\]
\underline{$\overline{DG}_1$} Let $d_j$ be the exponent of $t_j$ in the monomial $t_1$ so that
$(d_1,d_2,\cdots,d_9)=(1,0,0,0,0,0,0,0,0)$.
Focusing on the unique sink $t_1$ in $\overline{DG}_0$, one can verify
$d_1>0$, $b_{1^+}^0=b_{3}^0=1>0$ and $d_{1^+}=d_3=0<d_1$. Thus the pair $t_1$ and $j=1$ satisfy (a) of (2).
Hence, the graph $\overline{DG}_1$ is as follows:
\[
\overline{DG}_1\ :\ 
t_1\rightarrow t_1\cdot A_{1}^{-1}=\frac{t_2}{t_3}.
\]
The new monomial $\frac{t_2}{t_3}$
is associated with $(1,0,0,1,0,1,2,1,1)$.

\vspace{2mm}

\hspace{-7mm}
\underline{$\overline{DG}_2$} Let $d_j$ be the exponent of $t_j$ in the unique sink 
$\frac{t_2}{t_3}$ of $\overline{DG}_1$ so that $(d_1,d_2,\cdots,d_9)=(0,1,-1,0,0,0,0,0,0)$.
The associated integer vector to this sink is
$(b_1,\cdots,b_9)=(1,0,0,1,0,1,2,1,1)$.
Since $d_2>0$, $b_{2^+}=b_{6}>0$ and $d_{2^+}=d_6=0<d_2$, we see that $\frac{t_2}{t_3}$ and $j=2$ satisfy (a). Hence,
\[
\overline{DG}_2\ :\
t_1\rightarrow \frac{t_2}{t_3}\rightarrow \frac{t_2}{t_3}\cdot A_{2}^{-1}=\frac{t_3t_5^2}{t_6}.
\]
To $\frac{t_3t_5^2}{t_6}$, we associate
\[
(1,1,0,1,0,0,2,1,1).
\]
\underline{$\overline{DG}_3$} Similarly, 
in $\overline{DG}_2$, the sink
$\frac{t_3t_5^2}{t_6}$ and $j=5$ satisfy (a) so that
\[
\overline{DG}_3\ :\
t_1\rightarrow \frac{t_2}{t_3}\rightarrow \frac{t_3t_5^2}{t_6}\rightarrow \frac{t_3t_5^2}{t_6}\cdot A_{5}^{-1}=\frac{t_3t_5}{t_7}.
\]
To the new monomial $\frac{t_3t_5}{t_7}$, we associate $(1,1,0,1,1,0,1,1,1)$.
\vspace{2mm}

\hspace{-7mm}
\underline{$\overline{DG}_4$} $\frac{t_3t_5}{t_7}$ is the unique sink in $\overline{DG}_3$ and the associated vector
and exponents $d_j$ of $t_j$ are as follows:
\[
(b_1,b_2\cdots,b_9)=(1,1,0,1,1,0,1,1,1)
\]
\[
(d_1,d_2,\cdots,d_9)=(0,0,1,0,1,0,-1,0,0).
\]
One can confirm that
\begin{itemize}
\item $d_3>0$, $b_{3^+}=b_5>0$ and $d_{3^+}=d_3$, $d_{3^{2+}}=d_{7}=-1$, $b_{3^{2+}}=b_{7}=1$. Thus,
$\frac{t_3t_5}{t_7}$ and $j=3$ satisfy (b) with $p=2$.
\item $d_5>0$, $b_{5^+}=b_7>0$ and $d_{5^+}=d_7=-1<d_5$. Thus, $\frac{t_3t_5}{t_7}$ and $j=5$ satisfy (a).
\end{itemize}
Therefore, we get the graph $\overline{DG}_4$:
\[
\begin{xy}
(-30,20) *{\overline{DG}_4 :}="0",
(-20,20) *{t_1}="1",
(-10,20) *{\frac{t_2}{t_3}}="2",
(0,20) *{\frac{t_3t_5^2}{t_6}}="3",
(10,20) *{\frac{t_3t_5}{t_7}}="4",
(10,10) *{\qquad \qquad \frac{t_3t_5}{t_7}\cdot A_{5}^{-1}=\frac{t_3t_6}{t_7^2}}="4-1",
(30,20) *{\frac{t_3t_5}{t_7}\cdot A_{3}^{-1}=\frac{t_4}{t_7}}="5",
\ar@{->} "1";"2"^{}
\ar@{->} "2";"3"^{}
\ar@{->} "3";"4"^{}
\ar@{->} "4";"5"^{}
\ar@{->} "4";"4-1"^{}
\end{xy}
\]
To $\frac{t_4}{t_7}$ and $\frac{t_3t_6}{t_7^2}$, we associate $(1,1,1,1,0,0,1,1,1)$ and $(1,1,0,1,2,0,0,1,1)$, respectively.
Repeating this argument, we get
\[
\begin{xy}
(-35,20) *{\overline{DG}_7 :}="0",
(-20,20) *{t_1}="1",
(-10,20) *{\frac{t_2}{t_3}}="2",
(0,20) *{\frac{t_3t_5^2}{t_6}}="3",
(10,20) *{\frac{t_3t_5}{t_7}}="4",
(10,10) *{\frac{t_3t_6}{t_7^2}}="4-1",
(10,0) *{\frac{t_3}{t_8}}="4-2",
(20,10) *{\frac{t_4t_6}{t_5t_7^2}}="5-1",
(20,0) *{\frac{t_4}{t_5t_8}}="5-2",
(20,20) *{\frac{t_4}{t_7}}="5",
(30,20) *{\frac{t_5}{t_9}}="6",
(40,20) *{\frac{t_6}{t_7t_9}}="7",
(50,20) *{\frac{t_7}{t_8t_9}}="8",
\ar@{->} "1";"2"^{1}
\ar@{->} "2";"3"^{2\ }
\ar@{->} "3";"4"^{5\ }
\ar@{->} "4";"5"^3
\ar@{->} "4";"4-1"^{5}
\ar@{->} "5";"6"^{4}
\ar@{->} "4-1";"5-1"^{3}
\ar@{->} "4-1";"4-2"^{6}
\ar@{->} "5-1";"5-2"^{6}
\ar@{->} "5-1";"7"_{4}
\ar@{->} "5-2";"8"_{4}
\ar@{->} "4-2";"5-2"^{3}
\ar@{->} "6";"7"^{5}
\ar@{->} "7";"8"^{6}
\end{xy}
\]
Here, the associated vectors to other monomials are as follows:
\[
\frac{t_5}{t_9} : (1,1,1,2,0,0,1,1,0),\quad
\frac{t_6}{t_7t_9} : (1,1,1,2,1,0,0,1,0),\quad
\frac{t_7}{t_8t_9} : (1,1,1,2,1,1,0,0,0),
\]
\[
\frac{t_3}{t_8} : (1,1,0,1,2,1,0,0,1),\quad
\frac{t_4t_6}{t_5t_7^2} : (1,1,1,1,1,0,0,1,1),\quad
\frac{t_4}{t_5t_8} : (1,1,1,1,1,1,0,0,1)
\]
and $M_1\overset{j}{\rightarrow} M_2$ implies 
$M_2=M_1\cdot A_j^{-1}$ for two monomials
$M_1$, $M_2$.
The unique sink $\frac{t_7}{t_8t_9}$ of $\overline{DG}_7$ and any $j\in[1,9]$ do not satisfy the condition of (2) so that
the algorithm stops and $\overline{DG}=\overline{DG}_7$.
Therefore,
$\Delta_{w_0\Lm_2,s_2\Lm_2}\circ \theta^-_{\mathbf{i}}(t_1,\cdots,t_9)$ is a
linear combination of monomials in $\overline{DG}$ with positive integer coefficients.
In fact, it holds
\begin{eqnarray*}\hspace{-2mm}
\Delta_{w_0\Lm_2,s_2\Lm_2}\circ \theta^-_{\mathbf{i}}(t_1,\cdots,t_9)
&=&t_1+\frac{t_2}{t_3}+\frac{t_3t_5^2}{t_6}+2\frac{t_3t_5}{t_7}+\frac{t_4}{t_7}+\frac{t_5}{t_9}+\frac{t_6}{t_7t_9}\\
& & +\frac{t_7}{t_8t_9}+
\frac{t_3t_6}{t_7^2}+\frac{t_3}{t_8}+\frac{t_4t_6}{t_5t_7^2}+\frac{t_4}{t_5t_8}.
\end{eqnarray*}

Next, let us compute the graph $\overline{DG}$ for $i=1$ and $i=3$.

\vspace{2mm}

\hspace{-7mm}
\underline{The graph $\overline{DG}$ for $i=1$} The graph $\overline{DG}_0$ has the unique
monomial $t_9$ and no arrow. Since $9^+=10>9$, we see that the sink $t_9$ and any $j\in[1,9]$ do not satisfy the condition of (2)
without compute $(b_1^0,b_2^0,\cdots,b_9^0)$.
Thus, it holds $\overline{DG}=\overline{DG}_0$
and
$\Delta_{w_0\Lm_1,s_1\Lm_1}\circ \theta^-_{\mathbf{i}}(t_1,\cdots,t_9)$ has only one monomial $t_9$.
\vspace{2mm}

\hspace{-7mm}
\underline{The graph $\overline{DG}$ for $i=3$} 
By a similar way to the case $i=1$,
the graph $\overline{DG}_0$ has the unique
monomial $t_8$ and no arrow and $\overline{DG}=\overline{DG}_0$.
Hence, $\Delta_{w_0\Lm_3,s_3\Lm_3}\circ \theta^-_{\mathbf{i}}(t_1,\cdots,t_9)$ has only one monomial $t_8$.

Therefore, considering the tropicalization and identifying $X_*((\mathbb{C}^{\times})^{9})=\mathbb{Z}^9$, one obtains an explicit form $\mathbb{B}_{\theta^-_{\mathbf{i}},\Phi^{h}_{\rm BK}}$
in (\ref{boldb}) :
\[
\mathbb{B}_{\theta^-_{\mathbf{i}},\Phi^{h}_{\rm BK}}
=
\left\{
(z_1,z_2,\cdots,z_9)\in\mathbb{Z}^9 \left| 
\begin{array}{l}
z_1\geq0,\ z_2-z_3\geq0,\ z_3+2z_5-z_6\geq0,\\
z_3+z_5-z_7\geq0,\ z_4-z_7\geq0,\\
z_5-z_9\geq0,\ z_6-z_7-z_9\geq0,\\
 z_7-z_8-z_9\geq0,
z_3+z_6-2z_7\geq0,\\
z_3-z_8\geq0,
z_4+z_6-z_5-2z_7\geq0,\\
z_4-z_5-z_8\geq0,\ z_8\geq0,\ z_9\geq0
\end{array}\right.
\right\}.
\]

\end{ex}

\begin{ex}\label{ex2}
Let $G$ be of type ${\rm D}_4$, $\mathbf{i}=(2,1,3,2,4,2,3,2,1,2,3,4)$ and $i=2$.
By $s_4s_3s_2s_1s_2s_3s_2s_4s_2s_3s_1\alpha_2=\alpha_2$,
we see that $\overline{DG}_0$ is the graph which has only one vertex $t_1$ and no arrow.
Following Lemma \ref{lem-3}, the monomial $t_1$ is associated with
\[
(b_1^0,b_2^0,\cdots,b_{12}^0)=(0,0,0,1,1,0,1,1,2,1,1,1).
\]
Following the algorithm, we obtain the following graph $\overline{DG}$. Here, $M_1\overset{j}{\rightarrow} M_2$ implies 
$M_2=M_1\cdot A_j^{-1}$ for two monomials
$M_1$, $M_2$:

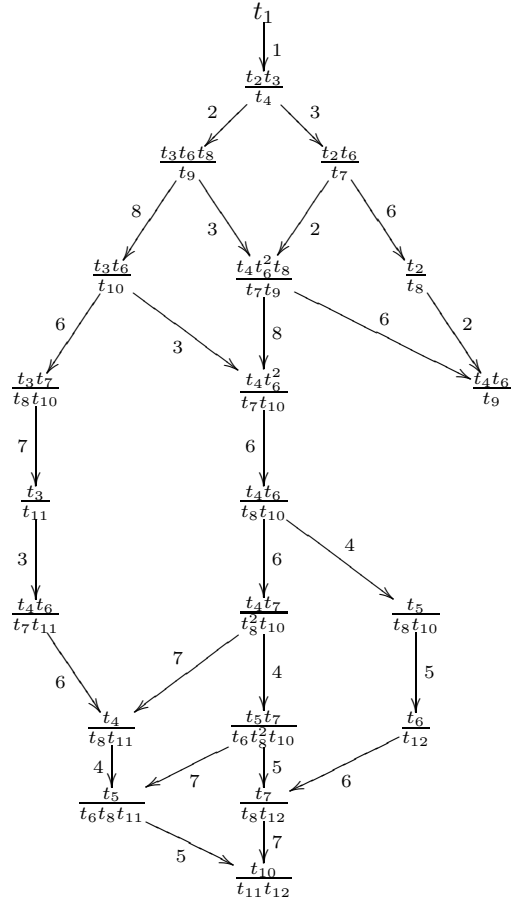
\[
\begin{xy}
(0,0)*{t_1}="1",
(0,-10) *{\frac{t_2t_3}{t_4}}="234",
(10,-20) *{\frac{t_2t_6}{t_7}}="267",
(-10,-20) *{\frac{t_3t_6t_8}{t_9}}="3689",
(20,-35) *{\frac{t_2}{t_8}}="28",
(0,-35) *{\frac{t_4t_6^2t_8}{t_7t_9}}="466879",
(-20,-35) *{\frac{t_3t_6}{t_{10}}}="3610",
(30,-50) *{\frac{t_4t_6}{t_9}}="469",
(0,-50) *{\frac{t_4t_6^2}{t_7t_{10}}}="466710",
(-30,-50) *{\frac{t_3t_7}{t_8t_{10}}}="37810",
(0,-65) *{\frac{t_4t_6}{t_8t_{10}}}="46810",
(-30,-65) *{\frac{t_3}{t_{11}}}="311",
(20,-80) *{\frac{t_5}{t_8t_{10}}}="5810",
(0,-80) *{\frac{t_4t_7}{t_8^2t_{10}}}="478810",
(-30,-80) *{\frac{t_4t_6}{t_7t_{11}}}="46711",
(20,-95) *{\frac{t_6}{t_{12}}}="612",
(0,-95) *{\frac{t_5t_7}{t_6t_8^2t_{10}}}="5768810",
(-20,-95) *{\frac{t_4}{t_8t_{11}}}="4811",
(0,-115) *{\frac{t_{10}}{t_{11}t_{12}}}="101112",
(0,-105) *{\frac{t_7}{t_8t_{12}}}="7812",
(-20,-105) *{\frac{t_5}{t_6t_8t_{11}}}="56811",
\ar@{->} "1";"234"^{1}
\ar@{->} "234";"267"^{3}
\ar@{->} "234";"3689"_{2}
\ar@{->} "267";"28"^{6}
\ar@{->} "267";"466879"^{2}
\ar@{->} "3689";"466879"_{3}
\ar@{->} "3689";"3610"_{8}
\ar@{->} "28";"469"^{2}
\ar@{->} "466879";"469"^{6}
\ar@{->} "466879";"466710"^{8}
\ar@{->} "3610";"466710"_{3}
\ar@{->} "3610";"37810"_{6}
\ar@{->} "466710";"46810"_{6}
\ar@{->} "37810";"311"_{7}
\ar@{->} "46810";"5810"^{4}
\ar@{->} "46810";"478810"^{6}
\ar@{->} "311";"46711"_{3}
\ar@{->} "5810";"612"^{5}
\ar@{->} "478810";"5768810"^{4}
\ar@{->} "478810";"4811"_{7}
\ar@{->} "46711";"4811"_{6}
\ar@{->} "612";"7812"^{6}
\ar@{->} "5768810";"7812"^{5}
\ar@{->} "5768810";"56811"^{7}
\ar@{->} "4811";"56811"_{4}
\ar@{->} "7812";"101112"^{7}
\ar@{->} "56811";"101112"_{5}
\end{xy}
\]
Therefore,
$\Delta_{w_0\Lm_2,s_2\Lm_2}\circ \theta^-_{\mathbf{i}}(t_1,\cdots,t_{12})$ is a
linear combination of monomials in $\overline{DG}$ with positive integer coefficients.
The graphs for $i=1,3,4$ are as follows:
\[
\overline{DG}
(i=1) : t_8\overset{8}{\rightarrow} \frac{t_9}{t_{10}},
\qquad
\overline{DG}(i=3) : t_{11},\qquad \overline{DG}(i=4) : t_{12}.
\]
Therefore, identifying $X_*((\mathbb{C}^{\times})^{12})=\mathbb{Z}^{12}$,
\[
\mathbb{B}_{\theta^-_{\mathbf{i}},\Phi^{h}_{\rm BK}}
=
\left\{
(z_1,z_2,\cdots,z_{12})\in\mathbb{Z}^{12} \left| 
\begin{array}{l}
z_1\geq0,\ z_2+z_3-z_4\geq0,\ z_3+z_6+z_8-z_9\geq0,\\
z_2+z_6-z_7\geq0,\ z_3+z_6-z_{10}\geq0,\\
z_4+2z_6+z_8-z_7-z_9\geq0,\ z_2-z_8\geq0,\\
z_3+ z_7-z_8-z_{10}\geq0,
z_4+2z_6-z_7-z_{10}\geq0,\\
z_4+z_6-z_9\geq0,
z_3-z_{11}\geq0,\\
z_4+z_6-z_8-z_{10}\geq0,\ z_4+z_6-z_7-z_{11}\geq0,\\
z_4+z_7-2z_8-z_{10}\geq0,\ z_5-z_8-z_{10}\geq0,\\
z_4-z_8-z_{11}\geq0,\ z_5+z_7-z_6-2z_8-z_{10}\geq0,\\
z_6-z_{12}\geq0,
z_5-z_6-z_8-z_{11}\geq0,\\ z_7-z_8-z_{12}\geq0,\ z_{10}-z_{11}-z_{12}\geq0,\\
z_8\ge0,\ z_9-z_{10}\geq0,\ 
z_{11}\ge0,\ z_{12}\geq0
\end{array} \right.
\right\}.
\]

\end{ex}

\section{Proof}

\subsection{Properties of $\mathbf{b}$-integers}

Let $\mathbf{i}=(i_1,\cdots,i_N)$ be a fixed reduced word of $w_0$.
In this subsection, the index $i\in I$ is arbitrary, that is, we do not assume the condition (\ref{mini-pro}).

\begin{prop}\label{lem-2}
Let $M=\prod_{l=1}^{N}t_l^{d_l}$ be a monomial in $\Delta_{w_0\Lm_i,s_i\Lm_i}\circ \theta^-_{\mathbf{i}}(t_1,\cdots,t_N)$
and 
$\{b_t\}_{t\in[1,N]}$ be integers
in Definition \ref{bint} determined from $M$.
If $b_{j^+}=0$ with some $j\in[1,N]$ such that $j^+\leq N$
then $M\cdot A_{j}^{-1}$ is not a monomial in
$\Delta_{w_0\Lm_i,s_i\Lm_i}\circ \theta^-_{\mathbf{i}}(t_1,\cdots,t_N)$.
\end{prop}

\nd
{\it Proof.}
We assume $M\cdot A_{j}^{-1}$ is a monomial in
$\Delta_{w_0\Lm_i,s_i\Lm_i}\circ \theta^-_{\mathbf{i}}(t_1,\cdots,t_N)$
and deduce a contradiction from this assumption.
Let $\pi=(\gamma_0,\gamma_1,\cdots,\gamma_N)$ be the $\mathbf{i}$-trail corresponding to $M$
with integers $c_l$ $(l=1,2,\cdots,N)$ such that
$\gamma_{l-1}-\gamma_l=c_l\alpha_{i_l}$.
We put $M\cdot A_{j}^{-1}=\prod_{l=1}^{N}t_l^{d_l'}$.
It is easy to see
\[
d_l'=d_l\quad (j^+<l\leq N),\quad d_{j^+}'=d_{j^+}-1.
\]
Let $\pi'=(\gamma_0',\gamma_1',\cdots,\gamma_N')$ be the $\mathbf{i}$-trail corresponding to $M\cdot A_{j}^{-1}$
with non-negative integers $c_l'$ ($l=1,2,\cdots,N$) such that $\gamma_{l-1}'-\gamma_l'=c_l'\alpha_{i_l}$.
For $l$ ($j^+<l\leq N$), let us show
\begin{equation}\label{lempr-1}
\text{if }\gamma_{l}=\gamma_{l}' \text{ then } \gamma_{l-1}=\gamma_{l-1}'. 
\end{equation}
By $\gamma_{l}=\gamma_{l}'$ and
\[
d_l=\frac{\gamma_{l-1}+\gamma_l}{2}(h_{i_l})=d_l'=\frac{\gamma_{l-1}'+\gamma_l'}{2}(h_{i_l}),
\]
it holds $\gamma_{l-1}(h_{i_l})=\gamma_{l-1}'(h_{i_l})$ so that
$c_l=\frac{\gamma_{l-1}-\gamma_l}{2}(h_{i_l})=\frac{\gamma_{l-1}'-\gamma_l'}{2}(h_{i_l})=c_l'$. Hence,
it follows $\gamma_{l-1}=\gamma_l+c_l\alpha_{i_l}=\gamma_l'+c_l'\alpha_{i_l}=\gamma_{l-1}'$ and we get
(\ref{lempr-1}). Since $\gamma_N=\gamma_N'(=-s_i\Lambda_i)$, it follows from (\ref{lempr-1}) that
\begin{equation}\label{lempr-2}
\gamma_l=\gamma_l'\ (j^+\leq l\leq N).
\end{equation}
Therefore, we see that
\[
d_{j^+}'=\frac{\gamma_{j^+-1}'+\gamma_{j^+}'}{2}(h_{i_j})
=\frac{\gamma_{j^+-1}'+\gamma_{j^+}}{2}(h_{i_j}).
\]
Combining this equality, $d_{j^+}=\frac{\gamma_{j^+-1}+\gamma_{j^+}}{2}(h_{i_j})$ and $d_{j^+}'=d_{j^+}-1$, one gets
\begin{equation}\label{lempr-3}
\gamma_{j^+-1}'(h_{i_j})=\gamma_{j^+-1}(h_{i_j})-2.
\end{equation}
The assumption $b_{j^+}=0$ and Lemma \ref{btct} induce
\begin{equation}\label{lempr-4}
\frac{\gamma_{j^+-1}-\gamma_{j^+}}{2}(h_{i_j})=
c_{j^+}=b_{j^+}=0.
\end{equation}
Therefore, using (\ref{lempr-2}), (\ref{lempr-3}) and (\ref{lempr-4}),
\[
c_{j^+}'=\frac{\gamma_{j^+-1}'-\gamma_{j^+}'}{2}(h_{i_j})=
\frac{\gamma_{j^+-1}'-\gamma_{j^+}}{2}(h_{i_j})=
\frac{\gamma_{j^+-1}-\gamma_{j^+}}{2}(h_{i_j})-1=-1,
\]
which contradicts the non-negativity of $c_{j^+}'$. \qed

By this proposition, if $b_{j^+}=0$ then $M\cdot A_j^{-1}$ is not a monomial in
$\Delta_{w_0\Lm_i,s_i\Lm_i}\circ \theta^-_{\mathbf{i}}(t_1,\cdots,t_N)$.
In several lemmas of next subsection, we consider the case $b_{j^+}>0$.
Before that, we see the following lemma:

\begin{lem}\label{lem-4}
We suppose that $M$ and $M':=M\cdot A_j^{-1}$ are monomials in 
$\Delta_{w_0\Lm_i,s_i\Lm_i}\circ \theta^-_{\mathbf{i}}(t_1,\cdots,t_N)$ with some $j\in[1,N]$ such that $j^+\leq N$.
Let $\{b_t\}_{t\in[1,N]}$ and $\{b_t'\}_{t\in[1,N]}$ be integers
in Definition \ref{bint} determined from $M$ and $M'$, respectively.
Then we have
\[
b_l'=b_{l}\quad \text{for}\ l\in[1,N]\setminus \{j,j^+\},
\]
\[
b_j'=b_{j}+1,\quad b_{j^+}'=b_{j^+}-1.
\]
\end{lem}

\nd
{\it Proof.} 
Let $\pi=(\gamma_0,\gamma_1,\cdots,\gamma_N)$ and
$\pi'=(\gamma_0',\gamma_1',\cdots,\gamma_N')$ be $\mathbf{i}$-trails 
corresponding to $M$ and $M'$ with integers $\{c_l\}_{l\in[1,N]}$, $\{c_l'\}_{l\in[1,N]}$
such that $\gamma_{l-1}-\gamma_l=c_l\alpha_{i_l}$ and $\gamma_{l-1}'-\gamma_l'=c_l'\alpha_{i_l}$, respectively.
Note that considering Proposition \ref{lem-2}, it holds
$b_{j^+}\neq0$ so that Lemma \ref{btct} induces $c_{j^+}=b_{j^+}>0$. If we take a pre-$\mathbf{i}$-trail
$\pi''=(\gamma_0'',\gamma_1'',\cdots,\gamma_N'')$ from $-w_0\Lambda_i$ to $-s_i\Lambda_i$
with integers $\{c''_l\}_{l\in[1,N]}$ satisfying $\gamma_{l-1}''-\gamma_l''=c_l''\alpha_{i_l}$ 
such that
\[
c_l''=c_{l}\quad \text{for}\ l\in[1,N]\setminus \{j,j^+\},
\]
\[
c_j''=c_{j}+1,\quad c_{j^+}''=c_{j^+}-1
\]
then Proposition \ref{trail-prop} yields $\prod_{l=1}^N t_l^{d_l(\pi'')}=M\cdot A_{j}^{-1}$.
Taking Lemma \ref{i-trail-lem} into account, we see that $\pi''=\pi'$ so that
\[
c_l'=c_{l}''\quad \text{for\ all}\ l\in[1,N].
\]
By Lemma \ref{btct}, it holds $c_l=b_l$, $c_l'=b_l'$, which yields our claim.\qed

\subsection{Proof of Theorem \ref{thm2}}

We consider the same setting as in subsection \ref{algo}. In particular,
the index $i\in I$ satisfies (\ref{mini-pro}).

\begin{lem}\label{lem-a1}
\begin{enumerate}
\item[(1)]
For $j\in[0,N]$ and $c_{j+1},\cdots,c_N\in\mathbb{Z}_{\geq0}$, we suppose that
$e^{c_{j+1}}_{i_{j+1}}\cdots e^{c_N}_{i_N}v_{-s_i\Lambda_i}\neq 0$.
For $t\in I$,
the value of
\[
\langle h_t,\ {\rm wt}(e^{c_{j+1}}_{i_{j+1}}\cdots e^{c_N}_{i_N}v_{-s_i\Lambda_i}) \rangle
\]
is in $\{2,1,0,-1,-2\}$.
\item[(2)] Let
$\pi=(\gamma_0,\gamma_1,\cdots,\gamma_N)$ be an $\mathbf{i}$-trail
from $-w_0\Lambda_i$ to $-s_i\Lambda_i$. For $j\in[1,N]$ and $t\in I$, the values of
\[
\gamma_{j-1}(h_t),\ \gamma_{j}(h_t) \text{ and }
d_j(\pi)=\frac{\gamma_{j-1}+\gamma_j}{2}(h_{i_j})
\]
are in $\{2,1,0,-1,-2\}$.
\end{enumerate}
\end{lem}

\nd
{\it Proof.}

\nd
(1) Our claim
is an easy consequence of (\ref{mini-pro}).

\vspace{2mm}

\nd
(2) Let $\{c_l\}$ be the integers for $\pi$ as in Definition \ref{pretrail}.
Since
\[
\gamma_j
=\gamma_N+\sum^{N-1}_{t=j}(\gamma_t-\gamma_{t+1}) 
= -s_i\Lambda_i + \sum^{N-1}_{t=j} c_{t+1}\alpha_{i_{t+1}} 
= {\rm wt}(e^{c_{j+1}}_{i_{j+1}}\cdots e^{c_N}_{i_N}v_{-s_i\Lambda_i}),
\]
\[
\gamma_{j-1}
= {\rm wt}(e^{c_{j}}_{i_{j}} e^{c_{j+1}}_{i_{j+1}}\cdots e^{c_N}_{i_N}v_{-s_i\Lambda_i}),
\]
our statement (2) follows by (1). \qed

In the following lemma and Lemma \ref{lem-a3}, we consider the condition $d_j>0$, $b_{j^+}>0$ and $d_{j^+}<d_j$,
which is same as (a) of (2) in Theorem \ref{thm2}.

\begin{lem}\label{lem-a2}
Let $M=\prod_{l=1}^{N}t_l^{d_l}$ be a monomial in $\Delta_{w_0\Lm_i,s_i\Lm_i}\circ \theta^-_{\mathbf{i}}(t_1,\cdots,t_N)$
and $\{b_l\}_{l\in[1,N]}$ be integers in Definition \ref{bint} determined from $M$.
We suppose that $d_j=2$ and $b_{j^+}>0$ for some $j\in[1,N]$ such that $j^+\leq N$.
Then $M\cdot A_{j}^{-1}$ is also a monomial in
$\Delta_{w_0\Lm_i,s_i\Lm_i}\circ \theta^-_{\mathbf{i}}(t_1,\cdots,t_N)$ and $d_{j^+}\leq1$.
\end{lem}

\nd
{\it Proof.}
Let
$\pi=(\gamma_0,\gamma_1,\cdots,\gamma_N)$ be the $\mathbf{i}$-trail corresponding to $M$
with integers $c_l$ $(l=1,2,\cdots,N)$ such that
$\gamma_{l-1}-\gamma_l=c_l\alpha_{i_l}$. By Lemma \ref{btct}, it holds $c_l=b_l$ for $l\in[1,N]$.
We obtain
\[
d_j=d_j(\pi)=\frac{\gamma_{j-1}+\gamma_j}{2}(h_{i_j})
\]
by Theorem \ref{trail-thm}. Thus, the assumption $d_j=2$ and Lemma \ref{lem-a1} (2) yield
\[
\gamma_{j-1}(h_{i_j})=\gamma_j(h_{i_j})=2
\]
so that 
\[
c_j=\frac{\gamma_{j-1}-\gamma_j}{2}(h_{i_j})=0.
\]
It holds
\begin{eqnarray}
\gamma_j(h_{i_j})&=&
\gamma_N(h_{i_j})+\sum^{N-1}_{l=j}(\gamma_l-\gamma_{l+1}) (h_{i_j})\nonumber\\
&=&
\gamma_N(h_{i_j})+\sum^{N-1}_{l=j^+-1}(\gamma_l-\gamma_{l+1}) (h_{i_j})
+\sum^{j^+-2}_{l=j}(\gamma_l-\gamma_{l+1}) (h_{i_j})
\nonumber\\
&=&
\gamma_{j^+-1}(h_{i_j})+\langle h_{i_j}, \sum^{j^+-1}_{l=j+1} c_l\alpha_{i_l}\rangle.\qquad \quad\label{a2-pr1}
\end{eqnarray}
Note that combining $\langle h_{i_j}, \sum^{j^+-1}_{l=j+1} c_l\alpha_{i_l}\rangle\leq 0$, $\gamma_j(h_{i_j})=2$
with
Lemma \ref{lem-a1} (2), it holds $\langle h_{i_j}, \sum^{j^+-1}_{l=j+1} c_l\alpha_{i_l}\rangle=0$ and
\begin{equation}\label{pr1-1}
\gamma_{j^+-1}(h_{i_j})=2.
\end{equation}
Considering $\langle h_{i_j}, \sum^{j^+-1}_{l=j+1} c_l\alpha_{i_l}\rangle=0$,
we see that for $l\in[j+1,j^+-1]$ such that $a_{i_j,i_l}<0$, it follows
\begin{equation}\label{pr1-2}
c_l=0.
\end{equation}
Using the assumption $b_{j^+}=c_{j^+}>0$, 
there is a pre-$\mathbf{i}$-trail
$\pi'=(\gamma_0',\gamma_1',\cdots,\gamma_N')$
from $-w_0\Lambda_i$ to $-s_i\Lambda_i$ with integers $c_l'$ $(l=1,2,\cdots,N)$ such that
$\gamma_{l-1}'-\gamma_l'=c_l'\alpha_{i_l}$ and
\[
c_l'=c_{l}\quad \text{for}\ l\in[1,N]\setminus \{j,j^+\},
\]
\[
c_j'=c_{j}+1,\quad c_{j^+}'=c_{j^+}-1.
\]
The definition of $\mathbf{i}$-trails implies $e_{i_1}^{c_1}e_{i_2}^{c_2}\cdots e_{i_N}^{c_N}v_{-s_i\Lambda_i}\neq0$.
Considering (\ref{pr1-2}), we see
$e_{i_1}^{c_1'}e_{i_2}^{c_2'}\cdots e_{i_N}^{c_N'}v_{-s_i\Lambda_i}=e_{i_1}^{c_1}e_{i_2}^{c_2}\cdots e_{i_N}^{c_N}v_{-s_i\Lambda_i}\neq0$ so that
$\pi'$ is an $\mathbf{i}$-trail.
Let $M'$ be the monomial
in $\Delta_{w_0\Lm_i,s_i\Lm_i}\circ \theta^-_{\mathbf{i}}(t_1,\cdots,t_N)$ corresponding to $\pi'$.
By Proposition \ref{trail-prop},
we have
\[
M'=MA_{j}^{-1}
\]
and $MA_{j}^{-1}$ is a monomial in $\Delta_{w_0\Lm_i,s_i\Lm_i}\circ \theta^-_{\mathbf{i}}(t_1,\cdots,t_N)$. 
Note that by $c_{j^+}=\frac{\gamma_{j^+-1}-\gamma_{j^+}}{2}(h_{i_j})>0$ and (\ref{pr1-1}), it holds
$\gamma_{j^+}(h_{i_j})\leq 0$ so that
\[
d_{j^+}=\frac{\gamma_{j^+-1}+\gamma_{j^+}}{2}(h_{i_j})\leq1.
\]
\qed

\begin{lem}\label{lem-a3}
Let $M=\prod_{l=1}^{N}t_l^{d_l}$ be a monomial in $\Delta_{w_0\Lm_i,s_i\Lm_i}\circ \theta^-_{\mathbf{i}}(t_1,\cdots,t_N)$
and $\{b_l\}_{l\in[1,N]}$ be integers in Definition \ref{bint} determined from $M$.
We suppose that $d_j=1$, $d_{j^+}<1$ and $b_{j^+}>0$. Then $M\cdot A_{j}^{-1}$ is also a monomial in
$\Delta_{w_0\Lm_i,s_i\Lm_i}\circ \theta^-_{\mathbf{i}}(t_1,\cdots,t_N)$.
\end{lem}

\nd
{\it Proof.}
Let
$\pi=(\gamma_0,\gamma_1,\cdots,\gamma_N)$ be the $\mathbf{i}$-trail corresponding to $M$
with integers $c_l$ $(l=1,2,\cdots,N)$ such that
$\gamma_{l-1}-\gamma_l=c_l\alpha_{i_l}$.
We see that
\[
\gamma_{j-1}(h_{i_j})=\gamma_{j}(h_{i_j})+c_j\alpha_{i_j}(h_{i_j})
=\gamma_{j}(h_{i_j})+2c_j
\geq \gamma_{j}(h_{i_j}).
\]
By $d_j=\frac{\gamma_{j-1}+\gamma_{j}}{2}(h_{i_j})=1$ and Lemma \ref{lem-a1} (2),
it holds either $(\gamma_{j-1}(h_{i_j}), \gamma_{j}(h_{i_j}))=(1,1)$ or 
$(\gamma_{j-1}(h_{i_j}), \gamma_{j}(h_{i_j}))=(2,0)$.
Just as in (\ref{a2-pr1}), it holds
\begin{equation}\label{a3-pr1}
\gamma_j(h_{i_j})
=\gamma_{j^+-1}(h_{i_j})+\langle h_{i_j}, \sum^{j^+-1}_{l=j+1} c_l\alpha_{i_l}\rangle.
\end{equation}
\nd
Let us prove our claim by division into cases.

\vspace{2mm}

\nd
\underline{Case 1. $(\gamma_{j-1}(h_{i_j}), \gamma_{j}(h_{i_j}))=(1,1)$}

\vspace{2mm}

In this case, we obtain 
\begin{equation}\label{a3-pr2}
c_j=\frac{\gamma_{j-1}-\gamma_{j}}{2}(h_{i_j})=0.
\end{equation}
Taking (\ref{a3-pr1}) and $\langle h_{i_j}, \sum^{j^+-1}_{l=j+1} c_l\alpha_{i_l}\rangle\leq0$
into account, the value of $\gamma_{j^+-1}(h_{i_j})$ is $1$ or $2$.

\vspace{2mm}

\nd
\underline{Case 1-1. $\gamma_{j^+-1}(h_{i_j})=1$}

\vspace{2mm}

In this case, the relation (\ref{a3-pr1}) yields $\langle h_{i_j}, \sum^{j^+-1}_{l=j+1} c_l\alpha_{i_l}\rangle=0$ so that
$c_l=0$ for $l\in[j+1,j^+-1]$ such that $a_{i_j,i_l}<0$. Since $b_{j^+}=c_{j^+}>0$, 
we can take a pre-$\mathbf{i}$-trail
$\pi'=(\gamma_0',\gamma_1',\cdots,\gamma_N')$
from $-w_0\Lambda_i$ to $-s_i\Lambda_i$
just as in the 
proof of Lemma \ref{lem-a2}
and similarly
verify $M\cdot A_{j}^{-1}$ is a monomial in
$\Delta_{w_0\Lm_i,s_i\Lm_i}\circ \theta^-_{\mathbf{i}}(t_1,\cdots,t_N)$.

\vspace{2mm}

\nd
\underline{Case 1-2. $\gamma_{j^+-1}(h_{i_j})=2$}

\vspace{2mm}

In this case, by (\ref{a3-pr1}), there uniquely exists $l\in[j+1,j^+-1]$ such that 
\begin{equation}\label{a3-pr3}
a_{i_j,i_l}=-1,\quad c_l=1.
\end{equation}
For $l'\in[j+1,j^+-1]\setminus\{l\}$ such that $a_{i_j,i_{l'}}<0$, it holds
\begin{equation}\label{a3-pr4}
c_{l'}=0.
\end{equation}
By $\gamma_{j^+-1}(h_{i_j})=2$ and $d_{j^+}=\frac{\gamma_{j^+-1}+\gamma_{j^+}}{2}(h_{i_j})<1$,
we have $\gamma_{j^+}(h_{i_j})=-2$ so that
\begin{equation}\label{a3-pr5}
c_{j^+}=\frac{\gamma_{j^+-1}-\gamma_{j^+}}{2}(h_{i_j})=2.
\end{equation}
Since $\pi$ is an $\mathbf{i}$-trail and it holds (\ref{a3-pr2}), (\ref{a3-pr3}), (\ref{a3-pr5}), we obtain
\begin{equation}\label{a3-pr55}
e_{i_1}^{c_1} \cdots 
e_{i_j}^{c_j}\cdots e_{i_l}^{c_l}\cdots e_{i_{j^+}}^{c_{j^+}}\cdots 
e_{i_N}^{c_N}v_{-s_i\Lambda_i}
=
e_{i_1}^{c_1} \cdots 
e_{i_j}^{0}\cdots e_{i_l}^{1}\cdots e_{i_{j^+}}^{2}\cdots 
e_{i_N}^{c_N}v_{-s_i\Lambda_i}
\neq0.
\end{equation}
If $e_{i_{l+1}}^{c_{l+1}}\cdots e_{i_{j^+-1}}^{c_{j^+-1}} e_{i_{j^+}}^{0}e_{i_{j^++1}}^{c_{j^++1}}
\cdots 
e_{i_N}^{c_N}v_{-s_i\Lambda_i}\neq0$ then
\[
\langle h_{i_j}, {\rm wt}(e_{i_{l+1}}^{c_{l+1}}\cdots e_{i_{j^+-1}}^{c_{j^+-1}} e_{i_{j^+}}^{0}e_{i_{j^++1}}^{c_{j^++1}}
\cdots 
e_{i_N}^{c_N}v_{-s_i\Lambda_i})
\rangle=-2
\]
by $-2=\gamma_{j^+}(h_{i_j})=\langle h_{i_j}, {\rm wt}(e_{i_{j^++1}}^{c_{j^++1}}
\cdots 
e_{i_N}^{c_N}v_{-s_i\Lambda_i})
\rangle$. Thus, the assumption (\ref{mini-pro}) and (\ref{a3-pr3}) yield
\begin{equation}\label{0-cond}
e_{i_l}^1e_{i_{l+1}}^{c_{l+1}}\cdots e_{i_{j^+-1}}^{c_{j^+-1}} e_{i_{j^+}}^{0}e_{i_{j^++1}}^{c_{j^++1}}
\cdots 
e_{i_N}^{c_N}v_{-s_i\Lambda_i}=0,
\end{equation}
otherwise $\langle h_{i_j}, e_{i_l}^1e_{i_{l+1}}^{c_{l+1}}\cdots e_{i_{j^+-1}}^{c_{j^+-1}} e_{i_{j^+}}^{0}e_{i_{j^++1}}^{c_{j^++1}}
\cdots 
e_{i_N}^{c_N}v_{-s_i\Lambda_i} \rangle=-3$.
In case of
\[
e_{i_{l+1}}^{c_{l+1}}\cdots e_{i_{j^+-1}}^{c_{j^+-1}} e_{i_{j^+}}^{0}e_{i_{j^++1}}^{c_{j^++1}}
\cdots 
e_{i_N}^{c_N}v_{-s_i\Lambda_i}=0,
\]
the equation (\ref{0-cond}) clearly holds.
Therefore,
\begin{equation}\label{a3-pr6}
e_{i_1}^{c_1} \cdots e_{i_{j-1}}^{c_{j-1}} 
e_{i_j}^{2}e_{i_{j+1}}^{c_{j+1}} 
\cdots e_{i_{l-1}}^{c_{l-1}} 
e_{i_l}^{1}e_{i_{l+1}}^{c_{l+1}} 
\cdots e_{i_{j^+-1}}^{c_{j^+-1}} 
 e_{i_{j^+}}^{0} 
e_{i_{j^++1}}^{c_{j^+ +1}} 
\cdots 
e_{i_N}^{c_N}v_{-s_i\Lambda_i}=0.
\end{equation}
Note that it holds $i_j=i_{j^+}$.
Using (\ref{a3-pr4}), (\ref{a3-pr55}), (\ref{a3-pr6})
and the Serre relation $e_{i_j}^2e_{i_l}-2e_{i_j}e_{i_l}e_{i_{j^+}}+e_{i_l}e_{i_{j^+}}^2=0$, we can verify
\begin{equation}\label{a3-pr7}
e_{i_1}^{c_1} \cdots 
e_{i_j}^{1}\cdots e_{i_l}^{1}\cdots e_{i_{j^+}}^{1}\cdots 
e_{i_N}^{c_N}v_{-s_i\Lambda_i}
\neq0.
\end{equation}
We define
$\pi'=(\gamma_0',\gamma_1',\cdots,\gamma_N')$ as the pre-$\mathbf{i}$-trail
from $-w_0\Lambda_i$ to $-s_i\Lambda_i$ with integers $c_t'$ $(t=1,2,\cdots,N)$ such that
$\gamma_{t-1}'-\gamma_t'=c_t'\alpha_{i_t}$ and
\[
c_t'=c_{t}\quad \text{for}\ t\in[1,N]\setminus \{j,j^+\},
\]
\[
c_j'=c_{j}+1=1,\quad c_{j^+}'=c_{j^+}-1=1.
\]
Then $\pi'$ is an $\mathbf{i}$-trail by (\ref{a3-pr7}).
Let $M'$ be the monomial
in $\Delta_{w_0\Lm_i,s_i\Lm_i}\circ \theta^-_{\mathbf{i}}(t_1,\cdots,t_N)$ corresponding to $\pi'$.
By Proposition \ref{trail-prop},
we have
\[
M'=MA_{j}^{-1}
\]
and $MA_{j}^{-1}$ is a monomial in $\Delta_{w_0\Lm_i,s_i\Lm_i}\circ \theta^-_{\mathbf{i}}(t_1,\cdots,t_N)$.

\vspace{2mm}

\nd
\underline{Case 2. $(\gamma_{j-1}(h_{i_j}), \gamma_{j}(h_{i_j}))=(2,0)$}

\vspace{2mm}

In this case, we obtain 
\begin{equation}\label{a3-pr8}
c_j=\frac{\gamma_{j-1}-\gamma_{j}}{2}(h_{i_j})=1.
\end{equation}
Taking (\ref{a3-pr1}) into account, the value of $\gamma_{j^+-1}(h_{i_j})$ is $0$, $1$ or $2$.

\vspace{2mm}

\nd
\underline{Case 2-1. $\gamma_{j^+-1}(h_{i_j})=0$}

\vspace{2mm}

In this case, the relation (\ref{a3-pr1}) yields
$c_l=0$ for $l\in[j+1,j^+-1]$ such that $a_{i_j,i_l}<0$. By the same argument as in the 
proof of Lemma \ref{lem-a2}, we see that $M\cdot A_{j}^{-1}$ is a monomial in
$\Delta_{w_0\Lm_i,s_i\Lm_i}\circ \theta^-_{\mathbf{i}}(t_1,\cdots,t_N)$.

\vspace{2mm}

\nd
\underline{Case 2-2. $\gamma_{j^+-1}(h_{i_j})=1$}

\vspace{2mm}

In this case, by (\ref{a3-pr1}), there uniquely exists $l\in[j+1,j^+-1]$ such that 
\[
a_{i_j,i_l}=-1,\quad c_l=1.
\]
For $l'\in[j+1,j^+-1]\setminus\{l\}$ such that $a_{i_j,i_{l'}}<0$, it holds
\[
c_{l'}=0.
\]
The assumption $b_{j^+}=c_{j^+}=\frac{\gamma_{j^+-1}-\gamma_{j^+}}{2}(h_{i_j})>0$ and Lemma \ref{lem-a1} yield
$\gamma_{j^+}(h_{i_j})=-1$ so that $c_{j^+}=1$. Note that considering
\[
-1=\gamma_{j^+}(h_{i_j})=
\langle
h_{i_j},
{\rm wt}(e^{c_{j^++1}}_{i_{j^++1}}\cdots e^{c_{N}}_{i_{N}}v_{-s_i\Lambda_i})
\rangle
\]
and (\ref{mini-pro}),
we obtain 
\[
e^{2}_{i_{j^+}}e^{c_{j^++1}}_{i_{j^++1}}\cdots e^{c_{N}}_{i_{N}}v_{-s_i\Lambda_i}=0,
\]
which yields
\begin{equation}\label{a3-pr9}
e_{i_1}^{c_1}
\cdots
e_{i_{j-1}}^{c_{j-1}}
e_{i_j}^{0}
e_{i_{j+1}}^{c_{j+1}}
\cdots e_{i_{l-1}}^{c_{l-1}}
e_{i_l}^{1}
e_{i_{l+1}}^{c_{l+1}}
\cdots e^{c_{j^+-1}}_{i_{j^+-1}}
e^{2}_{i_{j^+}}e^{c_{j^++1}}_{i_{j^++1}}\cdots e^{c_{N}}_{i_{N}}v_{-s_i\Lambda_i}=0.
\end{equation}
Since $\pi$ is an $\mathbf{i}$-trail,
\begin{equation}\label{a3-pr10}
e_{i_1}^{c_1} \cdots 
e_{i_j}^{c_j}\cdots e_{i_l}^{c_l}\cdots e_{i_{j^+}}^{c_{j^+}}\cdots 
e_{i_N}^{c_N}v_{-s_i\Lambda_i}
=
e_{i_1}^{c_1} \cdots 
e_{i_j}^{1}\cdots e_{i_l}^{1}\cdots e_{i_{j^+}}^{1}\cdots 
e_{i_N}^{c_N}v_{-s_i\Lambda_i}
\neq0.
\end{equation}
Using the Serre relation and (\ref{a3-pr9}), (\ref{a3-pr10}), it follows
\[
e_{i_1}^{c_1} \cdots 
e_{i_j}^{2}\cdots e_{i_l}^{1}\cdots e_{i_{j^+}}^{0}\cdots 
e_{i_N}^{c_N}v_{-s_i\Lambda_i}
\neq0.
\]
By a similar argument to Case 1-2, we see that
$M\cdot A_{j}^{-1}$ is a monomial in
$\Delta_{w_0\Lm_i,s_i\Lm_i}\circ \theta^-_{\mathbf{i}}(t_1,\cdots,t_N)$.

\vspace{2mm}

\nd
\underline{Case 2-3. $\gamma_{j^+-1}(h_{i_j})=2$}

\vspace{2mm}

Considering the assumption $d_{j^+}=\frac{\gamma_{j^+-1}+\gamma_{j^+}}{2}(h_{i_j})<1$,
we see that $\gamma_{j^+}(h_{i_j})=-2$ so that 
\begin{equation}\label{a3-pr11}
c_{j^+}=\frac{\gamma_{j^+-1}-\gamma_{j^+}}{2}(h_{i_j})=2.
\end{equation}
By (\ref{a3-pr1}), we obtain
\[
\langle h_{i_j}, \sum^{j^+-1}_{l=j+1} c_l\alpha_{i_l}\rangle=-2
\]
and can consider the following three cases:

\vspace{2mm}

\nd
\underline{Case 2-3-1. $a_{i_j,i_{l_1}}=a_{i_j,i_{l_2}}=-1$, $c_{l_1}=c_{l_2}=1$ with some $l_1$, $l_2$ ($j<l_1<l_2<j^+$)}

\vspace{2mm}

For $m\in[j+1,j^+-1]\setminus\{l_1,l_2\}$ such that $a_{i_j,i_m}\neq0$, it holds $c_m=0$.
Since $\pi$ is an $\mathbf{i}$-trail, one obtains
\begin{eqnarray}
0&\neq&
e_{i_1}^{c_1} \cdots 
e_{i_j}^{c_j}\cdots e_{i_{l_1}}^{c_{l_1}}\cdots e_{i_{l_2}}^{c_{l_2}}\cdots e_{i_{j^+}}^{c_{j^+}}\cdots 
e_{i_N}^{c_N}v_{-s_i\Lambda_i}\nonumber\\
&=&
e_{i_1}^{c_1} \cdots 
e_{i_j}^{1}\cdots e_{i_{l_1}}^{1}\cdots e_{i_{l_2}}^{1}\cdots e_{i_{j^+}}^{2}\cdots 
e_{i_N}^{c_N}v_{-s_i\Lambda_i}. \label{a3-pr12}
\end{eqnarray}
By
\begin{equation}\label{a3-prb1}
-2=\gamma_{j^+}(h_{i_j})
=\langle h_{i_j}, {\rm wt}(e_{i_{j^++1}}^{c_{j^++1}}\cdots 
e_{i_N}^{c_N}v_{-s_i\Lambda_i}) \rangle
\end{equation}
and (\ref{mini-pro}),
we obtain $e_{i_{l_2}}^{1}
e_{i_{l_2+1}}^{c_{l_2+1}}
\cdots 
e_{i_{j^+-1}}^{c_{j^+-1}}
e_{i_{j^+}}^{0}e_{i_{j^++1}}^{c_{j^++1}}\cdots 
e_{i_N}^{c_N}v_{-s_i\Lambda_i}=0$ so that
\begin{equation}\label{a3-pr13}
e_{i_j^+}^2
e_{i_{l_2}}^{1}
e_{i_{l_2+1}}^{c_{l_2+1}}
\cdots 
e_{i_{j^+-1}}^{c_{j^+-1}}
e_{i_{j^+}}^{0}e_{i_{j^++1}}^{c_{j^++1}}\cdots 
e_{i_N}^{c_N}v_{-s_i\Lambda_i}=0.
\end{equation}
Combining (\ref{a3-pr12}), (\ref{a3-pr13}) with the Serre relation $e_{i_{j^+}}^2e_{i_{l_2}}-2e_{i_{j^+}}e_{i_{l_2}}e_{i_{j^+}}+e_{i_{l_2}}e_{i_{j^+}}^2=0$, it follows
\begin{multline}\label{a3-pr14}
e_{i_1}^{c_1} \cdots 
e_{i_{j-1}}^{c_{j-1}}
e_{i_j}^{1}
e_{i_{j+1}}^{c_{j+1}}
\cdots 
e_{i_{l_1-1}}^{c_{l_1-1}}
e_{i_{l_1}}^{1}\\
e_{i_{l_1+1}}^{c_{l_1+1}}
\cdots 
e_{i_{l_2-1}}^{c_{l_2-1}}
e_{i_{j^+}}^{1} 
e_{i_{l_2}}^{1}
e_{i_{l_2+1}}^{c_{l_2+1}}
\cdots 
e_{i_{j^+-1}}^{c_{j^+-1}}
e_{i_{j^+}}^{1}e_{i_{j^++1}}^{c_{j^++1}}\cdots 
e_{i_N}^{c_N}v_{-s_i\Lambda_i} \neq0. 
\end{multline}
Similarly, it follows by (\ref{a3-prb1}) that
\[
-1
=\langle h_{i_j}, {\rm wt}(e_{i_{l_2}}^{1}
e_{i_{l_2+1}}^{c_{l_2+1}}
\cdots e_{i_{j^+-1}}^{c_{j^+-1}}e_{i_{j^+}}^{1}e_{i_{j^++1}}^{c_{j^++1}}\cdots 
e_{i_N}^{c_N}v_{-s_i\Lambda_i}) \rangle.
\]
Hence,
one obtains $e_{i_{j^+}}^{2}e_{i_{l_2}}^{1}e_{i_{l_2+1}}^{c_{l_2+1}}\cdots e_{i_{j^+-1}}^{c_{j^+-1}}e_{i_{j^+}}^{1}e_{i_{j^++1}}^{c_{j^++1}}\cdots 
e_{i_N}^{c_N}v_{-s_i\Lambda_i}=0$. 
In particular, we get
\begin{multline*}
e_{i_1}^{c_1} \cdots 
e_{i_{j-1}}^{c_{j-1}}
e_{i_j}^{0}
e_{i_{j+1}}^{c_{j+1}}
\cdots 
e_{i_{l_1-1}}^{c_{l_1-1}}
e_{i_{l_1}}^{1}\\
e_{i_{l_1+1}}^{c_{l_1+1}}
\cdots 
e_{i_{l_2-1}}^{c_{l_2-1}}
e_{i_{j^+}}^{2} 
e_{i_{l_2}}^{1}
e_{i_{l_2+1}}^{c_{l_2+1}}
\cdots 
e_{i_{j^+-1}}^{c_{j^+-1}}
e_{i_{j^+}}^{1}e_{i_{j^++1}}^{c_{j^++1}}\cdots 
e_{i_N}^{c_N}v_{-s_i\Lambda_i} =0. 
\end{multline*}
Consequently, it follows from (\ref{a3-pr14}) and the Serre relation $e_{i_j}^2e_{i_{l_1}}-2e_{i_j}e_{i_{l_1}}e_{i_{j^+}}+e_{i_{l_1}}e_{i_{j^+}}^2=0$ that
\begin{multline}\label{a3-pr15}
e_{i_1}^{c_1} \cdots 
e_{i_{j-1}}^{c_{j-1}}
e_{i_j}^{2}
e_{i_{j+1}}^{c_{j+1}}
\cdots 
e_{i_{l_1-1}}^{c_{l_1-1}}
e_{i_{l_1}}^{1}\\
e_{i_{l_1+1}}^{c_{l_1+1}}
\cdots 
e_{i_{l_2-1}}^{c_{l_2-1}} 
e_{i_{l_2}}^{1}
e_{i_{l_2+1}}^{c_{l_2+1}}
\cdots 
e_{i_{j^+-1}}^{c_{j^+-1}}
e_{i_{j^+}}^{1}e_{i_{j^++1}}^{c_{j^++1}}\cdots 
e_{i_N}^{c_N}v_{-s_i\Lambda_i} \neq0. 
\end{multline}
We define
$\pi'=(\gamma_0',\gamma_1',\cdots,\gamma_N')$ as the pre-$\mathbf{i}$-trail
from $-w_0\Lambda_i$ to $-s_i\Lambda_i$ with integers $c_l'$ $(l=1,2,\cdots,N)$ such that
$\gamma_{l-1}'-\gamma_l'=c_l'\alpha_{i_l}$ and
\[
c_l'=c_{l}\quad \text{for}\ l\in[1,N]\setminus \{j,j^+\},
\]
\[
c_j'=c_{j}+1=2,\quad c_{j^+}'=c_{j^+}-1=1.
\]
Then $\pi'$ is an $\mathbf{i}$-trail by (\ref{a3-pr15}).
Let $M'$ be the monomial
in $\Delta_{w_0\Lm_i,s_i\Lm_i}\circ \theta^-_{\mathbf{i}}(t_1,\cdots,t_N)$ corresponding to $\pi'$.
By Proposition \ref{trail-prop},
we have
\[
M'=MA_{j}^{-1}
\]
and $MA_{j}^{-1}$ is a monomial in $\Delta_{w_0\Lm_i,s_i\Lm_i}\circ \theta^-_{\mathbf{i}}(t_1,\cdots,t_N)$. 

\vspace{2mm}

\nd
\underline{Case 2-3-2. $a_{i_j,i_{l}}=-1$ and $c_{l}=2$ with some $l$ ($j<l<j^+$)}

\vspace{2mm}

For $m\in[j+1,j^+-1]\setminus\{l\}$ such that $a_{i_j,i_m}\neq0$, it holds $c_m=0$.
We see that
\begin{eqnarray*}
0&\neq&
e_{i_1}^{c_1} \cdots 
e_{i_j}^{c_j}\cdots e_{i_{l}}^{2}\cdots e_{i_{j^+}}^{c_{j^+}}\cdots 
e_{i_N}^{c_N}v_{-s_i\Lambda_i}\nonumber\\
&=&
e_{i_1}^{c_1} \cdots 
e_{i_j}^{1}\cdots e_{i_{l}}^{1}e_{i_{l}}^{1}\cdots e_{i_{j^+}}^{2}\cdots 
e_{i_N}^{c_N}v_{-s_i\Lambda_i}. 
\end{eqnarray*}
One can show 
$MA_{j}^{-1}$ is a monomial in $\Delta_{w_0\Lm_i,s_i\Lm_i}\circ \theta^-_{\mathbf{i}}(t_1,\cdots,t_N)$
by a similar way to Case 2-3-1.

\vspace{2mm}

\nd
\underline{Case 2-3-3. $a_{i_j,i_{l}}=-2$ and $c_{l}=1$ with some $l$ ($j<l<j^+$)}

\vspace{2mm}

For $m\in[j+1,j^+-1]\setminus\{l\}$ such that $a_{i_j,i_m}\neq0$, it holds $c_m=0$.
Since $\pi$ is an $\mathbf{i}$-trail,
\begin{eqnarray}
0&\neq&
e_{i_1}^{c_1} \cdots 
e_{i_j}^{c_j}\cdots e_{i_{l}}^{c_{l}}\cdots e_{i_{j^+}}^{c_{j^+}}\cdots 
e_{i_N}^{c_N}v_{-s_i\Lambda_i}\nonumber\\
&=&
e_{i_1}^{c_1} \cdots 
e_{i_j}^{1}\cdots e_{i_{l}}^{1}\cdots e_{i_{j^+}}^{2}\cdots 
e_{i_N}^{c_N}v_{-s_i\Lambda_i}. \label{a3-pr16}
\end{eqnarray}

\nd
Because of the Serre relation
\[
e_{i_j}^3e_{i_l}-3e_{i_j}^2e_{i_l}e_{i_j}+3e_{i_j}e_{i_l}e_{i_j}^2-e_{i_l}e_{i_j}^3=0
\]
and (\ref{mini-pro}),
it holds
\[
e_{i_j}^2e_{i_l}e_{i_j}=e_{i_j}e_{i_l}e_{i_j}^2
\]
on $V(-w_0\Lambda_i)$. Hence, (\ref{a3-pr16}) means
\[
e_{i_1}^{c_1} \cdots 
e_{i_j}^{2}\cdots e_{i_{l}}^{1}\cdots e_{i_{j^+}}^{1}\cdots 
e_{i_N}^{c_N}v_{-s_i\Lambda_i}\neq0.
\]
Therefore, by a similar argument to Case 2-3-1, we see that $M\cdot A_j^{-1}$ is a monomial
in $\Delta_{w_0\Lambda_i,s_i\Lambda_i}\circ \theta^-_{\mathbf{i}}(t_1,\cdots,t_N)$. \qed

In the following lemma and Lemma \ref{lem-a5}, we consider the condition (b) of (2) in Theorem \ref{thm2}.
Let us recall the notation defined in (\ref{jmp}).

\begin{lem}\label{lem-a4}
Let $M=\prod_{l=1}^{N}t_l^{d_l}$ be a monomial in $\Delta_{w_0\Lm_i,s_i\Lm_i}\circ \theta^-_{\mathbf{i}}(t_1,\cdots,t_N)$
and $\{b_l\}_{l\in[1,N]}$ be integers in Definition \ref{bint} determined from $M$.
We suppose that $d_j=d_{j^+}=1$, $b_{j^+}>0$, $b_j=0$ and there exists $p\in\mathbb{Z}_{\geq2}$ such that
\[
d_{j^{m+}}=b_{j^{m+}}=0 \quad(m=2,3,\cdots,p-1)
\]
and $d_{j^{p+}}=-1$, $b_{j^{p+}}=1$.
Then $M\cdot A_{j}^{-1}$ is also a monomial in
$\Delta_{w_0\Lm_i,s_i\Lm_i}\circ \theta^-_{\mathbf{i}}(t_1,\cdots,t_N)$.
\end{lem}

\nd
{\it Proof.} Let
$\pi=(\gamma_0,\gamma_1,\cdots,\gamma_N)$ be the $\mathbf{i}$-trail corresponding to $M$
with integers $c_l$ $(l=1,2,\cdots,N)$ such that
$\gamma_{l-1}-\gamma_l=c_l\alpha_{i_l}$.
The assumption
$d_j=\frac{\gamma_{j-1}+\gamma_{j}}{2}(h_{i_j})=1$ and 
$b_j=c_j=\frac{\gamma_{j-1}-\gamma_{j}}{2}(h_{i_j})=0$ mean
$\gamma_{j-1}(h_{i_j})=\gamma_{j}(h_{i_j})=1$. Just as in (\ref{a2-pr1}), it holds
\begin{equation}\label{a4-pr1}
1=\gamma_j(h_{i_j})
=\gamma_{j^+-1}(h_{i_j})+\langle h_{i_j}, \sum^{j^+-1}_{l=j+1} c_l\alpha_{i_l}\rangle.
\end{equation}
It follows from the assumption
$d_{j^+}=\frac{\gamma_{j^+-1}+\gamma_{j^+}}{2}(h_{i_j})=1$, $b_{j^+}=c_{j^+}=\frac{\gamma_{j^+-1}-\gamma_{j^+}}{2}(h_{i_j})>0$
and Lemma \ref{lem-a1} (2)
that $\gamma_{j^+-1}(h_{i_j})=2$, $\gamma_{j^+}(h_{i_j})=0$ and
\begin{equation}\label{a4-pr2}
c_{j^+}=1.
\end{equation}
In particular, combining $\gamma_{j^+-1}(h_{i_j})=2$ with (\ref{a4-pr1}), it holds
$\langle h_{i_j}, \sum^{j^+-1}_{l=j+1} c_l\alpha_{i_l}\rangle=-1$ so that there uniquely exists $l\in[j+1,j^+-1]$ such that
$a_{i_j,i_l}=-1$ and
\begin{equation}\label{a4-pr3}
c_{l}=1.
\end{equation}
For other $l'\in[j+1,j^+-1]\setminus\{l\}$ such that
$a_{i_j,i_{l'}}\neq0$, we get
\begin{equation}\label{a4-prime}
c_{l'}=0.
\end{equation}
Using the assumption $d_{j^{m+}}=\frac{\gamma_{j^{m+}-1}+\gamma_{j^{m+}}}{2}(h_{i_j})=0$  
and
$b_{j^{m+}}=c_{j^{m+}}=\frac{\gamma_{j^{m+}-1}-\gamma_{j^{m+}}}{2}(h_{i_j})=0$, it holds 
\begin{equation}\label{a4-pr4}
\gamma_{j^{m+}-1}(h_{i_j})=\gamma_{j^{m+}}(h_{i_j})=0
\end{equation}
for $m=2,3,\cdots,p-1$. Using
$d_{j^{p+}}=\frac{\gamma_{j^{p+}-1}+\gamma_{j^{p+}}}{2}(h_{i_j})=-1$, $b_{j^{p+}}=c_{j^{p+}}=\frac{\gamma_{j^{p+}-1}-\gamma_{j^{p+}}}{2}(h_{i_j})=1$,
we obtain
\begin{equation}\label{a4-pr5}
\gamma_{j^{p+}-1}(h_{i_j})=0,\quad \gamma_{j^{p+}}(h_{i_j})=-2.
\end{equation}
Just as in (\ref{a2-pr1}), one can prove
\[
\gamma_{j^{(m-1)+}}(h_{i_j})
=\gamma_{j^{m+}-1}(h_{i_j})+\langle h_{i_j}, \sum^{j^{m+}-1}_{t=j^{(m-1)+}+1} c_t\alpha_{i_t}\rangle
\]
for $m=2,3,\cdots,p$.
In conjunction with $\gamma_{j^+}(h_{i_j})=0$, (\ref{a4-pr4}) and (\ref{a4-pr5}), it holds
\[
\langle h_{i_j}, \sum^{j^{m+}-1}_{t=j^{(m-1)+}+1} c_t\alpha_{i_t}\rangle=0 \quad (m=2,3,\cdots,p).
\]
Hence, for $t\in[j^+,j^{p+}]$ such that $a_{i_j,i_t}<0$, we see that
\begin{equation}\label{a4-pr6}
c_t=0.
\end{equation}
Since $\pi$ is an $\mathbf{i}$-trail, 
\begin{eqnarray}
& &e_{i_1}^{c_1} \cdots 
e_{i_j}^{c_j}\cdots e_{i_l}^{c_l}\cdots e_{i_{j^+}}^{c_{j^+}}\cdots 
e_{i_{j^{p+}}}^{c_{j^{p+}}}\cdots 
e_{i_N}^{c_N}v_{-s_i\Lambda_i} \nonumber \\
&=&
e_{i_1}^{c_1} \cdots 
e_{i_j}^{0}\cdots e_{i_l}^{1}\cdots e_{i_{j^+}}^{1}\cdots 
e_{i_{j^{p+}}}^{1}\cdots 
e_{i_N}^{c_N}v_{-s_i\Lambda_i}
\neq0. \label{a4-pr7}
\end{eqnarray}
Note that it holds $i_j=i_{j^+}=i_{j^{p+}}$.
Considering the Serre relation 
$e_{i_l}e_{i_{j^+}}e_{i_{j^{p+}}}-2e_{i_j}e_{i_l}e_{i_{j^{p+}}} +e_{i_j}^2e_{i_l}=0$
with $c_{j^{m+}}=0$ ($m=2,3,\cdots,p-1$), (\ref{a4-prime}) and (\ref{a4-pr6}), we see that
\begin{eqnarray}
& &e_{i_1}^{c_1} \cdots 
e_{i_j}^{0}\cdots e_{i_l}^{1}\cdots e_{i_{j^+}}^{1}\cdots 
e_{i_{j^{p+}}}^{1}\cdots 
e_{i_N}^{c_N}v_{-s_i\Lambda_i} \nonumber \\
&=&
2e_{i_1}^{c_1} \cdots 
e_{i_j}^{1}\cdots e_{i_l}^{1}\cdots e_{i_{j^+}}^{0}\cdots 
e_{i_{j^{p+}}}^{1}\cdots 
e_{i_N}^{c_N}v_{-s_i\Lambda_i}\nonumber\\
& &-
e_{i_1}^{c_1} \cdots 
e_{i_j}^{2}\cdots e_{i_l}^{1}\cdots e_{i_{j^+}}^{0}\cdots 
e_{i_{j^{p+}}}^{0}\cdots 
e_{i_N}^{c_N}v_{-s_i\Lambda_i}.
\label{a4-pr77}
\end{eqnarray}
Let us show the second term of right-hand side in (\ref{a4-pr77}) is $0$.
Since it is clear in case of the vector 
$e_{i_{l+1}}^{c_{l+1}}\cdots 
e_{i_{j^+-1}}^{c_{j^+-1}}
e_{i_{j^+}}^{0}e_{i_{j^++1}}^{c_{j^++1}}
\cdots 
e_{i_{j^{p+}-1}}^{c_{j^{p+}-1}}
e_{i_{j^{p+}}}^{0}e_{i_{j^{p+}+1}}^{c_{j^{p+}+1}}\cdots 
e_{i_N}^{c_N}v_{-s_i\Lambda_i}$ is $0$, we may assume
this vector is not $0$.
Because of
$-2=\gamma_{j^{p+}}(h_{i_j})=
\langle h_{i_j}, {\rm wt}(e_{i_{j^{p+}+1}}^{c_{j^{p+}+1}}\cdots 
e_{i_N}^{c_N}v_{-s_i\Lambda_i})
\rangle$,
$c_{j^{m+}}=0$ for $m=2,3,\cdots,p-1$,
(\ref{a4-prime}) and
(\ref{a4-pr6}), one gets
\[
\langle h_{i_j}, {\rm wt}(e_{i_{l+1}}^{c_{l+1}}\cdots 
e_{i_{j^+-1}}^{c_{j^+-1}}
e_{i_{j^+}}^{0}e_{i_{j^++1}}^{c_{j^++1}}
\cdots 
e_{i_{j^{p+}-1}}^{c_{j^{p+}-1}}
e_{i_{j^{p+}}}^{0}e_{i_{j^{p+}+1}}^{c_{j^{p+}+1}}\cdots 
e_{i_N}^{c_N}v_{-s_i\Lambda_i})
\rangle=-2
\]
so that
Lemma \ref{lem-a1} (1) yields
\[
e_{i_l}^{1}e_{i_{l+1}}^{c_{l+1}}\cdots 
e_{i_{j^+-1}}^{c_{j^+-1}}
e_{i_{j^+}}^{0}e_{i_{j^++1}}^{c_{j^++1}}
\cdots 
e_{i_{j^{p+}-1}}^{c_{j^{p+}-1}}
e_{i_{j^{p+}}}^{0}e_{i_{j^{p+}+1}}^{c_{j^{p+}+1}}\cdots 
e_{i_N}^{c_N}v_{-s_i\Lambda_i}=0,
\]
which means
the second term of right-hand side in (\ref{a4-pr77}) is $0$.
Thus, using (\ref{a4-pr7}), (\ref{a4-pr77}), it follows
\[
e_{i_1}^{c_1} \cdots 
e_{i_j}^{1}\cdots e_{i_l}^{1}\cdots e_{i_{j^+}}^{0}\cdots 
e_{i_{j^{p+}}}^{1}\cdots 
e_{i_N}^{c_N}v_{-s_i\Lambda_i}
\neq0. 
\]
Therefore, taking
pre-$\mathbf{i}$-trail
$\pi'=(\gamma_0',\gamma_1',\cdots,\gamma_N')$
from $-w_0\Lambda_i$ to $-s_i\Lambda_i$ with integers $c_t'$ $(t=1,2,\cdots,N)$ such that
$\gamma_{t-1}'-\gamma_t'=c_t'\alpha_{i_t}$ and
\[
c_t'=c_{t}\quad \text{for}\ t\in[1,N]\setminus \{j,j^+\},
\]
\[
c_j'=c_{j}+1=1,\quad c_{j^+}'=c_{j^+}-1=0,
\]
we can verify $\pi'$ is an
$\mathbf{i}$-trail.
Let $M'$ be the monomial
in $\Delta_{w_0\Lm_i,s_i\Lm_i}\circ \theta^-_{\mathbf{i}}(t_1,\cdots,t_N)$ corresponding to $\pi'$.
Proposition \ref{trail-prop} implies
$M'=MA_{j}^{-1}$
and $MA_{j}^{-1}$ is a monomial in $\Delta_{w_0\Lm_i,s_i\Lm_i}\circ \theta^-_{\mathbf{i}}(t_1,\cdots,t_N)$. \qed

\begin{lem}\label{lem-a5}
Let $M=\prod_{l=1}^{N}t_l^{d_l}$ be a monomial in $\Delta_{w_0\Lm_i,s_i\Lm_i}\circ \theta^-_{\mathbf{i}}(t_1,\cdots,t_N)$
and $\{b_l\}_{l\in[1,N]}$ be integers in Definition \ref{bint} determined from $M$.
We suppose that $d_j=d_{j^+}=1$, $b_{j^+}>0$, $b_j>0$ and there exists $p\in\mathbb{Z}_{\geq2}$ such that
\[
d_{j^{m+}}=b_{j^{m+}}=0 \quad(m=2,3,\cdots,p-1)
\]
and $d_{j^{p+}}=-1$, $b_{j^{p+}}=1$.
Then $M\cdot A_{j}^{-1}$ is also a monomial in
$\Delta_{w_0\Lm_i,s_i\Lm_i}\circ \theta^-_{\mathbf{i}}(t_1,\cdots,t_N)$.
\end{lem}

Only difference between Lemma \ref{lem-a4} and \ref{lem-a5} is the value of $b_j$.

\vspace{3mm}

\nd
{\it Proof.}

Let
$\pi=(\gamma_0,\gamma_1,\cdots,\gamma_N)$ be the $\mathbf{i}$-trail corresponding to $M$
with integers $c_l$ $(l=1,2,\cdots,N)$ such that
$\gamma_{l-1}-\gamma_l=c_l\alpha_{i_l}$.
The assumption
$d_j=\frac{\gamma_{j-1}+\gamma_{j}}{2}(h_{i_j})=1$,
$b_j=c_j=\frac{\gamma_{j-1}-\gamma_{j}}{2}(h_{i_j})>0$ and Lemma \ref{lem-a1} (2) mean
$\gamma_{j-1}(h_{i_j})=2$, $\gamma_{j}(h_{i_j})=0$ and $b_j=1$. Just as in (\ref{a2-pr1}), it holds
\begin{equation}\label{a5-pr1}
0=\gamma_j(h_{i_j})
=\gamma_{j^+-1}(h_{i_j})+\langle h_{i_j}, \sum^{j^+-1}_{l=j+1} c_l\alpha_{i_l}\rangle.
\end{equation}
By the assumption $d_{j^+}=\frac{\gamma_{j^+-1}+\gamma_{j^+}}{2}(h_{i_j})=1$,
$\frac{\gamma_{j^+-1}-\gamma_{j^+}}{2}(h_{i_j})=c_{j^+}=b_{j^+}>0$ and Lemma \ref{lem-a1} (2),
we obtain $\gamma_{j^+-1}(h_{i_j})=2$, $\gamma_{j^+}(h_{i_j})=0$ and
\begin{equation}\label{a5-pr2}
c_{j^+}=1.
\end{equation}
By the same way as in (\ref{a4-pr4}), (\ref{a4-pr5}), it holds
\[
\gamma_{j^{m+}-1}(h_{i_j})=\gamma_{j^{m+}}(h_{i_j})=0 \quad (m=2,3,\cdots,p-1)
\]
and
\begin{equation}\label{a5-pr3}
\gamma_{j^{p+}-1}(h_{i_j})=0,\quad \gamma_{j^{p+}}(h_{i_j})=-2.
\end{equation}
Just as in
(\ref{a4-pr6}), we see that for any $t\in[j^+,j^{p+}]$ such that $a_{i_j,i_t}<0$, it follows
\begin{equation}\label{a5-pr2}
c_{t}=0.
\end{equation}
Since we know $\gamma_{j^+-1}(h_{i_j})=2$ and (\ref{a5-pr1}), 
it follows $\langle h_{i_j}, \sum^{j^+-1}_{l=j+1} c_l\alpha_{i_l}\rangle=-2$.
Hence, there are following three patterns:

\vspace{3mm}

\nd
\underline{Case 1. $a_{i_j,i_{l_1}}=a_{i_j,i_{l_2}}=-1$, $c_{l_1}=c_{l_2}=1$ with some $l_1$, $l_2$ ($j<l_1<l_2<j^+)$}
\vspace{2mm}

For $l'\in[j+1,j^+-1]\setminus\{l_1,l_2\}$ such that $a_{i_j,i_{l'}}<0$, we obtain
\begin{equation}\label{a5-pr4}
c_{l'}=0.
\end{equation}
Since $\pi$ is an $\mathbf{i}$-trail, 
\begin{eqnarray}
& &
e_{i_1}^{c_1}\cdots e_{i_{l_2}}^{c_{l_2}}\cdots e_{i_{j^+}}^{c_{j^+}}
\cdots e_{i_{j^{p+}}}^{c_{j^{p+}}}\cdots e_{i_N}^{c_{i_N}}v_{-s_i\Lambda_i}\nonumber\\
&=&e_{i_1}^{c_1}\cdots e_{i_{l_2}}^{1}\cdots e_{i_{j^+}}^{1}
\cdots e_{i_{j^{p+}}}^{1}\cdots e_{i_N}^{c_{i_N}}v_{-s_i\Lambda_i}\neq0. \label{a5-pr5}
\end{eqnarray}
By (\ref{a5-pr3}), we obtain
\[
-2=\gamma_{j^{p+}}(h_{i_j})
=\langle
h_{i_j}, {\rm wt}(e_{i_{j^{p+}+1}}^{c_{j^{p+}+1}}\cdots e_{i_N}^{c_{i_N}}v_{-s_i\Lambda_i})
\rangle.
\]
Taking this equation, $c_{j^{m+}}=0$ $(m=2,3,\cdots,p-1)$, (\ref{a5-pr2}) and (\ref{a5-pr4}) into account, we see that
\[
e_{i_{l_2}}^{1}
e_{i_{l_2+1}}^{c_{l_2+1}}
\cdots 
e_{i_{j^+-1}}^{c_{j^+-1}}
e_{i_{j^+}}^{0}
e_{i_{j^++1}}^{c_{j^++1}}
\cdots
e_{i_{j^{p+}-1}}^{c_{j^{p+}-1}}
e_{i_{j^{p+}}}^{0}
e_{i_{j^{p+}+1}}^{c_{j^{p+}+1}}
\cdots e_{i_N}^{c_{i_N}}v_{-s_i\Lambda_i}=0,
\]
otherwise, it holds
\[\langle h_{i_j},{\rm wt}(e_{i_{l_2}}^{1}e_{i_{l_2+1}}^{c_{l_2+1}}\cdots 
e_{i_{j^+-1}}^{c_{j^+-1}}
e_{i_{j^+}}^{0}
e_{i_{j^++1}}^{c_{j^++1}}
\cdots
e_{i_{j^{p+}-1}}^{c_{j^{p+}-1}}
e_{i_{j^{p+}}}^{0}
e_{i_{j^{p+}+1}}^{c_{j^{p+}+1}}
\cdots e_{i_N}^{c_{i_N}}v_{-s_i\Lambda_i})\rangle=-3,
\]
which contradicts Lemma \ref{lem-a1} (1). In particular, one gets
\begin{equation}\label{a5-pr6}
e_{i_j}^2
e_{i_{l_2}}^{1}e_{i_{l_2+1}}^{c_{l_2+1}}\cdots 
e_{i_{j^+-1}}^{c_{j^+-1}}
e_{i_{j^+}}^{0}
e_{i_{j^++1}}^{c_{j^++1}}
\cdots
e_{i_{j^{p+}-1}}^{c_{j^{p+}-1}}
e_{i_{j^{p+}}}^{0}
e_{i_{j^{p+}+1}}^{c_{j^{p+}+1}}
\cdots e_{i_N}^{c_{i_N}}v_{-s_i\Lambda_i}=0.
\end{equation}
Note that $i_j=i_{j^+}=i_{j^{p+}}$.
Hence, (\ref{a5-pr2}), (\ref{a5-pr4}), (\ref{a5-pr5}), (\ref{a5-pr6}) and the Serre relation
$e_{i_{j}}^2e_{i_{l_2}}-2 e_{i_{j^+}}e_{i_{l_2}}e_{i_{j^{p+}}} + e_{i_{l_2}}e_{i_{j^+}}e_{i_{j^{p+}}}=0$
mean
\begin{eqnarray*}
0&\neq&e_{i_1}^{c_1}\cdots e_{i_{l_2}}^{c_{l_2}}\cdots e_{i_{j^+}}^{c_{j^+}}
\cdots e_{i_{j^{p+}}}^{c_{j^{p+}}}\cdots e_{i_N}^{c_{i_N}}v_{-s_i\Lambda_i}\\
&=&e_{i_1}^{c_1}\cdots e_{i_{l_2}}^{1}\cdots e_{i_{j^+}}^{1}
\cdots e_{i_{j^{p+}}}^{1}\cdots e_{i_N}^{c_{i_N}}v_{-s_i\Lambda_i}\\
&=&
2e_{i_1}^{c_1}\cdots
e_{i_{j^+}}^1e_{i_{l_2}}^{1}\cdots 
e_{i_{j^+}}^{0}
\cdots e_{i_{j^{p+}}}^{1}\cdots e_{i_N}^{c_{i_N}}v_{-s_i\Lambda_i}.
\end{eqnarray*}
Since $c_{j^{p+}}=1$ and $c_{l_2}=1$, the above vector is equal to
\begin{equation}\label{a5-pr7}
2e_{i_1}^{c_1}\cdots
e_{i_{l_2-1}}^{c_{l_2-1}}
e_{i_{j^+}}e_{i_{l_2}}^{c_{l_2}}\cdots 
e_{i_{j^+-1}}^{c_{j^+-1}}
e_{i_{j^+}}^{0}
e_{i_{j^++1}}^{c_{j^++1}}
\cdots e_{i_{j^{p+}}}^{c_{j^{p+}}}\cdots e_{i_N}^{c_{i_N}}v_{-s_i\Lambda_i} (\neq0).
\end{equation}
Considering
$0=\gamma_{j^+}(h_{i_j})=\langle h_{i_j}, {\rm wt}(e_{i_{j^++1}}^{c_{j^++1}}\cdots e_{i_N}^{c_{i_N}}v_{-s_i\Lambda_i})\rangle$, (\ref{a5-pr4}) and $c_{l_2}=1$, it holds
\[
\langle h_{i_j},
{\rm wt}(e_{i_{l_2}}^{c_{l_2}}\cdots e_{i_{j^+-1}}^{c_{j^+-1}}
e_{i_{j^+}}^{0}e_{i_{j^++1}}^{c_{j^++1}}
\cdots e_{i_N}^{c_{i_N}}v_{-s_i\Lambda_i})\rangle=-1
\]
so that $e_{i_j^+}^2e_{i_{l_2}}^{c_{l_2}}\cdots e_{i_{j^+-1}}^{c_{j^+-1}}
e_{i_{j^+}}^{0}e_{i_{j^++1}}^{c_{j^++1}}
\cdots e_{i_N}^{c_{i_N}}v_{-s_i\Lambda_i}=0$. Hence,
we have
\[
e_{i_{l_1}}^{1}e_{i_{l_1+1}}^{c_{l_1+1}}\cdots 
e_{i_{l_2-1}}^{c_{l_2-1}}e_{i_j^+}^2e_{i_{l_2}}^{c_{l_2}}\cdots e_{i_{j^+-1}}^{c_{j^+-1}}
e_{i_{j^+}}^{0}e_{i_{j^++1}}^{c_{j^++1}}
\cdots e_{i_N}^{c_{i_N}}v_{-s_i\Lambda_i}=0.
\]
Thus, applying the Serre relation
$e_{i_j}^{2}e_{i_{l_1}}-2e_{i_j}e_{i_{l_1}}e_{i_{j^+}}+e_{i_{l_1}}e_{i_{j^+}}^2=0$,
the vector in (\ref{a5-pr7}) is transformed as follows:
\begin{eqnarray}
0&\neq&
2e_{i_1}^{c_1}\cdots e_{i_j}^{c_j}\cdots e_{i_{l_1}}^{c_{l_1}}\cdots 
e_{i_{l_2-1}}^{c_{l_2-1}}
e_{i_{j^+}}e_{i_{l_2}}^{c_{l_2}}\cdots e_{i_{j^+}}^{0}
\cdots e_{i_{j^{p+}}}^{c_{j^{p+}}}\cdots e_{i_N}^{c_{i_N}}v_{-s_i\Lambda_i}\nonumber\\
&=&
2e_{i_1}^{c_1}\cdots e_{i_j}^{1}\cdots e_{i_{l_1}}^{1}\cdots 
e_{i_{l_2-1}}^{c_{l_2-1}}
e_{i_{j^+}}e_{i_{l_2}}^{c_{l_2}}\cdots e_{i_{j^+}}^{0}
\cdots e_{i_{j^{p+}}}^{c_{j^{p+}}}\cdots e_{i_N}^{c_{i_N}}v_{-s_i\Lambda_i}\nonumber\\
&=&
e_{i_1}^{c_1}\cdots e_{i_j}^{2}\cdots e_{i_{l_1}}^{1}\cdots e_{i_{l_2-1}}^{c_{l_2-1}}e_{i_{l_2}}^{c_{l_2}}\cdots e_{i_{j^+}}^{0}
\cdots e_{i_{j^{p+}}}^{c_{j^{p+}}}\cdots e_{i_N}^{c_{i_N}}v_{-s_i\Lambda_i}\nonumber\\
&=&
e_{i_1}^{c_1}\cdots e_{i_j}^{c_j+1}\cdots e_{i_{l_1}}^{c_{l_1}}\cdots e_{i_{l_2}}^{c_{l_2}}\cdots e_{i_{j^+}}^{c_{j^+}-1}
\cdots e_{i_{j^{p+}}}^{c_{j^{p+}}}\cdots e_{i_N}^{c_{i_N}}v_{-s_i\Lambda_i}. \label{a5-pr8}
\end{eqnarray}
One defines
$\pi'=(\gamma_0',\gamma_1',\cdots,\gamma_N')$ as the pre-$\mathbf{i}$-trail
from $-w_0\Lambda_i$ to $-s_i\Lambda_i$ with integers $c_l'$ $(l=1,2,\cdots,N)$ such that
$\gamma_{l-1}'-\gamma_l'=c_l'\alpha_{i_l}$ and
\[
c_l'=c_{l}\quad \text{for}\ l\in[1,N]\setminus \{j,j^+\},
\]
\[
c_j'=c_{j}+1,\quad c_{j^+}'=c_{j^+}-1=0.
\]
Then $\pi'$ is an $\mathbf{i}$-trail by (\ref{a5-pr8}).
Let $M'$ be the monomial
in $\Delta_{w_0\Lm_i,s_i\Lm_i}\circ \theta^-_{\mathbf{i}}(t_1,\cdots,t_N)$ corresponding to $\pi'$.
By Proposition \ref{trail-prop},
we have
$M'=MA_{j}^{-1}$
and $MA_{j}^{-1}$ is a monomial in $\Delta_{w_0\Lm_i,s_i\Lm_i}\circ \theta^-_{\mathbf{i}}(t_1,\cdots,t_N)$. 

\vspace{3mm}

\nd
\underline{Case 2. $a_{i_j,i_{l}}=-1$, $c_{l}=2$ with some $l$ ($j<l<j^+$)}

\vspace{2mm}

By a similar way to Case 1, we see that
$MA_{j}^{-1}$ is a monomial in $\Delta_{w_0\Lm_i,s_i\Lm_i}\circ \theta^-_{\mathbf{i}}(t_1,\cdots,t_N)$.

\vspace{3mm}

\nd
\underline{Case 3. $a_{i_j,i_{l}}=-2$, $c_{l}=1$ with some $l$ ($j<l<j^+$)}

\vspace{2mm}

For $l'\in[j+1,j^+-1]\setminus\{l\}$ such that $a_{i_j,i_{l'}}<0$, it holds $c_{l'}=0$.
Since $\pi$ is an $\mathbf{i}$-trail, 
\begin{eqnarray*}
& &e_{i_1}^{c_1}\cdots e_{i_{j}}^{c_j}\cdots e_{i_{l}}^{c_{l}}\cdots e_{i_{j^+}}^{c_{j^+}}
\cdots e_{i_{j^{p+}}}^{c_{j^{p+}}}\cdots e_{i_N}^{c_{i_N}}v_{-s_i\Lambda_i}\\
&=&
e_{i_1}^{c_1}\cdots e_{i_{j}}^{1}\cdots e_{i_{l}}^{1}\cdots e_{i_{j^+}}^{1}
\cdots e_{i_{j^{p+}}}^{1}\cdots e_{i_N}^{c_{i_N}}v_{-s_i\Lambda_i}\neq0.
\end{eqnarray*}
Because of the Serre relation
\[
e_{i_j}^3e_{i_l}-3e_{i_j}^2e_{i_l}e_{i_j}+3e_{i_j}e_{i_l}e_{i_j}^2-e_{i_l}e_{i_j}^3=0
\]
and (\ref{mini-pro}),
it holds
\[
e_{i_j}^2e_{i_l}e_{i_j}=e_{i_j}e_{i_l}e_{i_j}^2
\]
on $V(-w_0\Lambda_i)$. Thus, using $e_{i_j}e_{i_l}e_{i_{j^+}}e_{i_{j^{p+}}}=e_{i_j}^2e_{i_l}e_{i_{j^{p+}}}$,
we get
\begin{eqnarray*}
0&\neq&e_{i_1}^{c_1}\cdots e_{i_{j}}^{1}\cdots e_{i_{l}}^{1}\cdots e_{i_{j^+}}^{1}
\cdots e_{i_{j^{p+}}}^{1}\cdots e_{i_N}^{c_{i_N}}v_{-s_i\Lambda_i}\\
&=&
e_{i_1}^{c_1}\cdots e_{i_{j}}^{2}\cdots e_{i_{l}}^{1}\cdots e_{i_{j^+}}^{0}
\cdots e_{i_{j^{p+}}}^{1}\cdots e_{i_N}^{c_{i_N}}v_{-s_i\Lambda_i}\\
&=&
e_{i_1}^{c_1}\cdots e_{i_{j}}^{c_{j}+1}\cdots e_{i_{l}}^{c_l}\cdots e_{i_{j^+}}^{c_{j^+}-1}
\cdots e_{i_{j^{p+}}}^{c_{j^{p+}}}\cdots e_{i_N}^{c_{i_N}}v_{-s_i\Lambda_i}.
\end{eqnarray*}
Therefore, taking
$\pi'=(\gamma_0',\gamma_1',\cdots,\gamma_N')$ as the pre-$\mathbf{i}$-trail
from $-w_0\Lambda_i$ to $-s_i\Lambda_i$ with integers $c_r'$ $(r=1,2,\cdots,N)$ such that
$\gamma_{r-1}'-\gamma_r'=c_r'\alpha_{i_r}$ and
\[
c_r'=c_{r}\quad \text{for}\ r\in[1,N]\setminus \{j,j^+\},
\]
\[
c_j'=c_{j}+1,\quad c_{j^+}'=c_{j^+}-1,
\]
one can verify
$\pi'$ is an $\mathbf{i}$-trail.
Let $M'$ be the monomial
in $\Delta_{w_0\Lm_i,s_i\Lm_i}\circ \theta^-_{\mathbf{i}}(t_1,\cdots,t_N)$ corresponding to $\pi'$.
By Proposition \ref{trail-prop},
we have $M'=MA_{j}^{-1}$
and $MA_{j}^{-1}$ is a monomial in $\Delta_{w_0\Lm_i,s_i\Lm_i}\circ \theta^-_{\mathbf{i}}(t_1,\cdots,t_N)$. \qed

\begin{lem}\label{Amplus}
Let $M=\prod_{l=1}^N t_l^{d_l}$ be a monomial in $\Delta_{w_0\Lambda_i,s_i\Lambda_i}\circ \theta^-_{\mathbf{i}}(t_1,\cdots,t_N)$.
We assume $d_m<0$ with some $m\in[1,N]$. Then there exists $r\in\mathbb{Z}_{\geq1}$ such that
$M\cdot A_{m^{r-}}$ is a monomial in $\Delta_{w_0\Lambda_i,s_i\Lambda_i}\circ \theta^-_{\mathbf{i}}(t_1,\cdots,t_N)$ and
putting $j:=m^{r-}$ and $M\cdot A_{m^{r-}}=\prod_{l=1}^N t_l^{d_l'}$, it follows $d_j'>0$, $b_{j^+}'>0$ and one of the following
holds:
\begin{enumerate}
\item[(a)] $d_{j^+}'< d_j'$,
\item[(b)] $d_{j^+}'= d_j'$ and there exists
$p\in\mathbb{Z}_{\geq2}$ such that
\[
d_{j^{s+}}'=b_{j^{s+}}'=0 \quad(s=2,3,\cdots,p-1)
\]
and $d_{j^{p+}}'=-1$, $b_{j^{p+}}'=1$.
\end{enumerate}
Here, $\{b_l'\}_{l\in[1,N]}$ are integers in Definition \ref{bint} determined from $M\cdot A_{m^{r-}}$.
\end{lem}
\nd
{\it Proof.} Let $\pi=(\gamma_0,\cdots,\gamma_N)$ be the $\mathbf{i}$-trail corresponding to $M$
with nonnegative integers $\{c_l\}_{l=1}^N$.
Let $\xi:={\rm max}\{z\in\mathbb{Z}_{\geq0} | m^{z-}\geq1\}$.
We may assume that
\begin{equation}\label{Amplus-pr1}
d_{m^-},d_{m^{2-}},\cdots, d_{m^{\xi-}}\geq0.
\end{equation}
Note that it holds $\gamma_{m-1}(h_{i_m})=\gamma_m(h_{i_m})+2c_m\geq \gamma_{m}(h_{i_m})$.
By the assumption $d_m=\frac{\gamma_{m-1}+\gamma_{m}}{2}(h_{i_m})<0$ and Lemma \ref{lem-a1} (2),
one of the following holds:
\begin{enumerate}
\item[(i)] $\gamma_{m-1}(h_{i_m})=\gamma_m(h_{i_m})=-2$,
\item[(ii)] $\gamma_{m-1}(h_{i_m})=\gamma_m(h_{i_m})=-1$,
\item[(iii)] $\gamma_{m-1}(h_{i_m})=0$, $\gamma_m(h_{i_m})=-2$.
\end{enumerate}
\nd
\underline{The proof for the case (i) : $\gamma_{m-1}(h_{i_m})=\gamma_m(h_{i_m})=-2$}

\vspace{2mm}

In this case, it holds
\begin{equation}\label{a6-pr1}
d_m=\frac{\gamma_{m-1}(h_{i_m})+\gamma_{m}(h_{i_m})}{2}=-2,\quad b_m=c_m=\frac{\gamma_{m-1}(h_{i_m})-\gamma_{m}(h_{i_m})}{2}=0.
\end{equation}
Considering
\[
-2=\gamma_{m-1}(h_{i_m})
=\langle h_{i_m}, {\rm wt}(e_{i_m}^{c_m}\cdots e_{i_N}^{c_N}v_{-s_i\Lambda_i})\rangle,
\]
if $m^-=0$ then $\langle h_{i_m}, {\rm wt}(e_{i_1}^{c_1}\cdots e_{i_m}^{c_m}\cdots e_{i_N}^{c_N}v_{-s_i\Lambda_i})\rangle\leq -2$,
which contradicts ${\rm wt}(e_{i_1}^{c_1}\cdots e_{i_m}^{c_m}\cdots e_{i_N}^{c_N}v_{-s_i\Lambda_i})=-w_0\Lambda_i$. Thus,
we get $m^-\geq1$. Just as in (\ref{a2-pr1}), it holds
\begin{equation}\label{a6-pr2}
\gamma_{m^-}(h_{i_m})=\gamma_{m-1}(h_{i_m})+\langle h_{i_m}, \sum_{l=m^-+1}^{m-1} c_l\alpha_{i_l} \rangle,
\end{equation}
which yields $\gamma_{m^-}(h_{i_m})\leq \gamma_{m-1}(h_{i_m})=-2$ so that $\gamma_{m^-}(h_{i_m})=-2$ by Lemma \ref{lem-a1} (2). 
It follows by (\ref{a6-pr2}) that
\begin{equation}\label{a6-pr3}
c_l=0
\end{equation}
for $l\in[m^-+1,m-1]$ such that $a_{i_m,i_l}<0$.
The inequality $d_{m^-}=\frac{\gamma_{m^--1}(h_{i_m}) +\gamma_{m^-}(h_{i_m})}{2}\leq0$ follows from $\gamma_{m^-}(h_{i_m})=-2$
and Lemma \ref{lem-a1} (2). By (\ref{Amplus-pr1}), one gets $\gamma_{m^--1}(h_{i_m})=2$, $d_{m^-}=0$ and
\begin{equation}\label{a6-pr4}
b_{m^-}=c_{m^-}=\frac{\gamma_{m^--1}(h_{i_m})-\gamma_{m^-}(h_{i_m})}{2}=2.
\end{equation}
Since $\pi$ is an $\mathbf{i}$-trail, it follows
$e_{i_1}^{c_1}\cdots e_{i_N}^{c_N}v_{-s_i\Lambda_i}\neq0$.
Let 
$\pi'=(\gamma_0',\gamma_1',\cdots,\gamma_N')$ be the pre-$\mathbf{i}$-trail
from $-w_0\Lambda_i$ to $-s_i\Lambda_i$ with integers $c_l'$ $(l=1,2,\cdots,N)$ such that
$\gamma_{l-1}'-\gamma_l'=c_l'\alpha_{i_l}$ and
\[
c_l'=c_{l}\quad \text{for}\ l\in[1,N]\setminus \{m^-,m\},
\]
\[
c_{m^-}'=c_{m^-}-1=1,\quad c_{m}'=c_{m}+1=1.
\]
One can verify
\[
0\neq e_{i_1}^{c_1}\cdots e_{i_N}^{c_N}v_{-s_i\Lambda_i} =
e_{i_1}^{c_1'}\cdots e_{i_N}^{c_N'}v_{-s_i\Lambda_i}
\]
by using (\ref{a6-pr3}) so that $\pi'$ is an $\mathbf{i}$-trail.
Let $M'$ be the monomial
in $\Delta_{w_0\Lm_i,s_i\Lm_i}\circ \theta^-_{\mathbf{i}}(t_1,\cdots,t_N)$ corresponding to $\pi'$.
By Proposition \ref{trail-prop},
we have $M'=MA_{m^-}$.
Let $j=m^-$, $M'=MA_{j}=\prod^N_{l=1}t_l^{d_l'}$ and $\{b_{l}'\}_{l\in[1,N]}$ be integers 
in Definition \ref{bint} determined from $M'$.
Then it holds $d_j'=d_{j}+1=d_{m^-}+1=1>0$, $b_{j^+}'=c_{j^+}'=c_{m}'=1>0$ and
$d_{j^+}'=d_{j^+}+1=d_{m}+1=-1$ so that (a) in our claim follows.

\vspace{2mm}

\nd
\underline{The proof for the case (ii) : $\gamma_{m-1}(h_{i_m})=\gamma_m(h_{i_m})=-1$}

\vspace{2mm}

In this case, we obtain $b_m=c_m=\frac{\gamma_{m-1}(h_{i_m})-\gamma_m(h_{i_m})}{2}=0$ and
$d_m=\frac{\gamma_{m-1}(h_{i_m})+\gamma_m(h_{i_m})}{2}=-1$. If $m^-=0$ then
by
\[
\gamma_{m-1}(h_{i_m})=\langle h_{i_m}, {\rm wt}(e_{i_m}^{c_m}\cdots e_{i_N}^{c_N}v_{-s_i\Lambda_i}) \rangle,
\]
one gets
\[
\langle h_{i_m}, {\rm wt}(e_{i_1}^{c_1}\cdots e_{i_m}^{c_m}\cdots e_{i_N}^{c_N}v_{-s_i\Lambda_i}) \rangle\leq -1,
\]
which contradicts ${\rm wt}(e_{i_1}^{c_1}\cdots e_{i_m}^{c_m}\cdots e_{i_N}^{c_N}v_{-s_i\Lambda_i})=-w_0\Lambda_i$.
Hence, it holds $m^-\geq1$. By the same argument as in (\ref{a6-pr2}), it follows
\begin{equation}\label{a6-pr5}
\gamma_{m^-}(h_{i_m})=\gamma_{m-1}(h_{i_m})+\langle h_{i_m}, \sum_{l=m^-+1}^{m-1} c_l\alpha_{i_l} \rangle,
\end{equation}
which yields $\gamma_{m^-}(h_{i_m})\leq \gamma_{m-1}(h_{i_m})$.
Now, we assumed $\gamma_{m-1}(h_{i_m})=-1$ so that $\gamma_{m^-}(h_{i_m})=-1$ or $\gamma_{m^-}(h_{i_m})=-2$.
\vspace{2mm}

\nd
\underline{Case 1. $\gamma_{m^-}(h_{i_m})=-1$}

\vspace{2mm}

Considering (\ref{a6-pr5}), for $l\in[m^-+1,m-1]$ such that $a_{i_m,i_l}<0$, it holds $c_l=0$.
Taking the assumption $0\leq d_{m^-}=\frac{\gamma_{m^--1}(h_{i_m})+\gamma_{m^-}(h_{i_m})}{2}$ in (\ref{Amplus-pr1})
and Lemma \ref{lem-a1} (2) into account, one obtains $\gamma_{m^--1}(h_{i_m})=1$ so that
$d_{m^-}=0$ and $b_{m^-}=c_{m^-}=\frac{\gamma_{m^--1}(h_{i_m})-\gamma_{m^-}(h_{i_m})}{2}=1$.
Let 
$\pi'=(\gamma_0',\gamma_1',\cdots,\gamma_N')$ be the pre-$\mathbf{i}$-trail
from $-w_0\Lambda_i$ to $-s_i\Lambda_i$ with integers $c_l'$ $(l=1,2,\cdots,N)$ such that
$\gamma_{l-1}'-\gamma_l'=c_l'\alpha_{i_l}$ and
\[
c_l'=c_{l}\quad \text{for}\ l\in[1,N]\setminus \{m^-,m\},
\]
\[
c_{m^-}'=c_{m^-}-1=0,\quad c_{m}'=c_{m}+1=1.
\]
By the same argument as in (i), the pre-$\mathbf{i}$-trail $\pi'$ is an $\mathbf{i}$-trail.
Let $M'$ be the monomial
in $\Delta_{w_0\Lm_i,s_i\Lm_i}\circ \theta^-_{\mathbf{i}}(t_1,\cdots,t_N)$ corresponding to $\pi'$.
By Proposition \ref{trail-prop},
we have $M'=MA_{m^-}$.
Let $j=m^-$, $M'=MA_{j}=\prod^N_{l=1}t_l^{d_l'}$ and $\{b_{l}'\}_{l\in[1,N]}$ be integers 
in Definition \ref{bint} determined from $M'$.
Then it holds $d_j'=d_{m^-}+1=1>0$, $b_{j^+}'=c_{j^+}'=c_{m}'=1>0$ and $d_{j^+}'=d_m+1=0$ so that (a) in our claim is satisfied.

\vspace{2mm}

\nd
\underline{Case 2. $\gamma_{m^-}(h_{i_m})=-2$}

\vspace{2mm}

Using (\ref{a6-pr5}), one can verify there exists $l\in[m^-+1,m-1]$ such that $c_l=1$ and $a_{i_m,i_l}=-1$.
For $l'\in[m^-+1,m-1]\setminus\{l\}$ such that $a_{i_m,i_{l'}}<0$, it holds $c_{l'}=0$.
The assumption $0\leq d_{m^-}=\frac{\gamma_{m^--1}(h_{i_m})+\gamma_{m^{-}}(h_{i_m})}{2}$ in (\ref{Amplus-pr1}) yields
$\gamma_{m^--1}(h_{i_m})=2$ so that $d_{m^-}=0$ and $b_{m^-}=c_{m^-}=\frac{\gamma_{m^{-}-1}(h_{i_m})-\gamma_{m^{-}}(h_{i_m})}{2}=2$.
Considering Lemma \ref{lem-a1} (1) and
\[
-1=\gamma_m(h_{i_m})=\langle h_{i_m}, {\rm wt}(e_{i_{m+1}}^{c_{m+1}}\cdots e_{i_N}^{c_N}v_{-s_i\Lambda_i})\rangle,
\]
we get $e_{i_m}^2e_{i_{m+1}}^{c_{m+1}}\cdots e_{i_N}^{c_N}v_{-s_i\Lambda_i}=0$. In particular,
\begin{equation}\label{a6-pr6}
e_{i_1}^{c_1}\cdots e_{i_{m^-}}^0 \cdots e_{i_l}^1 \cdots e_{i_m}^2e_{i_{m+1}}^{c_{m+1}}\cdots e_{i_N}^{c_N}v_{-s_i\Lambda_i}=0.
\end{equation}
The definition of $\mathbf{i}$-trail implies
\begin{eqnarray*}
0&\neq& e_{i_1}^{c_1}\cdots e_{i_{m^-}}^{c_{m^-}} \cdots e_{i_l}^{c_l} \cdots e_{i_m}^{c_m}e_{i_{m+1}}^{c_{m+1}}\cdots e_{i_N}^{c_N}v_{-s_i\Lambda_i}\\
&=& e_{i_1}^{c_1}\cdots e_{i_{m^-}}^{2} \cdots e_{i_l}^{1} \cdots e_{i_m}^{0}e_{i_{m+1}}^{c_{m+1}}\cdots e_{i_N}^{c_N}v_{-s_i\Lambda_i}.
\end{eqnarray*}
Combining this with (\ref{a6-pr6}) and the Serre relation $e_{i_{m^-}}^2e_{i_l}-2e_{i_{m^-}}e_{i_l}e_{i_m}+e_{i_l}e_{i_m}^2=0$, 
we see that
\begin{eqnarray*}
0&\neq&
e_{i_1}^{c_1}\cdots e_{i_{m^-}}^{1} \cdots e_{i_l}^{1} \cdots e_{i_m}^{1}e_{i_{m+1}}^{c_{m+1}}\cdots e_{i_N}^{c_N}v_{-s_i\Lambda_i}\\
&=&
e_{i_1}^{c_1}\cdots e_{i_{m^-}}^{c_{m^-}-1} \cdots e_{i_l}^{c_l} \cdots e_{i_m}^{c_m+1}e_{i_{m+1}}^{c_{m+1}}\cdots e_{i_N}^{c_N}v_{-s_i\Lambda_i}.
\end{eqnarray*}
Hence, the pre-$\mathbf{i}$-trail $\pi'=(\gamma_0',\gamma_1',\cdots,\gamma_N')$
from $-w_0\Lambda_i$ to $-s_i\Lambda_i$ with integers $c_t'$ $(t=1,2,\cdots,N)$ such that
$\gamma_{t-1}'-\gamma_t'=c_t'\alpha_{i_t}$ and
\[
c_t'=c_{t}\quad \text{for}\ t\in[1,N]\setminus \{m^-,m\},
\]
\[
c_{m^-}'=c_{m^-}-1=1,\quad c_{m}'=c_{m}+1=1
\]
is an $\mathbf{i}$-trail.
Let $M'$ be the monomial
in $\Delta_{w_0\Lm_i,s_i\Lm_i}\circ \theta^-_{\mathbf{i}}(t_1,\cdots,t_N)$ corresponding to $\pi'$.
By Proposition \ref{trail-prop},
we have $M'=MA_{m^-}$.
Let $j=m^-$, $M'=MA_{j}=\prod^N_{r=1}t_r^{d_r'}$ and $\{b_{r}'\}_{r\in[1,N]}$ be integers 
in Definition \ref{bint} determined from $M'$.
Then it holds $d_j'=d_{m^-}+1=1>0$, $b_{j^+}'=c_{j^+}'=c_{m}'=1>0$ and $d_{j^+}'=d_m+1=0$ so that (a) in our claim is satisfied.

\vspace{2mm}

\nd
\underline{The proof for the case (iii) : $\gamma_{m-1}(h_{i_m})=0$, $\gamma_m(h_{i_m})=-2$}

\vspace{2mm}

In this case, we see that $d_m=\frac{\gamma_{m-1}(h_{i_m})+\gamma_m(h_{i_m})}{2}=-1$ and
$c_m=\frac{\gamma_{m-1}(h_{i_m})-\gamma_m(h_{i_m})}{2}=1$.
By $\gamma_m(h_{i_m})=-2$, it holds
\begin{equation}\label{Amplus-pr2}
\langle
h_{i_m}, {\rm wt}(e_{i_{m+1}}^{c_{m+1}}\cdots e_{i_N}^{c_N}v_{-s_i\Lambda_i})
\rangle=-2.
\end{equation}

\nd
\underline{Case 1. $\gamma_{m^{l-}}(h_{i_m})=0$ for all $l\in[1,\xi]$}

\vspace{2mm}

\nd
\underline{Case 1-1. $\xi\geq1$ and $\gamma_{m^{l-}-1}(h_{i_m})\neq0$ for some $l\in[1,\xi]$}

\vspace{2mm}

\nd
In this case, by $\gamma_{m^{l-}-1}(h_{i_m})
= \gamma_{m^{l-}}(h_{i_m})+2c_{m^{l-}}
\geq \gamma_{m^{l-}}(h_{i_m})=0$ and Lemma \ref{lem-a1} (2),
it holds
\begin{equation}\label{Amplus-pr3}
\gamma_{m^{l-}-1}(h_{i_m})=2.
\end{equation}
We take $l$ as the smallest one which satisfies (\ref{Amplus-pr3}),
that is, $\gamma_{m^{l'-}-1}(h_{i_m})=0$ if $0\leq l'<l$. 
Hence it follows $c_{m^{l-}}=\frac{\gamma_{m^{l-}-1}-\gamma_{m^{l-}}}{2}(h_{i_m})=1$
and $d_{m^{l-}}=\frac{\gamma_{m^{l-}-1}+\gamma_{m^{l-}}}{2}(h_{i_m})=1$.
By a similar argument to (\ref{a6-pr2}), we have
\[
\gamma_{m^{l-}}(h_{i_m})=\gamma_{m^{(l-1)-}-1}(h_{i_m})+\langle h_{i_m}, \sum_{p=m^-+1}^{m^{(l-1)-}-1} c_p\alpha_{i_p} \rangle.
\]
Thus, for any $p\in[m^{l-}+1,m^{(l-1)-}-1]$ such that $a_{i_m,i_p}<0$ we have $c_p=0$.
Since $\pi$ is an $\mathbf{i}$-trail, it holds $e_{i_1}^{c_1}\cdots e_{i_N}^{c_N}v_{-s_i\Lambda_i}\neq0$.
Hence, the pre-$\mathbf{i}$-trail $\pi'=(\gamma_0',\cdots,\gamma_N')$ from $-w_0\Lambda_i$ to $-s_i\Lambda_i$
with integers $\{c_s'\}_{s=1}^N$ such that
$\gamma_{s-1}'-\gamma_s'=c_s'\alpha_{i_s}$
and
\[
c_{s}'=c_s\quad (s\neq m^{l-}, m^{(l-1)-}),
\]
\[
c_{m^{l-}}'=c_{m^{l-}}-1=0,\quad c_{m^{(l-1)-}}'=c_{m^{(l-1)-}}+1
\]
is an $\mathbf{i}$-trail by $0\neq e_{i_1}^{c_1}\cdots e_{i_N}^{c_N}v_{-s_i\Lambda_i}=e_{i_1}^{c_1'}\cdots e_{i_N}^{c_N'}v_{-s_i\Lambda_i}$.
Let $M'$ be the monomial in $\Delta_{w_0\Lambda_i,s_i\Lambda_i}\circ \theta^-_{\mathbf{i}}(t_1,\cdots,t_N)$
corresponding to $\pi'$. By Proposition \ref{trail-prop},
we have $M'=M\cdot A_{m^{l-}}$.
Let $j=m^{l-}$, $M'=MA_{j}=\prod^N_{r=1}t_r^{d_r'}$ and $\{b_{r}'\}_{r\in[1,N]}$ be integers 
in Definition \ref{bint} determined from $M'$.
Then it holds $d_j'=d_{m^{l-}}+1=2>0$, $b_{j^+}'=c_{j^+}'=c_{m^{(l-1)-}}'>0$ and
\[
d_{j^+}'=d_{m^{(l-1)-}}+1=
\frac{\gamma_{m^{(l-1)-}-1}+\gamma_{m^{(l-1)-}}}{2}(h_{i_m})+1=
\begin{cases}
0 & \text{if }l=1,\\
1 & \text{if }l>1
\end{cases}
\]
so that (a) in our claim is satisfied.

\vspace{2mm}

\nd
\underline{Case 1-2. $\xi=0$ or $\gamma_{m^{l-}-1}(h_{i_m})=0$ for all $l\in[1,\xi]$}

\vspace{2mm}

One can prove
\[
\gamma_{m^{l-}}(h_{i_m})=\gamma_{m^{(l-1)-}-1}(h_{i_m})+\langle h_{i_m}, \sum_{p=m^{l-}+1}^{m^{(l-1)-}-1} c_p\alpha_{i_p} \rangle \quad (l\in[1,\xi]).
\]
Since we assumed $\gamma_{m^{l-}-1}(h_{i_m})=\gamma_{m^{l-}}(h_{i_m})=0$ for $l\in[1,\xi]$, it follows
\begin{equation}\label{Amplus-pr6}
\langle h_{i_m}, \sum_{p=m^{l-}+1}^{m^{(l-1)-}-1} c_p\alpha_{i_p}\rangle=0 \quad (l\in[1,\xi]).
\end{equation}
In the case $\xi=0$, we do not need this argument.
Considering
\[
0\leq-w_0\Lambda_i(h_{i_m})=
\gamma_{0}(h_{i_m})=\gamma_{m^{\xi-}-1}(h_{i_m})+\langle h_{i_m}, \sum_{p=1}^{m^{\xi-}-1} c_p\alpha_{i_p} \rangle
\]
and $\gamma_{m^{\xi-}-1}(h_{i_m})=0$, it holds
\begin{equation}\label{Amplus-pr7}
\langle h_{i_m}, \sum_{p=1}^{m^{\xi-}-1} c_p\alpha_{i_p} \rangle=0
\end{equation}
so that $-w_0\Lambda_i(h_{i_m})=0$.
Thus, by (\ref{Amplus-pr6}) and (\ref{Amplus-pr7}),
for any $p\in[1,m-1]$ such that $a_{i_m,i_p}<0$ it follows $c_p=0$.
Hence, there exists $c\in\mathbb{C}^{\times}$ such that
\begin{eqnarray*}
cv_{-w_0\Lambda_i}
&=&e_{i_1}^{c_1}\cdots e_{i_m}^{c_m}\cdots e_{i_N}^{c_N}v_{-s_i\Lambda_i}\\
&=&e_{i_m}^{c_m}e_{i_1}^{c_1}\cdots e_{i_{m-1}}^{c_{m-1}}e_{i_{m+1}}^{c_{m+1}}\cdots e_{i_N}^{c_N}v_{-s_i\Lambda_i}\\
&=&e_{i_m}^{1}e_{i_1}^{c_1}\cdots e_{i_{m-1}}^{c_{m-1}}e_{i_{m+1}}^{c_{m+1}}\cdots e_{i_N}^{c_N}v_{-s_i\Lambda_i}.
\end{eqnarray*}
The above description implies $f_{i_m}v_{-w_0\Lambda_i}\neq0$, which contradicts $-w_0\Lambda_i(h_{i_m})=0$.

\vspace{2mm}

\nd
\underline{Case 2. $\gamma_{m^{l-}}(h_{i_m})\neq0$ for some $l\in[1,\xi]$}

\vspace{2mm}

In this case, we take the smallest $l\in[1,\xi]$ such that $\gamma_{m^{l-}}(h_{i_m})\neq0$.
First, we suppose that there exists some $1\leq p<l$ such that $\gamma_{m^{p-}-1}(h_{i_m})\neq 0$.
We take $p$ as the smallest one satisfying $1\leq p<l$ and $\gamma_{m^{p-}-1}(h_{i_m})\neq 0$.
By $\gamma_{m^{p-}}(h_{i_m})=0$ and 
$\gamma_{m^{p-}}(h_{i_m})
\leq \gamma_{m^{p-}}(h_{i_m})+c_{m^{p-}}\alpha_{i_m}(h_{i_m})
= \gamma_{m^{p-}-1}(h_{i_m})$, it holds $\gamma_{m^{p-}-1}(h_{i_m})=2$, $c_{m^{p-}}=1$
and $d_{m^{p-}}=\frac{\gamma_{m^{p-}-1}(h_{i_m})+\gamma_{m^{p-}}(h_{i_m})}{2}=1$.
One can also verify that
\[
d_{m^{(p-1)-}}=\frac{\gamma_{m^{(p-1)-}-1}(h_{i_m})+\gamma_{m^{(p-1)-}}(h_{i_m})}{2}=
\begin{cases}
-1 & \text{if }p=1,\\
0 & \text{if }p>1,
\end{cases}\]
\[
c_{m^{(p-1)-}}=\frac{\gamma_{m^{(p-1)-}-1}(h_{i_m})-\gamma_{m^{(p-1)-}}(h_{i_m})}{2}=
\begin{cases}
1 & \text{if }p=1,\\
0 & \text{if }p>1.
\end{cases}
\]
It follows from
\[
\gamma_{m^{p-}}(h_{i_m})=\gamma_{m^{(p-1)-}-1}(h_{i_m})+\langle h_{i_m}, \sum_{t=m^{p-}+1}^{m^{(p-1)-}-1} c_t\alpha_{i_t} \rangle
\]
and $\gamma_{m^{p-}}(h_{i_m})=\gamma_{m^{(p-1)-}-1}(h_{i_m})=0$ that $c_t=0$ for $t\in[m^{p-}+1,m^{(p-1)-}-1]$ such that
$a_{i_m,i_t}<0$. Just as in Case 1-1, we see that $M':=M\cdot A_{m^{p-}}$ is a monomial
 in $\Delta_{w_0\Lambda_i,s_i\Lambda_i}\circ \theta^-_{\mathbf{i}}(t_1,\cdots,t_N)$.
Let $j=m^{p-}$, $M':=MA_{j}=\prod^N_{r=1}t_r^{d_r'}$ and $\{b_{r}'\}_{r\in[1,N]}$ be integers 
in Definition \ref{bint} determined from $M'$. By a similar way to Case 1-1,
it holds $d_j'=d_{m^{p-}}+1=2>0$, $b_{j^+}'=b_{j^+}+1>0$ and
\[
d_{j^+}'=d_{m^{(p-1)-}}+1=
\begin{cases}
0 & \text{if }p=1,\\
1 & \text{if }p>1
\end{cases}
\]
so that (a) in our claim is satisfied.

Thus, we may assume $\gamma_{m^{p-}-1}(h_{i_m})= 0$ for all $p$ ($1\leq p<l$) so that
\begin{equation}\label{Amplus-pr8}
c_{m^{p-}}=\frac{\gamma_{m^{p-}-1}(h_{i_m})-\gamma_{m^{p-}}(h_{i_m})}{2}=0 \quad (1\leq p<l),
\end{equation}
\begin{equation}\label{Amplus-pr8a}
d_{m^{p-}}=\frac{\gamma_{m^{p-}-1}(h_{i_m})+\gamma_{m^{p-}}(h_{i_m})}{2}=0 \quad (1\leq p<l),
\end{equation}
\begin{equation}\label{Amplus-pr9}
c_{m^{(l-1)-}}=\frac{\gamma_{m^{(l-1)-}-1}(h_{i_m})-\gamma_{m^{(l-1)-}}(h_{i_m})}{2}=
\begin{cases}
1 & \text{if }l=1,\\
0 & \text{if }l>1,
\end{cases}
\end{equation}
\begin{equation}\label{Amplus-pr10}
d_{m^{(l-1)-}}=\frac{\gamma_{m^{(l-1)-}-1}(h_{i_m})+\gamma_{m^{(l-1)-}}(h_{i_m})}{2}=
\begin{cases}
-1 & \text{if }l=1,\\
0 & \text{if }l>1.
\end{cases}
\end{equation}

\vspace{2mm}

\nd
\underline{Case 2-1. $\gamma_{m^{l-}}(h_{i_m})>0$}

\vspace{2mm}

By $\gamma_{m-1}(h_{i_m})=0$
and $\gamma_{m-1}(h_{i_m})\geq\gamma_{m^-}(h_{i_m})$,
it holds $l\geq2$.
Combining $\gamma_{m^{(l-1)-}-1}(h_{i_m})\geq\gamma_{m^{l-}}(h_{i_m})>0$
and $\gamma_{m^{(l-1)-}}(h_{i_m})=0$, it follows
$\gamma_{m^{(l-1)-}-1}(h_{i_m})=2$,
$c_{m^{(l-1)-}}=
\frac{\gamma_{m^{(l-1)-}-1}(h_{i_m})-\gamma_{m^{(l-1)-}}(h_{i_m})}{2}=1$ and
$d_{m^{(l-1)-}}=
\frac{\gamma_{m^{(l-1)-}-1}(h_{i_m})+\gamma_{m^{(l-1)-}}(h_{i_m})}{2}=1$. Note that
by
\[
\gamma_{m^{(l-1)-}}(h_{i_m})=\gamma_{m^{(l-2)-}-1}(h_{i_m})+\langle h_{i_m}, \sum_{t=m^{(l-1)-}+1}^{m^{(l-2)-}-1} c_t\alpha_{i_t} \rangle
\]
and $\gamma_{m^{(l-1)-}}(h_{i_m})=\gamma_{m^{(l-2)-}-1}(h_{i_m})=0$,
for any $t\in[m^{(l-1)-},m^{(l-2)-}]$ such that $a_{i_m,i_t}<0$, it holds $c_t=0$.
One obtains $M\cdot A_{m^{(l-1)-}}$ is a monomial
 in $\Delta_{w_0\Lambda_i,s_i\Lambda_i}\circ \theta^-_{\mathbf{i}}(t_1,\cdots,t_N)$
by a similar way to the end of Case 1-1.
We see that
\[
d_{m^{(l-2)-}}=
\frac{\gamma_{m^{(l-2)-}-1}(h_{i_m})+\gamma_{m^{(l-2)-}}(h_{i_m})}{2}=
\begin{cases}
-1 & \text{if }l=2,\\
0 & \text{if }l>2.
\end{cases}
\]
Let $j=m^{(l-1)-}$, $M':=MA_{j}=\prod^N_{r=1}t_r^{d_r'}$ and $\{b_{r}'\}_{r\in[1,N]}$ be integers 
in Definition \ref{bint} determined from $M'$. By a similar way to Case 1-1,
it holds $d_j'=d_{m^{(l-1)-}}+1=2>0$, $b_{j^+}'=b_{j^+}+1>0$ and
\[
d_{j^+}'=d_{m^{(l-2)-}}+1=
\begin{cases}
0 & \text{if }l=2,\\
1 & \text{if }l>2
\end{cases}
\]
so that (a) in our claim is satisfied.

\vspace{2mm}

\nd
\underline{Case 2-2. $\gamma_{m^{l-}}(h_{i_m})<0$}

\vspace{2mm}

Lemma \ref{lem-a1} (2) means $\gamma_{m^{l-}}(h_{i_m})=-2$ or $-1$.

\vspace{2mm}

\nd
\underline{Case 2-2-1. $\gamma_{m^{l-}}(h_{i_m})=-2$}

\vspace{2mm}

In the case $\gamma_{m^{l-}}(h_{i_m})=-2$, by (\ref{Amplus-pr1}), we see that
$0\leq d_{m^{l-}}=\frac{\gamma_{m^{l-}-1}(h_{i_m})+\gamma_{m^{l-}}(h_{i_m})}{2}$
and
$\gamma_{m^{l-}-1}(h_{i_m})=2$ so that $c_{m^{l-}}=\frac{\gamma_{m^{l-}-1}(h_{i_m})-\gamma_{m^{l-}}(h_{i_m})}{2}=2$
and $d_{m^{l-}}=0$.
Since 
\begin{equation}\label{221-1}
\gamma_{m^{l-}}(h_{i_m})=\gamma_{m^{(l-1)-}-1}(h_{i_m})+\langle h_{i_m}, \sum_{s=m^{l-}+1}^{m^{(l-1)-}-1} c_s\alpha_{i_s} \rangle
\end{equation}
and
$\gamma_{m^{(l-1)-}-1}(h_{i_m})=0$, there are three patterns:
\begin{enumerate}
\item[$(1)$] there exist $l_1$, $l_2\in[m^{l-},m^{(l-1)-}]$ such that $l_1< l_2$, $a_{i_{m},i_{l_1}}=a_{i_{m},i_{l_2}}=-1$
and $c_{l_1}=c_{l_2}=1$,
\item[$(2)$] there exists $l_1\in[m^{l-},m^{(l-1)-}]$ such that $a_{i_{m},i_{l_1}}=-1$ and $c_{l_1}=2$,
\item[$(3)$] there exists $l_1\in[m^{l-},m^{(l-1)-}]$ such that $a_{i_{m},i_{l_1}}=-2$ and $c_{l_1}=1$.
\end{enumerate}
In the case (1), by (\ref{221-1}), for $s\in[m^{l-}+1,m^{(l-1)-}-1]\setminus\{l_1,l_2\}$ such that $a_{i_m,i_s}<0$, we have $c_s=0$.
Considering
\[
\gamma_{m^{r-}}(h_{i_m})=\gamma_{m^{(r-1)-}-1}(h_{i_m})+\langle h_{i_m}, \sum_{s=m^{r-}+1}^{m^{(r-1)-}-1} c_s\alpha_{i_s} \rangle\quad
(r\in[1,l-1])
\]
and $\gamma_{m^{r-}}(h_{i_m})=\gamma_{m^{(r-1)-}-1}(h_{i_m})=0$ ($r\in[1,l-1]$),
for $s\in[m^{(l-1)-}+1,m-1]\setminus\{l_1,l_2\}$ such that $a_{i_m,i_s}<0$, it holds $c_s=0$. Thus, it holds
\begin{equation}\label{Amplus-pr10-2}
c_s=0 \quad (s\in[m^{l-}+1,m-1]\setminus\{l_1,l_2\}\text{ such that }a_{i_m,i_s}<0).
\end{equation}
Since $\pi$ is an $\mathbf{i}$-trail, it holds
\begin{eqnarray}
0&\neq&e_{i_1}^{c_1}\cdots e_{i_{m^{l-}}}^{c_{m^{l-}}}\cdots e_{i_{l_1}}^{c_{l_1}}\cdots e_{i_{l_2}}^{c_{l_2}}\cdots e_{i_{m^{(l-1)-}}}^{c_{m^{(l-1)-}}}
\cdots e_{i_N}^{c_N}v_{-s_i\Lambda_i}\nonumber\\
&=&e_{i_1}^{c_1}\cdots e_{i_{m^{l-}}}^{2}\cdots e_{i_{l_1}}^{1}\cdots e_{i_{l_2}}^{c_{l_2}}\cdots e_{i_{m^{(l-1)-}}}^{c_{m^{(l-1)-}}}
\cdots e_{i_N}^{c_N}v_{-s_i\Lambda_i}. \label{Amplus-pr11}
\end{eqnarray}
As assumed
\[
-2=\gamma_{m^{l-}}(h_{i_m})
=\langle
h_{i_m},
{\rm wt}(e_{i_{m^{l-}+1}}^{c_{m^{l-}+1}}\cdots e_{i_{l_1}}^{c_{l_1}}\cdots e_{i_N}^{c_N}v_{-s_i\Lambda_i})
\rangle,
\]
in conjunction with (\ref{Amplus-pr10-2}), one obtains
$\langle
h_{i_m},
{\rm wt}(e_{i_{l_1+1}}^{c_{l_1+1}}\cdots e_{i_N}^{c_N}v_{-s_i\Lambda_i})
\rangle=-1$.
Hence, $e_{i_{m^{l-}}}^2e_{i_{l_1+1}}^{c_{l_1+1}}\cdots e_{i_N}^{c_N}v_{-s_i\Lambda_i}=0$, otherwise, we get
$\langle
h_{i_m},
{\rm wt}(e_{i_{m^{l-}}}^2e_{i_{l_1+1}}^{c_{l_1+1}}\cdots e_{i_N}^{c_N}v_{-s_i\Lambda_i})
\rangle=3$, which contradicts Lemma \ref{lem-a1} (1). In particular,
\begin{equation}\label{Amplus-pr12}
e_{i_1}^{c_1}
\cdots
e_{i_{m^{l-}-1}}^{c_{m^{l-}-1}}
e_{i_{m^{l-}}}^0
e_{i_{m^{l-}+1}}^{c_{m^{l-}+1}}
\cdots 
e_{i_{l_1-1}}^{c_{l_1-1}}
e_{i_{l_1}}^1
e_{i_{m^{l-}}}^2e_{i_{l_1+1}}^{c_{l_1+1}}\cdots e_{i_N}^{c_N}v_{-s_i\Lambda_i}=0.
\end{equation}
Considering the Serre relation
$e_{i_{m^{l-}}}^{2} e_{i_{l_1}}^{1}-2e_{i_{m^{l-}}}^1e_{i_{l_1}}^{1}e_{i_{m^{l-}}}^1+e_{i_{l_1}}^1e_{i_{m^{l-}}}^2=0$
with (\ref{Amplus-pr10-2}), (\ref{Amplus-pr11}) and (\ref{Amplus-pr12}), we see that
\[
e_{i_1}^{c_1}
\cdots
e_{i_{m^{l-}-1}}^{c_{m^{l-}-1}}
e_{i_{m^{l-}}}^1
e_{i_{m^{l-}+1}}^{c_{m^{l-}+1}}
\cdots 
e_{i_{l_1-1}}^{c_{l_1-1}}
e_{i_{l_1}}^1
e_{i_{m^{l-}}}^1e_{i_{l_1+1}}^{c_{l_1+1}}\cdots e_{i_N}^{c_N}v_{-s_i\Lambda_i}\neq0.
\]
By $c_{m^{l-}}=2$, $c_{l_1}=1$, $c_{l_2}=1$ and $c_m=1$, the above vector equals
\begin{equation}\label{Amplus-pr13}
e_{i_1}^{c_1}
\cdots
e_{i_{m^{l-}}}^{c_{m^{l-}}-1}
\cdots e_{i_{l_1}}^{c_{l_1}}
e_{i_{m^{l-}}}^1e_{i_{l_1+1}}^{c_{l_1+1}}\cdots 
e_{i_{l_2-1}}^{c_{l_2-1}}e_{i_{l_2}}^1e_{i_{l_2+1}}^{c_{l_2+1}}
\cdots 
e_{i_{m-1}}^{c_{m-1}}
e_{i_m}^1e_{i_{m+1}}^{c_{m+1}} \cdots 
e_{i_N}^{c_N}v_{-s_i\Lambda_i}\neq0.
\end{equation}
One can verify by
(\ref{Amplus-pr2}) and (\ref{Amplus-pr8}) that
\[
e_{i_{l_2}}^1e_{i_{l_2+1}}^{c_{l_2+1}}
\cdots 
e_{i_{m-1}}^{c_{m-1}}
e_{i_m}^0e_{i_{m+1}}^{c_{m+1}} \cdots 
e_{i_N}^{c_N}v_{-s_i\Lambda_i}=0,
\]
otherwise, $\langle h_{i_m}, {\rm wt}(e_{i_{l_2}}^1e_{i_{l_2+1}}^{c_{l_2+1}}
\cdots 
e_{i_{m-1}}^{c_{m-1}}
e_{i_m}^0e_{i_{m+1}}^{c_{m+1}} \cdots 
e_{i_N}^{c_N}v_{-s_i\Lambda_i}) \rangle=-3$, which is absurd. In particular,
\begin{equation}\label{Amplus-pr14}
e_{i_1}^{c_1}
\cdots
e_{i_{m^{l-}}}^{c_{m^{l-}}-1}
\cdots e_{i_{l_1}}^{c_{l_1}}
e_{i_{m^{l-}}}^2e_{i_{l_1+1}}^{c_{l_1+1}}\cdots 
e_{i_{l_2-1}}^{c_{l_2-1}}e_{i_{l_2}}^1e_{i_{l_2+1}}^{c_{l_2+1}}
\cdots 
e_{i_{m-1}}^{c_{m-1}}
e_{i_m}^0e_{i_{m+1}}^{c_{m+1}} \cdots 
e_{i_N}^{c_N}v_{-s_i\Lambda_i}=0.
\end{equation}
Using the Serre relation $e_{i_{m^{l-}}}^2e_{i_{l_2}}^1-2e_{i_{m^{l-}}}^1e_{i_{l_2}}^1e_{i_m}^1+e_{i_{l_2}}^1e_{i_{m^{(l-1)-}}}e_{i_m}=0$
with (\ref{Amplus-pr10-2}), (\ref{Amplus-pr13}) and (\ref{Amplus-pr14}), it can be shown that
\begin{eqnarray*}
0&\neq&e_{i_1}^{c_1}
\cdots
e_{i_{m^{l-}}}^{c_{m^{l-}}-1}
\cdots e_{i_{l_1}}^{c_{l_1}}
e_{i_{l_1+1}}^{c_{l_1+1}}\cdots 
e_{i_{l_2-1}}^{c_{l_2-1}}e_{i_{l_2}}^1e_{i_{l_2+1}}^{c_{l_2+1}}
\cdots \\
& &
\qquad e_{i_{m^{(l-1)-}-1}}^{c_{m^{(l-1)-}-1}}
e_{i_{m^{(l-1)-}}}^{c_{m^{(l-1)-}}+1}
e_{i_{m^{(l-1)-}+1}}^{c_{m^{(l-1)-}+1}}
\cdots
e_{i_{m-1}}^{c_{m-1}}
e_{i_m}^1e_{i_{m+1}}^{c_{m+1}} \cdots 
e_{i_N}^{c_N}v_{-s_i\Lambda_i}\\
&=&
e_{i_1}^{c_1}
\cdots
e_{i_{m^{l-}}}^{c_{m^{l-}}-1}
\cdots e_{i_{l_1}}^{c_{l_1}}
e_{i_{l_1+1}}^{c_{l_1+1}}\cdots 
e_{i_{l_2-1}}^{c_{l_2-1}}e_{i_{l_2}}^{c_{l_2}}e_{i_{l_2+1}}^{c_{l_2+1}}
\cdots \\
& &\qquad
e_{i_{m^{(l-1)-}-1}}^{c_{m^{(l-1)-}-1}}
e_{i_{m^{(l-1)-}}}^{c_{m^{(l-1)-}}+1}
e_{i_{m^{(l-1)-}+1}}^{c_{m^{(l-1)-}+1}}
\cdots
e_{i_{m-1}}^{c_{m-1}}
e_{i_m}^{c_m}e_{i_{m+1}}^{c_{m+1}} \cdots 
e_{i_N}^{c_N}v_{-s_i\Lambda_i}.
\end{eqnarray*}
In the 4th line of above vector, we understand
$e_{i_{m^{(l-1)-}}}^{c_{m^{(l-1)-}}+1}
e_{i_{m^{(l-1)-}+1}}^{c_{m^{(l-1)-}+1}}
\cdots
e_{i_{m-1}}^{c_{m-1}}
e_{i_m}^{c_m}=e_{i_m}^{c_m+1}$ in case of $l=1$.
Consequently,
the pre-$\mathbf{i}$-trail
$\pi'=(\gamma_0',\gamma_1',\cdots,\gamma_N')$
from $-w_0\Lambda_i$ to $-s_i\Lambda_i$ with integers $c_r'$ $(r=1,2,\cdots,N)$ such that
$\gamma_{r-1}'-\gamma_r'=c_r'\alpha_{i_r}$ and
\[
c_r'=c_{r}\quad \text{for}\ r\in[1,N]\setminus \{m^{l-},m^{(l-1)-}\},
\]
\[
c_{m^{l-}}'=c_{m^{l-}}-1=1,\quad c_{m^{(l-1)-}}'=c_{m^{(l-1)-}}+1
\]
is an $\mathbf{i}$-trail. Let $M'$
be the monomial in $\Delta_{w_0\Lambda_i,s_i\Lambda_i}\circ \theta^-_{\mathbf{i}}(t_1,\cdots,t_N)$
corresponding to $\pi'$. By Proposition \ref{trail-prop},
we get $M'=M\cdot A_{m^{l-}}$.
Let $j:=m^{l-}$, $M'=MA_{j}=\prod^N_{r=1}t_r^{d_r'}$ and $\{b_{r}'\}_{r\in[1,N]}$ be integers 
in Definition \ref{bint} determined from $M'$. We see that $d_{j}'=d_j+1=d_{m^{l-}}+1=1>0$ and
 $b_{j^+}'=c_{m^{(l-1)-}}'=c_{m^{(l-1)-}}+1>0$.
Note that (\ref{Amplus-pr8}) and (\ref{Amplus-pr8a}) imply $b_{j^{r+}}'=c_{j^{r+}}'=c_{j^{r+}}=0$ and
$d_{j^{r+}}'=d_{j^{r+}}=0$ for $r\in[2,l-1]$. It follows by (\ref{Amplus-pr10}) that
\begin{equation}\label{Amplus-pr16}
d_{j^+}'=d_{j^+}+1
=d_{m^{(l-1)-}}+1
=
\begin{cases}
0 & \text{if }l=1,\\
1 & \text{if }l>1.
\end{cases}
\end{equation}
Since we know $d_{j^{l+}}=d_m=-1$ and $b_{j^{l+}}=c_{j^{l+}}=c_m=1$, it holds
\begin{equation}\label{Amplus-pr17}
d_{j^{l+}}'
=\begin{cases}
d_{j^{l+}}+1=0 & \text{if }l=1,\\
d_{j^{l+}}=-1 & \text{if }l>1,
\end{cases}\qquad
b_{j^{l+}}'
=\begin{cases}
b_{j^{l+}}+1=2 & \text{if }l=1,\\
b_{j^{l+}}=1 & \text{if }l>1.
\end{cases}
\end{equation}
In this way, if $l=1$ then (a) in our claim holds and if $l>1$ then putting $p:=l$, the condition (b) in our claim holds.
One can similarly show our claim in the case (2).

Next, we assume (3). 
By the same argument as in (\ref{Amplus-pr10-2}),
for $s\in[m^{l-}+1,m-1]\setminus\{l_1\}$ such that $a_{i_m,i_s}<0$, it holds 
\begin{equation}\label{Amplus-pr15}
c_s=0.
\end{equation}
Since $\pi$ is an $\mathbf{i}$-trail,
\begin{eqnarray*}
0&\neq & e_{i_1}^{c_1}\cdots e_{i_{m^{l-}}}^{c_{m^{l-}}}\cdots
e_{i_{l_1}}^{c_{l_1}}\cdots e_{i_m}^{c_m}\cdots e_{i_N}^{c_N}
v_{-s_i\Lambda_i}\\
&=& e_{i_1}^{c_1}\cdots e_{i_{m^{l-}}}^{2}\cdots
e_{i_{l_1}}^{1}\cdots e_{i_m}^{1}\cdots e_{i_N}^{c_N}
v_{-s_i\Lambda_i}.
\end{eqnarray*}
Thus, considering the Serre relation $e_{i_{m^{l-}}}^{2}e_{i_{l_1}}^{1}e_{i_m}^{1}
=e_{i_{m^{l-}}}^{1}e_{i_{l_1}}^{1}e_{i_{m^{(l-1)-}}}^{1}e_{i_m}^{1}$ on $V(-w_0\Lambda_i)$ and (\ref{Amplus-pr15}), the following holds:
\begin{eqnarray*}
0&\neq& 
 e_{i_1}^{c_1}\cdots e_{i_{m^{l-}-1}}^{c_{m^{l-}-1}}
e_{i_{m^{l-}}}^{1}
e_{i_{m^{l-}+1}}^{c_{m^{l-}+1}}
\cdots
e_{i_{l_1-1}}^{c_{l_1-1}}
e_{i_{l_1}}^{1}
e_{i_{l_1+1}}^{c_{l_1+1}}
\cdots
e_{i_{m^{(l-1)-}}}^{c_{m^{(l-1)-}}+1}\cdots
 e_{i_m}^{1}\cdots e_{i_N}^{c_N}v_{-s_i\Lambda_i}\\
&=&
e_{i_1}^{c_1}\cdots e_{i_{m^{l-}-1}}^{c_{m^{l-}-1}}
e_{i_{m^{l-}}}^{c_{m^{l-}}-1}
e_{i_{m^{l-}+1}}^{c_{m^{l-}+1}}
\cdots
e_{i_{l_1-1}}^{c_{l_1-1}}
e_{i_{l_1}}^{c_{l_1}}
e_{i_{l_1+1}}^{c_{l_1+1}}
\cdots
e_{i_{m^{(l-1)-}}}^{c_{m^{(l-1)-}}+1}\cdots
 e_{i_m}^{c_{m}}\cdots e_{i_N}^{c_N}v_{-s_i\Lambda_i}.
\end{eqnarray*}
Here, in the case $l=1$, we understand $e_{i_{m^{(l-1)-}}}^{c_{m^{(l-1)-}}+1}\cdots
 e_{i_m}^{c_{m}}=e_{i_m}^{c_m+1}$.
By the same argument as in the end of case (1), we see that $M':=M\cdot A_{m^{l-}}$ is a monomial
in $\Delta_{w_0\Lambda_i,s_i\Lambda_i}\circ \theta^-_{\mathbf{i}}(t_1,\cdots,t_N)$
and defining
$j:=m^{l-}$, $M'=MA_{j}=\prod^N_{r=1}t_r^{d_r'}$ and $\{b_{r}'\}_{r\in[1,N]}$ as integers 
in Definition \ref{bint} determined from $M'$, we also see that $d_{j}'>0$,
 $b_{j^+}'>0$, $b_{j^{r+}}'=d_{j^{r+}}'=0$ for $r\in[2,l-1]$
 and (\ref{Amplus-pr16}), (\ref{Amplus-pr17}) hold.
Therefore, if $l=1$ then (a) in our claim holds and
if $l>1$ then setting $p:=l$, the condition (b) in our claim holds.

\vspace{2mm}

\nd
\underline{Case 2-2-2. $\gamma_{m^{l-}}(h_{i_m})=-1$}

\vspace{2mm}

In the case $\gamma_{m^{l-}}(h_{i_m})=-1$, by (\ref{Amplus-pr1}), we see that
$0\leq d_{m^{l-}}=\frac{\gamma_{m^{l-}-1}(h_{i_m})+\gamma_{m^{l-}}(h_{i_m})}{2}$
and
$\gamma_{m^{l-}-1}(h_{i_m})=1$ so that $c_{m^{l-}}=\frac{\gamma_{m^{l-}-1}(h_{i_m})-\gamma_{m^{l-}}(h_{i_m})}{2}=1$ and $d_{m^{l-}}=0$.
Considering $\gamma_{m^{(l-1)-}-1}(h_{i_m})=0$ and (\ref{221-1}), we see that
there exists $l_1\in[m^{l-}+1,m^{(l-1)-}-1]$ such that $a_{i_{m},i_{l_1}}=-1$ and $c_{l_1}=1$.
For $l'\in[m^{l-}+1,m^{(l-1)-}-1]\setminus\{l_1\}$ such that $a_{i_{m},i_{l'}}<0$, it holds $c_{l'}=0$.
By (\ref{Amplus-pr8}) and the same argument as in (\ref{Amplus-pr10-2}),
we see that
for $s\in[m^{l-}+1,m-1]\setminus\{l_1\}$ such that $a_{i_m,i_s}\neq0$, it holds 
\begin{equation}\label{Amplus-pr18}
c_s=0.
\end{equation}
The definition of $\mathbf{i}$-trail implies
\begin{eqnarray}
0&\neq&e_{i_1}^{c_1}\cdots e_{i_{m^{l-}}}^{c_{m^{l-}}}\cdots e_{i_{l_1}}^{c_{l_1}}
\cdots 
e_{i_m}^{c_m}\cdots
e_{i_N}^{c_N}v_{-s_i\Lambda_i}\nonumber\\
&=&e_{i_1}^{c_1}\cdots e_{i_{m^{l-}}}^{1}\cdots e_{i_{l_1}}^{1}
\cdots 
e_{i_m}^{1}\cdots
e_{i_N}^{c_N}v_{-s_i\Lambda_i}. \label{Amplus-pr19}
\end{eqnarray}
The assumption
\[
-2=\gamma_{m}(h_{i_m})
=\langle
h_{i_m},
{\rm wt}(e_{i_{m+1}}^{c_{m+1}}\cdots
e_{i_N}^{c_N}v_{-s_i\Lambda_i})
\rangle,
\]
and (\ref{Amplus-pr18}) yield
$e_{i_{l_1}}^{1} e_{i_{l_1+1}}^{c_{l_1+1}}
\cdots
e_{i_{m-1}}^{c_{m-1}}e_{i_m}^0
e_{i_{m+1}}^{c_{m+1}}\cdots
e_{i_N}^{c_N}v_{-s_i\Lambda_i}=0$. In particular,
\begin{equation}\label{Amplus-pr20}
e_{i_1}^{c_1}\cdots
e_{i_{m^{l-}-1}}^{c_{m^{l-}-1}}
e_{i_{m^{l-}}}^2
e_{i_{m^{l-}+1}}^{c_{m^{l-}+1}}
\cdots
e_{i_{l_1-1}}^{c_{l_1-1}}
e_{i_{l_1}}^{1} e_{i_{l_1+1}}^{c_{l_1+1}}
\cdots
e_{i_{m-1}}^{c_{m-1}}e_{i_m}^0
e_{i_{m+1}}^{c_{m+1}}\cdots
e_{i_N}^{c_N}v_{-s_i\Lambda_i}=0.
\end{equation}
It follows by 
the Serre relation 
$e_{i_{m^{l-}}}^{2}e_{i_{l_1}}^{1}-2e_{i_{m^{l-}}}^{1}e_{i_{l_1}}^{1}e_{i_m}^{1}+e_{i_{l_1}}^{1}e_{i_{m^{(l-1)-}}}e_{i_m}^{1}=0$, (\ref{Amplus-pr18})
(\ref{Amplus-pr19}) and (\ref{Amplus-pr20}) that
\begin{eqnarray*}
0&\neq&
e_{i_1}^{c_1}\cdots
e_{i_{m^{l-}-1}}^{c_{m^{l-}-1}}
e_{i_{m^{l-}}}^0
e_{i_{m^{l-}+1}}^{c_{m^{l-}+1}}
\cdots
e_{i_{l_1-1}}^{c_{l_1-1}}
e_{i_{l_1}}^{1} e_{i_{l_1+1}}^{c_{l_1+1}}
\cdots\\
& &\quad
e_{i_{m^{(l-1)-}}}^{c_{m^{(l-1)-}}+1}
\cdots
e_{i_{m-1}}^{c_{m-1}}e_{i_m}^1
e_{i_{m+1}}^{c_{m+1}}\cdots
e_{i_N}^{c_N}v_{-s_i\Lambda_i}\\
&=&
e_{i_1}^{c_1}\cdots
e_{i_{m^{l-}-1}}^{c_{m^{l-}-1}}
e_{i_{m^{l-}}}^{c_{m^{l-}}-1}
e_{i_{m^{l-}+1}}^{c_{m^{l-}+1}}
\cdots
e_{i_{l_1-1}}^{c_{l_1-1}}
e_{i_{l_1}}^{c_{l_1}} e_{i_{l_1+1}}^{c_{l_1+1}}
\cdots \\
& &\quad
e_{i_{m^{(l-1)-}}}^{c_{m^{(l-1)-}}+1}
\cdots
e_{i_{m-1}}^{c_{m-1}}e_{i_m}^{c_m}
e_{i_{m+1}}^{c_{m+1}}\cdots
e_{i_N}^{c_N}v_{-s_i\Lambda_i},
\end{eqnarray*}
where in the case $l=1$, we understand 
$e_{i_{m^{(l-1)-}}}^{c_{m^{(l-1)-}}+1}
\cdots e_{i_{m-1}}^{c_{m-1}}e_{i_m}^{c_m}=e_{i_m}^{c_m+1}$.
Hence,
the pre-$\mathbf{i}$-trail
$\pi'=(\gamma_0',\gamma_1',\cdots,\gamma_N')$
from $-w_0\Lambda_i$ to $-s_i\Lambda_i$ with integers $c_r'$ $(r=1,2,\cdots,N)$ such that
$\gamma_{r-1}'-\gamma_r'=c_r'\alpha_{i_r}$ and
\[
c_r'=c_{r}\quad \text{for}\ r\in[1,N]\setminus \{m^{l-},m^{(l-1)-}\},
\]
\[
c_{m^{l-}}'=c_{m^{l-}}-1=0,\quad c_{m^{(l-1)-}}'=c_{m^{(l-1)-}}+1
\]
is an $\mathbf{i}$-trail. Let $M'$
be the monomial in $\Delta_{w_0\Lambda_i,s_i\Lambda_i}\circ \theta^-_{\mathbf{i}}(t_1,\cdots,t_N)$
corresponding to $\pi'$. By Proposition \ref{trail-prop},
we get $M'=M\cdot A_{m^{l-}}$.
Let $j:=m^{l-}$, $M'=MA_{j}=\prod^N_{r=1}t_r^{d_r'}$ and $\{b_{r}'\}_{r\in[1,N]}$ be integers 
in Definition \ref{bint} determined from $M'$. We see that $d_{j}'=d_j+1=d_{m^{l-}}+1=1>0$ and
$b_{j^+}'=c_{m^{(l-1)-}}'=c_{m^{(l-1)-}}+1>0$.
The equations $b_{j^{r+}}'=c_{j^{r+}}'=c_{j^{r+}}=0$ and
$d_{j^{r+}}'=d_{j^{r+}}=0$ for $r\in[2,l-1]$ follow from (\ref{Amplus-pr8}) and (\ref{Amplus-pr8a}).
By (\ref{Amplus-pr10}),
\[
d_{j^+}'=d_{j^+}+1
=d_{m^{(l-1)-}}+1
=
\begin{cases}
0 & \text{if }l=1,\\
1 & \text{if }l>1.
\end{cases}
\]
Since we know $d_{j^{l+}}=d_m=-1$ and $b_{j^{l+}}=c_{j^{l+}}=c_m=1$, it holds
\[
d_{j^{l+}}'
=\begin{cases}
0 & \text{if }l=1,\\
-1 & \text{if }l>1,
\end{cases}\qquad
b_{j^{l+}}'
=\begin{cases}
2 & \text{if }l=1,\\
1 & \text{if }l>1.
\end{cases}
\]
In this way, if $l=1$ then (a) in our claim holds and if $l>1$ then putting $p:=l$, the condition (b) in our claim holds. \qed

\vspace{3mm}

\nd
[Proof of Theorem \ref{thm2}]

Let $\mathcal{M}$ be the set of monomials 
appearing in
$\Delta_{w_0\Lm_i,s_i\Lm_i}\circ \theta^-_{\mathbf{i}}(t_1,\cdots,t_N)$. For a directed graph $D$,
let $V(D)$ be the set of vertices in $D$. First, let us show
\begin{enumerate}
\item 
well-definedness of
each vector $(b_1,\cdots,b_N)$ due to Theorem \ref{thm2} associated to each monomial $M$ in $V(\overline{DG}_l)$,
that is, if there exists $M'$, $M''\in V(\overline{DG}_{l-1})$ associated with $(b_1',\cdots,b_N')$, $(b_1'',\cdots,b_N'')$,
respectively and $j'$, $j''\in[1,N]$ such that both pairs $(M',j')$ and $(M'',j'')$ satisfy the
condition (2) of Theorem \ref{thm2} and $M'\cdot A_{j'}^{-1}=M''\cdot A_{j''}^{-1}$ then
\[
(b_1',\cdots,b_{j'}'+1,\cdots,b_{j'^+}'-1,\cdots,b_N')
=(b_1'',\cdots,b_{j''}''+1,\cdots,b_{j''^+}''-1,\cdots,b_N''),
\]
\item each vector $(b_1,\cdots,b_N)$ associated to each monomial $M$ in $V(\overline{DG}_l)$ coincides
with the $N$-tuple of integers in Definition \ref{bint} determined from $M$,
\item $V(\overline{DG}_l)\subset \mathcal{M}$,
\end{enumerate}
for $l=0,1,\cdots,r$ if the algorithm will stop at some step $r$ and for $l\in\mathbb{Z}_{\geq0}$ if the algorithm
will not stop (after, we will prove the algorithm will stop at some step $r$).
We will show the above 1, 2 and 3
by the induction on $l$ simultaneously. If $l=0$ then the above 1 and 2 are clear by
(1) of Theorem \ref{thm2}. 
It follows from Proposition \ref{high-prop} that
$V(\overline{DG}_0)=\{t_k\}\subset\mathcal{M}$ so that 3 holds.
Thus, let us assume 1, 2 and 3 hold for some $l\geq0$. We may assume that
$\overline{DG}_{l}$ has a sink $M$ associated with the integer vector $(b_1,\cdots,b_N)$ of Definition \ref{bint} determined from $M$
such that $M$ and some $j\in[1,N]$
satisfy the condition (2) in Theorem \ref{thm2}, otherwise, the algorithm stops at the step $l$ so that $l=r$.
By the induction assumption $V(\overline{DG}_l)\subset \mathcal{M}$,
it holds $M\in \mathcal{M}$.
If the sink $M$ and $j$ satisfy (a) of (2), that is, $d_j>0$, $b_{j^+}>0$ and $d_{j^+}<d_j$ then
in case of $d_j=2$, it follows from Lemma \ref{lem-a2} that $M\cdot A_{j}^{-1}\in\mathcal{M}$.
In case of $d_j=1$, Lemma \ref{lem-a3} implies $M\cdot A_{j}^{-1}\in\mathcal{M}$.
If the sink $M$ and $j$ satisfy (b) of (2), that is, $d_j>0$, $b_{j^+}>0$, $d_{j^+}=d_j$ and there exists
$p\in\mathbb{Z}_{\geq2}$ such that $d_{j^{m+}}=b_{j^{m+}}=0$ ($m=2,3,\cdots,p-1$) and
$d_{j^{p+}}=-1$, $b_{j^{p+}}=1$. Note that it holds $d_j\neq2$, otherwise, $d_j=d_{j^+}=2$, which contradicts Lemma \ref{lem-a2}.
Thus, we get $d_j=1$. In case of $b_j=0$ (resp. $b_j>0$), Lemma \ref{lem-a4} (resp. Lemma \ref{lem-a5}) means
$M\cdot A_{j}^{-1}\in\mathcal{M}$. In this way, one can prove
$V(\overline{DG}_{l+1})\subset \mathcal{M}$.
Defining a vector $(b_1',\cdots,b_N')$ by
$b'_s=b_s$ $(s\in[1,N]\setminus\{j,j^+\})$, $b_j'=b_j+1$, $b'_{j^+}=b_{j^+}-1$,
it coincides with $N$-tuple of integers of Definition \ref{bint}
determined from $M\cdot A_{j}^{-1}$ by Lemma \ref{lem-4}.
In particular, the associated vector $(b_1',\cdots,b_N')$ to $M\cdot A_{j}^{-1}$ is well-defined. In this way,
1, 2 and 3 hold for all $l$.

Next, let us prove the algorithm in Theorem \ref{thm2} will stop at some step $r\in\mathbb{Z}_{\geq0}$.
For each monomial $M\in \mathcal{M}$ associated with integers $(b_1,\cdots,b_N)$, we set
\[
L(M):=\sum_{t=1}^N t b_t. 
\]
One can easily verify that if $j\in[1,N]$, $j^+\leq N$ and $M\cdot A_j^{-1}\in\mathcal{M}$ then
\[
L(M\cdot A_j^{-1})=L(M)+j-j^+ <L(M)
\]
by Lemma \ref{lem-4}.
If the algorithm will not stop then there exists an infinite sequence
\[
t_k\rightarrow t_k A_{j_1}^{-1}\rightarrow t_k A_{j_1}^{-1}A_{j_2}^{-1}\rightarrow
 t_k A_{j_1}^{-1}A_{j_2}^{-1}A_{j_3}^{-1}\rightarrow\cdots
\]
such that $t_k A_{j_1}^{-1}A_{j_2}^{-1}\cdots A_{j_t}^{-1} \in V(\overline{DG}_t)\subset\mathcal{M}\quad (t=0,1,2,\cdots)$.
By $L(t_k)>L(t_kA_{j_1}^{-1})>L(t_kA_{j_1}^{-1}A_{j_2}^{-1})>\cdots$, monomials appearing in this sequence are all different.
Since we know $\mathcal{M}$ is a finite set, it is absurd. Therefore, the algorithm will stop at some step $r\in\mathbb{Z}_{\geq0}$ and
\[
V(\overline{DG})=V(\overline{DG}_r)\subset\mathcal{M}.
\]
Let us show $\mathcal{M}\subset V(\overline{DG})$. For any $M=\prod_{j=1}^N t_j^{d_j}\in \mathcal{M}$,
we prove $M\in V(\overline{DG})$. Our claim is evident when $M=t_k$ by
$M=t_k\in V(\overline{DG}_0)\subset V(\overline{DG})$, so we assume that $M\neq t_k$.
Considering Lemma \ref{lem2}, it follows $d_m<0$ with some $m\in[1,N]$.
By Lemma \ref{Amplus}, there exists $j_1\in[1,N]$ such that $M\cdot A_{j_1}\in\mathcal{M}$
and putting $M\cdot A_{j_1}=\prod_{l=1}t_l^{d_l'}$ and defining $(b_1',\cdots,b_N')$ as $N$-tuple of integers
of Definition \ref{bint} determined from $M\cdot A_{j_1}$, it holds
$d_{j_1}'>0$, $b_{j_1^+}'>0$ and either 
\begin{enumerate}
\item[$(a)$] $d_{j_1^+}'<d_{j_1}'$ or
\item[$(b)$] $d_{j_1^+}'=d_{j_1}'$ and there exists $p\in\mathbb{Z}_{\geq2}$ such that
$d_{j_1^{m+}}'=b_{j_1^{m+}}'=0$ $(m=2,3,\cdots,p-1)$ and $d_{j_1^{p+}}'=-1$, $b_{j_1^{p+}}'=1$,
\end{enumerate}
holds. Using Lemma \ref{lem-4}, one obtains $L(M\cdot A_{j_1})>L(M)$.
Let us prove
\begin{equation}\label{thm2-pr1}
\text{if }M\cdot A_{j_1}\in V(\overline{DG}_{l_1})\text{ with some }l_1\in\mathbb{Z}_{\geq0}\text{ then }M\in
V(\overline{DG}_{l_1+1}).
\end{equation}
In fact, if $M\cdot A_{j_1}\in V(\overline{DG}_{l_1})$ then one can take the smallest
$l_1'\in[0,l_1]$ such that $M\cdot A_{j_1}\in V(\overline{DG}_{l_1'})$.
By the definition of algorithm, $M\cdot A_{j_1}$ is a sink of $\overline{DG}_{l_1'}$ and
$M=M\cdot A_{j_1}\cdot A_{j_1}^{-1}\in V(\overline{DG}_{l_1'+1})\subset V(\overline{DG}_{l_1+1})$,
which implies (\ref{thm2-pr1}).
Repeating this argument, since $\mathcal{M}$ is finite, there exist $j_2,\cdots,j_s\in[1,N]$ such that
$M\cdot A_{j_1}\cdots A_{j_m}\in\mathcal{M}$,
\begin{equation}\label{thm2-pr2}
\text{if }M\cdot A_{j_1}\cdots A_{j_m}\in V(\overline{DG}_{l_m})\text{ with some }l_m\in\mathbb{Z}_{\geq0}\text{ then }M\cdot A_{j_1}\cdots
A_{j_{m-1}}\in
V(\overline{DG}_{l_m+1})
\end{equation}
for $m=2,3,\cdots,s$ and $M\cdot A_{j_1}\cdots A_{j_s}=t_k$.
Using the fact $t_k=M\cdot A_{j_1}\cdots A_{j_s}\in V(\overline{DG}_0)$,
(\ref{thm2-pr1}) and (\ref{thm2-pr2}), we see that
$M=t_k A_{j_s}^{-1}\cdots A_{j_1}^{-1}\in V(\overline{DG}_s)\subset V(\overline{DG})$. Therefore,
one gets $\mathcal{M}=V(\overline{DG})$.
\qed

\section{Type ${\rm G_2}$-case and minuscule cases}

\subsection{Minuscule case}

In this subsection, we take $i\in I$ such that
for any weight $\mu$ of $V(-w_0\Lambda_i)$ and $t\in I$, it holds
\[
\langle
h_t,\mu
\rangle\in\{1,0,-1\}
\]
and will prove the graph $\overline{DG}$ generated by the algorithm of Theorem \ref{thm2} coincides
with the graph $DG$ defined in Theorem 5.1 of \cite{KKN}.

Let $M=\prod_{l=1}^N t_l^{d_l}$ be a monomial in $\Delta_{w_0\Lambda_i,s_i\Lambda_i}\circ \theta^-_{\mathbf{i}}(t_1,\cdots,t_N)$
with a vector $(b_1,b_2,\cdots,b_N)$ in Definition \ref{bint}
and $\gamma=(\gamma_0,\cdots,\gamma_N)$ be the $\mathbf{i}$-trail corresponding to $M$.
First, we remark that the conditions $d_j>0$, $d_{j^+}=d_{j}$ ($j\in[1,N]$) imply $d_{j^+}=d_{j}=1$ so that
\[
\gamma_{j^+-1}(h_{i_j})=\gamma_{j^+}(h_{i_{j}})=1
\]
by Lemma 5.5 of \cite{KKN} 
and $b_{j^+}=c_{j^+}=\frac{\gamma_{j^+-1}(h_{i_j})-\gamma_{j^+}(h_{i_j})}{2}=0$. Therefore,
the condition (b) of (2) with $d_j>0$, $b_{j^+}>0$ in Theorem \ref{thm2} does not hold.
Hence, let us prove the condition (a) of (2) with $d_j>0$, $b_{j^+}>0$ is equivalent to
\begin{equation}\label{algo-prev}
d_j=1,\ d_{j^+}\neq1,
\end{equation}
which is the
necessary and sufficient condition
for the existence of
an arrow $M\rightarrow M\cdot A_j^{-1}$ in the graph of Theorem 5.1 of \cite{KKN}.
It is clear if the condition (a) of Theorem \ref{thm2} holds then (\ref{algo-prev}) holds
because $d_j\in\{-1,0,1\}$ (Lemma 5.5 of \cite{KKN}). Next, we suppose (\ref{algo-prev}).
By Lemma 5.7 of \cite{KKN}, it follows $d_{j^+}=0$. Using Proposition 5.8 in \cite{KKN},
we get
\[
b_{j^+}=c_{j^+}=1.
\]
Thus, the condition (a) of (2) in Theorem \ref{thm2} holds. \qed

\subsection{Type ${\rm G}_2$ case}

In this subsection, 
we take $G$ is of type ${\rm G}_2$ with Cartan matrix $(a_{i,j})_{i,j\in I}$ such that $a_{2,1}=-3$
and $a_{1,2}=-1$. 
Let us prove that the algorithm works for this case.

\vspace{3mm}

\nd
\underline{In the case $\mathbf{i}=(1,2,1,2,1,2)$}

\vspace{2mm}

By the results in \cite{KN},
\begin{eqnarray*}
& &\Delta_{w_0\Lm_1,s_1\Lm_1}\circ \theta^-_{\mathbf{i}}(t_1,\cdots,t_6)\\
&=&t_1+\frac{t_2^3}{t_3}+3\frac{t_2^2}{t_4}+3\frac{t_2t_3}{t_4^2}+\frac{t_3^2}{t_4^3}
+2\frac{t_3}{t_5}+\frac{t_4^3}{t_5^2}+3\frac{t_2t_4}{t_5}+3\frac{t_2}{t_6}
+3\frac{t_3}{t_4t_6}+3\frac{t_4^2}{t_5t_6}+3\frac{t_4}{t_6^2}+\frac{t_5}{t_6^3}.
\end{eqnarray*}
\[
\Delta_{w_0\Lm_2,s_2\Lm_2}\circ \theta^-_{\mathbf{i}}(t_1,\cdots,t_6)=t_6.
\]
Applying our algorithm to $\Delta_{w_0\Lm_1,s_1\Lm_1}$, we get the following graph $\overline{DG}$:
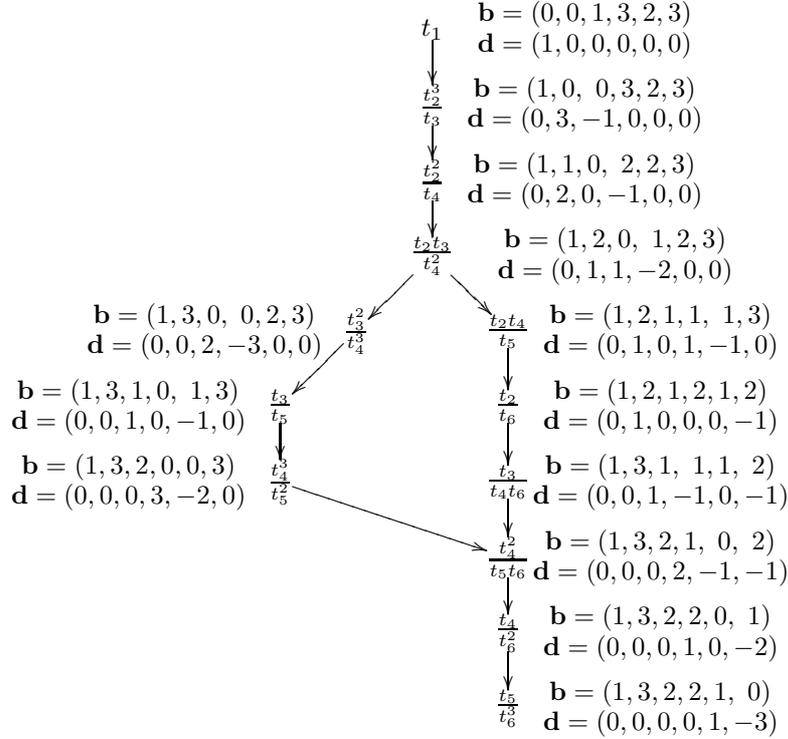
\[
\begin{xy}
(100,100) *{t_1}="1",
(120,102) *{\textbf{b}=(0,0,1,3,2,3)}="11",
(120,98) *{\textbf{d}=(1,0,0,0,0,0)}="111",
(100,90)*{\frac{t_2^3}{t_3}}="2",
(120,92) *{\textbf{b}=(1,0,\ 0,3,2,3)}="22",
(120,88) *{\textbf{d}=(0,3,-1,0,0,0)}="222",
(100,80)*{\frac{t_2^2}{t_4}}="3",
(120,82)*{\textbf{b}=(1,1,0,\ 2,2,3)}="33",
(120,78)*{\textbf{d}=(0,2,0,-1,0,0)}="333",
(100,70)*{\frac{t_2t_3}{t_4^2}}="4",
(124,72)*{\textbf{b}=(1,2,0,\ 1,2,3)}="44",
(124,68)*{\textbf{d}=(0,1,1,-2,0,0)}="444",
(90,60)*{\frac{t_3^2}{t_4^3}}="4-1",
(70,62)*{\textbf{b}=(1,3,0,\ 0,2,3)}="4-14-1",
(70,58)*{\textbf{d}=(0,0,2,-3,0,0)}="4-14-14-1",
(110,60)*{\frac{t_2t_4}{t_5}}="4-2",
(130,62)*{\textbf{b}=(1,2,1,1,\ 1,3)}="4-24-2",
(130,58)*{\textbf{d}=(0,1,0,1,-1,0)}="4-24-24-24",
(80,50)*{\frac{t_3}{t_5}}="5-1",
(60,52)*{\textbf{b}=(1,3,1,0,\ 1,3)}="5-15-1",
(60,48)*{\textbf{d}=(0,0,1,0,-1,0)}="5-15-15-1",
(110,50)*{\frac{t_2}{t_6}}="5-2",
(130,52)*{\textbf{b}=(1,2,1,2,1,2)}="5-25-2",
(130,48)*{\textbf{d}=(0,1,0,0,0,-1)}="5-25-25-2",
(80,40)*{\frac{t_4^3}{t_5^2}}="6-1",
(60,42)*{\textbf{b}=(1,3,2,0,0,3)}="6-16-1",
(60,38)*{\textbf{d}=(0,0,0,3,-2,0)}="6-16-16-1",
(110,40)*{\frac{t_3}{t_4t_6}}="6-2",
(130,42)*{\textbf{b}=(1,3,1,\ 1,1,\ 2)}="6-26-2",
(130,38)*{\textbf{d}=(0,0,1,-1,0,-1)}="6-26-26-2",
(110,30)*{\frac{t_4^2}{t_5t_6}}="7",
(130,32)*{\textbf{b}=(1,3,2,1,\ 0,\ 2)}="77",
(130,28)*{\textbf{d}=(0,0,0,2,-1,-1)}="777",
(110,20)*{\frac{t_4}{t_6^2}}="8",
(130,22)*{\textbf{b}=(1,3,2,2,0,\ 1)}="88",
(130,18)*{\textbf{d}=(0,0,0,1,0,-2)}="888",
(110,10)*{\frac{t_5}{t_6^3}}="9",
(130,12)*{\textbf{b}=(1,3,2,2,1,\ 0)}="99",
(130,8)*{\textbf{d}=(0,0,0,0,1,-3)}="999",
\ar@{->} "1";"2"^{}
\ar@{->} "2";"3"^{}
\ar@{->} "3";"4"^{}
\ar@{->} "4";"4-1"^{}
\ar@{->} "4";"4-2"^{}
\ar@{->} "4-1";"5-1"^{}
\ar@{->} "4-2";"5-2"^{}
\ar@{->} "5-1";"6-1"^{}
\ar@{->} "5-2";"6-2"^{}
\ar@{->} "6-1";"7"^{}
\ar@{->} "6-2";"7"^{}
\ar@{->} "7";"8"^{}
\ar@{->} "8";"9"^{}
\end{xy}
\]
Hence, our algorithm works.
The algorithm for $\Delta_{w_0\Lm_2,s_2\Lm_2}$ clearly works.
We also see that the graphs $\overline{DG}$ coincide with graphs $DG$ in subsection 6.1 of \cite{KKN}.

\vspace{3mm}

\nd
\underline{In the case $\mathbf{i}=(2,1,2,1,2,1)$}

\vspace{2mm}

As computed in \cite{KKN}, it holds
\[
\Delta_{w_0\Lm_2,s_2\Lm_2}\circ \theta^-_{\mathbf{i}}(t_1,\cdots,t_6)
=
t_1+\frac{t_2}{t_3}+\frac{t_3^2}{t_4}+2\frac{t_3}{t_5}+\frac{t_4}{t_5^2}+\frac{t_5}{t_6}.
\] 
Applying our algorithm to $\Delta_{w_0\Lm_2,s_2\Lm_2}$, one gets the following graph $\overline{DG}$:
\[
t_1\rightarrow
\frac{t_2}{t_3}\rightarrow
\frac{t_3^2}{t_4}\rightarrow
\frac{t_3}{t_5}\rightarrow
\frac{t_4}{t_5^2}\rightarrow
\frac{t_5}{t_6}.
\] 
Therefore, our algorithm works. Clearly, it holds
\[
\Delta_{w_0\Lm_1,s_1\Lm_1}\circ \theta^-_{\mathbf{i}}(t_1,\cdots,t_6)=t_6
\]
and our algorithm works. We also see that the graphs $\overline{DG}$ coincide with graphs $DG$ in subsection 6.1 of \cite{KKN}.


\end{document}